\documentclass{memo-l}
\newcommand{\R}{\mathbb{ R}}
\newcommand{\Rp}{\mathbb{ R}^{+}}
\newcommand{\Z}{\mathbb{ Z}}
\newcommand{\N}{\mathbb{ N}}
\newcommand{\C}{\mathbb{ C}}
\newcommand{\Half}{\mathbb{ H}}
\newcommand{\Q}{\mathbb{ Q}}

\newcommand{\by}{{\bf y}}
\newcommand{\uo}{{\bf u}_0}
\newcommand{\f}{{\bf f}}
\newcommand{\tB}{\tilde B}
\newcommand{\cH}{\mathcal{ H}}
\newcommand{\hf}{\mathcal{ H}^{f}}
\newcommand{\hs}{\mathcal{ H}^{s}}
\newcommand{\bq}{{\bf q}}
\newcommand{\bp}{{\bf p}}
\newcommand{\bx}{{\bf x}}
\newcommand{\p}{\psi}

\newcommand{\rk}{{\rm rk}}

\newcommand{\de}{ \delta }
\newcommand{\ep}{ \epsilon }

\newcommand{\Tmd}{ ( T^m)^\prime }
\newcommand{\Tnd}{ ( T^n )^\prime }

\newcommand{\Wt}{W_{z_0} (\psi) }

\newcommand{\Ja}{Jarn\'{\i}k }

\newcommand{\bu}{{\bf u}}
\newcommand{\QED}{{\hspace*{\fill}\raisebox{-1ex}{$ \#$}}}

\newtheorem{thkh}{Khintchine's Theorem (1924) \!\!\!\!\! }
\newtheorem{thjar}{Jarn\'{\i}k's Theorem (1931) \!\!\!\!\! }
\newtheorem{thgros}{ Theorem (Khintchine--Groshev) \!\!\!\!\!  }
\newtheorem{thbs}{Baker--Schmidt Theorem (1970) \!\!\!\!\! }
\newtheorem{thcs}{Theorem  \!\!\!\!\!\!  }
\newtheorem{thdel}{Theorem \ref{THM3}${}^{\prime}$ \!\!\!\!\!\!\!\! }

\newtheorem{thdet}{Theorem DV (1997) \!\!\!\!\!\!\! }
\newtheorem{problem}{Problem}
\newtheorem{lemma}{Lemma}
\newtheorem{theorem}{Theorem}
\newtheorem{corollary}{Corollary}
\newtheorem{proposition}{Proposition}
\newtheorem{adeq}{Lemma on adequate families of sets \!\!\!\!\!\!\! }

\newcommand{\B}{\mathcal{ B}}
\newcommand{\cR}{\mathcal{ R}}
\newcommand{\ka}{\kappa}
\newcommand{\ra}{R_{\alpha}}
\newcommand{\De}{\Delta}
\newcommand{\ma}{\beta_\alpha}
\newcommand{\dist}{\mathrm{dist}\,}
\renewcommand{\r}{\rho}
\newcommand{\La}{\Lambda}
\newcommand{\la}{\lambda}
\newcommand{\g}{\gamma}

\newcommand{\gb}{G_{\mbox{\tiny B}}}
\newcommand{\go}{G_{\mbox{\tiny $\Omega$}}}

\newcommand{\gbn}{G_{\mbox{\tiny $B_{n-1}$}}}

\newcommand{\ga}{\gamma}

\newcommand{\vbn}{V_{\mbox{\tiny $B_{n-1}$}}}

\newcommand{\ktk}{u_{t_1+i}}

\newcommand{\kd}{u_{t_2}}
\newcommand{\kdd}{u_{t_2+1}}

\newcommand{\be}{\begin{eqnarray*}}
\newcommand{\ee}{\end{eqnarray*}}
\newcommand{\kt}{u_{t_1}}
\newcommand{\ktt}{u_{t_1+1}}

\newcommand{\ve}{\varepsilon}
\renewcommand{\d}{\delta}

\newcommand{\fp}{f(\psi(}
\newcommand{\x}{\mathbf{x}}
\newcommand{\q}{\mathbf{q}}

\newcommand{\fN}{\mathcal{F}_N}
\newcommand{\br}{\mathbf{r}}

\begin{document}

\title[Measure theoretic laws for lim--sup sets]{Measure theoretic
  laws for lim--sup sets}

\author[V Beresnevich]{Victor Beresnevich\footnote{This work has been
 partially
supported by INTAS Project 00-429}
}
\address{V Beresnevich\\ Institute of Mathematics\\
Academy of Sciences of Belarus\\
220072\\ Surganova 11\\ Minsk\\ Belarus}
\email{beresnevich@im.bas-net.by}

\author[D Dickinson]{Detta Dickinson}
\address{D Dickinson\\ Department of Mathematics\\National
University of Ireland\\
 Maynooth\\  Co. Kildare\\  Republic of Ireland}
 \email{ddickinson@maths.may.ie}

\author[S Velani]{ Sanju Velani\footnote{Royal Society University Research
Fellow} }
\address{S Velani\\ Department of Mathematics\\ University
of York\\
 Heslington\\ York\\ YO10 5DD\\ England}
\email{slv3@york.ac.uk}

\maketitle
\newpage
\thispagestyle{empty}
\vspace*{2in}
\centerline{\Large{\textbf{Dedicated to Bridget, Ayesha and Iona}}}
\newpage
\thispagestyle{empty}
\vspace*{2cm}
 \section*{Abstract}
 Given a compact metric space $(\Omega,d)$ equipped with a
non-atomic, probability measure $m$ and a positive decreasing
function $\p$, we consider a natural class of $\limsup$ subsets
$\La(\psi)$ of $\Omega$. The classical $\limsup$ set $W(\psi)$ of
`$\p$--approximable' numbers in the theory of metric Diophantine
approximation fall within this class. We establish sufficient
conditions (which are also necessary under some natural
assumptions) for the $m$--measure of $\La(\psi)$ to be either
positive or full in $\Omega$ and for the Hausdorff $f$-measure to
be infinite. The classical theorems of  Khintchine-Groshev and
Jarn\'{\i}k concerning  $W(\psi)$  fall into our general
framework. The main results provide a unifying treatment of
numerous problems in metric Diophantine approximation including
those for real, complex and $p$-adic fields associated with both
independent and dependent quantities. Applications also include
those to  Kleinian groups and rational maps. Compared to previous
works our framework allows us to successfully remove many
unnecessary conditions and strengthen fundamental results such as
Jarn\'{\i}k's theorem and the Baker-Schmidt theorem. In
particular, the strengthening of Jarn\'{\i}k's theorem opens up
the Duffin-Schaeffer conjecture for Hausdorff measures.

\vskip 1cm
\textit{Mathematics Subject Classification: 11J83; 11J13, 11K60, 28A78, 28A80}

\newpage

\thispagestyle{empty}


\centerline{\Large{\textbf{Contents}}}

\vspace*{10mm}

\contentsline {section}{\tocsection {}{{\bf Section 1}}
\hspace{1ex} {{\bf Introduction}}}{1} {\parskip=1.4ex {\small
\contentsline {subsection}{\tocsubsection {}{\hspace{25mm}
\S1.1}{Background: the basic example}}{1} \contentsline
{subsection}{\tocsubsection {}{\hspace{25mm} \S1.2}{The general
setup and fundamental problems}}{5} }} \contentsline
{section}{\tocsection {}{{\bf Section 2}}\hspace{1ex}{{\bf
Ubiquity and conditions on the general setup}}}{8} {\parskip=1.4ex
{\small \contentsline {subsection}{\tocsubsection {}{\hspace{25mm}
\S2.1}{Upper and lower sequences and the sets $J^u_l(n)$}}{8}
\contentsline {subsection}{\tocsubsection {}{\hspace{25mm}
\S2.2}{The conditions on the measure and the space}}{9}
\contentsline {subsection}{\tocsubsection {}{\hspace{25mm}
\S2.3}{The intersection conditions}}{10} \contentsline
{subsection}{\tocsubsection {}{\hspace{25mm} \S2.4}{The ubiquitous
systems}}{10} \contentsline {subsection}{\tocsubsection
{}{\hspace{25mm} \S2.5}{A remark on related systems}}{13} }}
\contentsline {section}{\tocsection {}{{\bf Section
3}}\hspace{1ex}{{\bf The statements of the main theorems}}}{14}
\contentsline {section}{\tocsection {}{{\bf Section
4}}\hspace{1ex}{{\bf Remarks and corollaries to Theorem 1\hbox
{}}}}{16} \contentsline {section}{\tocsection {}{{\bf Section
5}}\hspace{1ex}{{\bf Remarks and corollaries to Theorem 2\hbox
{}}}}{18} \contentsline {section}{\tocsection {}{{\bf Section
6}}\hspace{1ex}{{\bf The classical results}}}{23} \contentsline
{section}{\tocsection {}{{\bf Section 7}}\hspace{1ex}{{\bf
Hausdorff measures and dimension }}}{24} \contentsline
{section}{\tocsection {}{{\bf Section 8}}\hspace{1ex}{{\bf
Positive and full $m$--measure sets}}}{26} \contentsline
{section}{\tocsection {}{{\bf Section 9}}\hspace{1ex}{{\bf Proof
of Theorem 1\hbox {}}}}{30} {\parskip=1.4ex {\small \contentsline
{subsection}{\tocsubsection {}{\hspace{25mm} \S9.1}{The subset
$A(\psi ,B)$ of $\Lambda (\psi )\cap B $}}{30} \contentsline
{subsection}{\tocsubsection {}{\hspace{25mm} \S9.2}{Proof of Lemma
8\hbox {} : quasi--independence on average \tmspace -\thinmuskip
{.1667em}\tmspace -\thinmuskip {.1667em} }}{34} }} \contentsline
{section}{\tocsection {}{{\bf Section 10}}{{\bf Proof of Theorem
2\hbox {}: $ 0\leq G < \infty $ }}}{37} {\parskip=1.4ex {\small
\contentsline {subsection}{\tocsubsection {}{\hspace{25mm}
\S10.1}{Preliminaries}}{37} \contentsline
{subsubsection}{\tocsubsubsection {}{10.1.1}{The sets $G_{\unhbox
\voidb@x \hbox {\chardef \@currsizeindex 2\relax \edef \begingroup
\@temptokena {{Beresnevich, Dickinson \& Velani}{\sectionrunhead
{}{10}{Proof of Theorem 2\hbox {}: $ 0\leq G < \infty $ }}}\xdef
{Beresnevich, Dickinson \& Velani}{\sectionrunhead {}{10}{Proof of
Theorem 2\hbox {}: $ 0\leq G < \infty $ }}{Beresnevich, Dickinson
\& Velani}\mark {Beresnevich, Dickinson \& Velani}\endgroup
{\relax \fontsize {10}{12}\selectfont }\begingroup \@temptokena
{{Beresnevich, Dickinson \& Velani}{\sectionrunhead {}{10}{Proof
of Theorem 2\hbox {}: $ 0\leq G < \infty $ }}}\xdef {Beresnevich,
Dickinson \& Velani}{\sectionrunhead {}{10}{Proof of Theorem
2\hbox {}: $ 0\leq G < \infty $ }}{Beresnevich, Dickinson \&
Velani}\mark {Beresnevich, Dickinson \& Velani}\endgroup B}}(n)$
and $G_{\unhbox \voidb@x \hbox {\chardef \@currsizeindex 2\relax
\edef \begingroup \@temptokena {{Beresnevich, Dickinson \&
Velani}{\sectionrunhead {}{10}{Proof of Theorem 2\hbox {}: $ 0\leq
G < \infty $ }}}\xdef {Beresnevich, Dickinson \&
Velani}{\sectionrunhead {}{10}{Proof of Theorem 2\hbox {}: $ 0\leq
G < \infty $ }}{Beresnevich, Dickinson \& Velani}\mark
{Beresnevich, Dickinson \& Velani}\endgroup {\relax \fontsize
{10}{12}\selectfont }\begingroup \@temptokena {{Beresnevich,
Dickinson \& Velani}{\sectionrunhead {}{10}{Proof of Theorem
2\hbox {}: $ 0\leq G < \infty $ }}}\xdef {Beresnevich, Dickinson
\& Velani}{\sectionrunhead {}{10}{Proof of Theorem 2\hbox {}: $
0\leq G < \infty $ }}{Beresnevich, Dickinson \& Velani}\mark
{Beresnevich, Dickinson \& Velani}\endgroup B}}(n,c)$.}}{37}
\contentsline {subsubsection}{\tocsubsubsection {}{10.1.2}{Working
on a subsequence of $u$.}}{39} \contentsline
{subsection}{\tocsubsection {}{\hspace{25mm} \S10.2}{The Cantor
set ${\bf K_{\eta }}$ }}{40} \contentsline
{subsubsection}{\tocsubsubsection {}{10.2.1}{Constructing the
first level ${\bf K(1)}$.}}{40} \contentsline
{subsubsection}{\tocsubsubsection {}{10.2.2}{Constructing the
second level {\bf K(2)}.}}{46} \contentsline
{subsubsection}{\tocsubsubsection {}{10.2.3}{Higher levels {\bf
K(n)} and the Cantor set ${\bf K_{\eta }}$.}}{49} \contentsline
{subsection}{\tocsubsection {}{\hspace{25mm} \S10.3}{A measure on
${\bf K_{\eta } }$}}{52} \contentsline
{subsubsection}{\tocsubsubsection {}{10.3.1}{Measure of a ball in
the Cantor construction.}}{53} \contentsline
{subsubsection}{\tocsubsubsection {}{10.3.2}{Measure of an
arbitrary ball $A$.}}{54} }}
 \contentsline {section}{\tocsection
{}{{\bf Section 11}}\hspace{1ex}{{\bf Proof of Theorem 2\hbox {}:
$G= \infty $ }}}{60} {\parskip=1.4ex {\small \contentsline
{subsection}{\tocsubsection {}{\hspace{25mm} \S11.1}{The Cantor
set {\bf K} and the measure $\mu $ }}{61} \contentsline
{subsection}{\tocsubsection {}{\hspace{25mm} \S11.2}{Completion of
the proof }}{62} }} \contentsline {section}{\tocsection {}{{\bf
Section 12}}\hspace{1ex}{{\bf Applications }}}{64} {\parskip=1.4ex
{\small \contentsline {subsection}{\tocsubsection {}{\hspace{25mm}
\S12.1}{Linear Forms }}{64} \contentsline
{subsection}{\tocsubsection {}{\hspace{25mm} \S12.2}{Algebraic
Numbers}}{66} \contentsline {subsection}{\tocsubsection
{}{\hspace{25mm} \S12.3}{Kleinian Groups }}{68} \contentsline
{subsection}{\tocsubsection {}{\hspace{25mm} \S12.4}{Rational
Maps}}{74} \contentsline {subsection}{\tocsubsection
{}{\hspace{25mm} \S12.5}{Diophantine approximation with
restrictions}}{78} \contentsline {subsection}{\tocsubsection
{}{\hspace{25mm} \S12.6}{Diophantine approximation in
$\mathbb{Q}_p$}}{79} \contentsline {subsection}{\tocsubsection
{}{\hspace{25mm} \S12.7}{Diophantine approximation on
manifolds}}{81} \contentsline {subsubsection}{\tocsubsubsection
{}{12.7.1}{Dual/linear approximation on manifolds.}}{82}
\contentsline {subsubsection}{\tocsubsubsection
{}{12.7.2}{Simultaneous approximation on manifolds.}}{84}
\contentsline {subsection}{\tocsubsection {}{\hspace{25mm}
\S12.8}{Sets of exact order}}{86} }} \contentsline
{chapter}{\tocchapter {Chapter}{}{{\bf Bibliography}}}{89}
\thispagestyle{empty}

\markboth{Beresnevich, Dickinson \& Velani}{\Section}

\mainmatter

\section{Introduction}

\subsection{Background: the basic example\label{beg}}
To set the scene for the  abstract framework considered in this
article we introduce a basic $\limsup$ set whose study has played
a central role in the development of the classical theory of
metric Diophantine approximation. Given a real, positive
decreasing function $\psi : \R^+ \to \R^+$, let $$ W (\p) := \{x
\in [0,1] :|x - p/q| < \p(q) \; {\rm for\ i.m.\ rationals\ }  p/q
\   (q>0) \}, $$ where `i.m.' means `infinitely many'. This is the
classical set of $\p$--well approximable numbers in the theory of
one dimensional Diophantine approximation. The fact that we have
restricted our attention to the unit interval rather than the real
line is purely for convenience.  It is natural to refer
 to the  function $\p$ as the {\em approximating function}.
It governs the `rate' at which points in the unit interval must be
approximated by rationals  in order to lie in $W(\p)$. It is not
difficult to see that $ W (\p) $ is a $\limsup$ set. For $n \in
\N$, let $$
 W(\p,n) :=  \!\!\!\!\!  \bigcup_{ { k^{n-1} < q  \leq k^n}} \bigcup_{ { 0 \leq  p  \leq q}}
\!\! B(p/q ,\p(q))  \cap [0,1]
$$ where  $k>1$ is
fixed and $B(c,r)$ is the open interval centred at $c$ of radius
$r$.  The set $ W(\psi)$ consists precisely  of points in the unit
interval that lie in infinitely many $ W(\p,n)$; that is
$$ W(\psi) = \limsup_{n \to \infty} W(\p,n) :=
\bigcap_{m=1}^{\infty} \bigcup_{n=m}^{\infty} W(\p,n) \ . $$

Investigating the measure theoretic properties of the set $W(\p)$
underpins the classical theory of metric Diophantine
approximation. We begin by considering the `size' of $W(\p)$
expressed in terms of the ambient measure $m$; i.e.
one--dimensional Lebesgue measure. On exploiting the $\limsup $
nature of $W(\p)$, a  straightforward application of the
convergence part of the Borel--Cantelli lemma from probability
theory yields that $$ m(W (\p)) = 0  \ \ \ \ \ {\rm if \ } \ \ \ \
\ \sum_{n=1}^{\infty} k^{2n} \p(k^n) \ < \ \infty \ . $$ Notice
that since $\psi$ is monotonic, the convergence/divergence
property of the above sum is equivalent to that of
$\sum_{r=1}^{\infty} r \, \p(r) $.

A natural problem now arises.  Under what conditions is  $m(W
(\p)) >0 $   ?  The  following fundamental result provides a
beautiful and simple criteria for the `size' of the set $W(\p)$
expressed in terms of Lebesgue measure.

\bigskip

\begin{thkh} If  $\p(r)$ is decreasing then
$$ m(W(\p)) =\left\{
\begin{array}{ll}
0 & {\rm if} \;\;\; \sum_{r=1}^{\infty} \; r \,  \p(r)  <\infty\;
,\\[1ex]
 1 & {\rm if} \;\;\; \sum_{r=1}^{\infty} \; r \,  \p(r)
 =\infty \; .
\end{array}\right.$$
\end{thkh}

Thus, in the divergence case,  which constitutes the main
substance of Khintchine's theorem,  not only do we have positive
Lebesgue measure but full Lebesgue measure. In fact, this  turns
out to be the case for all the examples considered in this paper.
Usually, there is a standard argument which allows one to deduce
full measure from positive measure -- such as the invariance of
the $\limsup$ set or some related set, under an ergodic
transformation. In any case, we shall prove a general result which
directly implies  the above full measure statement.  It is worth
mentioning that in Khintchine's original statement the stronger
hypothesis that $r^2 \p(r)$ is decreasing was assumed.

Returning to the convergence  case,  we cannot obtain any further
information regarding the `size' of $W(\p)$ in terms of Lebesgue
measure  --- it is always zero. Intuitively, the `size' of $W(\p)$
should decrease as the rate of approximation governed by the
function $\p$ increases. In short, we require a more delicate
notion of `size' than simply Lebesgue measure. The appropriate
notion of `size' best suited for describing the finer measure
theoretic structures of $W(\p)$ is that of generalized Hausdorff
measures. The Hausdorff $f$--measure $\mathcal{ H}^{f}$ with
respect to a dimension function $f$ is  a natural generalization
of Lebesgue measure.  So as not to interrupt the flow of this
background/motivation  exposition we referee the reader to
\S\ref{hmd} for the  standard definition of $\mathcal{ H}^{f}$ and
further comments regarding Hausdorff measures and dimension.

Again on exploiting the $\limsup $ nature of $W(\p)$, a
straightforward covering argument provides a simple convergence
condition under which $\mathcal{ H}^{f}(W(\p)) = 0 $. Thus, in
view of  the development of the Lebesgue theory it is  natural to
ask for conditions under which $\mathcal{ H}^{f}(W(\p))$  is
strictly positive.

The  following fundamental result provides a beautiful and simple
criteria for the `size' of the set $W(\p)$ expressed in terms of
Hausdorff measures.

\begin{thjar}  \label{main}
 Let $f$ be a dimension function such that
$r^{-1} \, f(r)\to \infty$ as $r\to 0 \, $ and  $r^{-1} \, f(r) $
is decreasing. Let $\p$ be  a real, positive  decreasing function.
Then $$ \hf\left(W(\p)\right)=\left\{\begin{array}{cl} 0 & {\rm \
if} \;\;\; \sum_{r=1}^{\infty}  \ \; r \,  f\left(\p(r)\right)
 <\infty \; ,\\[1ex]
\infty & {\rm \ if} \;\;\; \sum_{r=1}^{\infty} \  \; r \,
f\left(\p(r)\right) =\infty \; .
\end{array}\right.$$
\end{thjar}

Clearly  the above theorem  can be regarded as the Hausdorff
measure version of Khintchine's theorem.  As with the latter, the
divergence part constitutes the main substance. Notice, that the
case when $ \hf $ is comparable to one--dimensional Lebesgue
measure $m$ (i.e. $ f(r)= r$) is excluded by the condition $r^{-1}
\, f(r)\to \infty$ as $r\to 0 \, $. Analogous to Khintchine's
original statement, in  Jarn\'{\i}k's original statement the
additional hypotheses that $r^{2}\p(r)$ is decreasing, $r^{2}\p(r)
\to 0 $ as $ r \to \infty $ and that $r^2f(\p(r)) $ is decreasing
were assumed. Thus, even in  the simple case when $f(r) = r^s $ $
(s\geq 0) $  and the approximating function is given by $\psi(r) =
r^{-\tau} \log r $ $(\tau > 2)$,  Jarn\'{\i}k's original statement
gives no information regarding the $s$--dimensional Hausdorff
measure  of $W(\p)$ at the critical exponent $s=2/\tau$ -- see
below. That this is the case is due to the fact that  $r^2f(\p(r))
$ is not decreasing. However, as we shall see these additional
hypotheses are unnecessary. More to the point, Jarn\'{\i}k's
theorem as stated  above is the precise Hausdorff measure version
of Khintchine's theorem. Of course, as with Khintchine's theorem
the question of removing the monotonicity condition on the
approximating function $\p$ now arises.  That is to say, it now
makes perfect sense to consider a generalized Duffin-Schaeffer
conjecture for Hausdorff measures -- for a detailed account
regarding the original Duffin-Schaeffer conjecture see
\cite{Har,Sp}. Briefly, let $\psi(n)$ be a sequence of
non-negative real numbers and consider the set
$\widetilde{W}(\psi)$ of $x \in [0,1] $ for which there exist
infinitely many rationals $p/q$ $(q\geq 1)$ such that $$ |x - p/q|
\; < \; \psi(q) \hspace{9mm}  {\rm with } \hspace{9mm} (p,q)=1 \ .
$$

\noindent {\em The Duffin-Schaeffer conjecture for Hausdorff
measures}\/: Let $f$ be a dimension function such that $r^{-1} \,
f(r)\to \infty$ as $r\to 0 \, $ and $r^{-1} \, f(r) $ is
decreasing. Let $\phi$ denote the Euler function. Then $$
\hf\left(\widetilde{W}(\p)\right)=  \infty \ \ \ \ {\rm if} \ \ \
\  \sum_{n=1}^\infty \;  f\left(\p(n)\right) \, \phi(n) \;\;
=\infty \ . $$ It is easy to show that $\hf(\widetilde{W}(\p))= 0
$ if the above sum converges. The higher dimensional
Duffin-Schaeffer conjecture corresponding to simultaneous
approximation is known to be true \cite{PV}. It is plausible that
the ideas in \cite{PV} together with those in this paper are
sufficient to prove the higher dimensional version of the above
conjecture. The first and last authors have shown that this is
indeed the case \cite{DS}.

Returning to Jarn\'{\i}k's theorem,  note that in the case when $
\hf $ is the standard $s$--dimensional Hausdorff measure $\hs$
(i.e. $ f(r)= r^s $), it follows from the definition of Hausdorff
dimension (see \S\ref{hmd}) that $$ \dim W(\p) \, =  \, \inf \{ s
: \mbox{$\sum_{r=1}^\infty $} \; r \, \p(r)^s  < \infty \} \; . $$

Previously, Jarn\'{\i}k (1929) and independently Besicovitch
(1934) had determined the Hausdorff dimension of the set
$W(r\mapsto r^{-\tau})$, usually denoted by $W(\tau)$, of
$\tau$--well approximable numbers. They proved that for $\tau >
2$, $\dim W(\tau) = 2/\tau $. Thus, as the `rate' of approximation
increases (i.e. as $\tau$ increases) the `size' of the set
$W(\tau)$ expressed in terms of Hausdorff dimension decreases. As
discussed earlier, this is in precise keeping with one's
intuition. Obviously, the dimension result implies that $$
\mathcal{ H}^s \left(W(\tau)\right)=\left\{\begin{array}{ll} 0 & \
\ \ \
{\rm if} \;\;\; s \; > \; 2/\tau   \\
\infty & \ \ \ \ {\rm if} \;\;\; s \; < \; 2/\tau
\end{array}\right. ,$$
but gives no information regarding the $s$--dimensional Hausdorff
measure of $W(\tau)$ at the critical value $s=\dim W(\tau)$.
Clearly, Jarn\'{\i}k's zero--infinity law implies the dimension
result and that for $\tau > 2$ $$ \cH^{2/\tau}(W(\tau)) \, = \,
\infty \ . $$ Furthermore, the `zero--infinity' law allows us to
discriminate between sets with the same dimension and even the
same $s$--dimensional Hausdorff measure. For example, with $\tau
\geq 2 $ and   $ 0 < \ep_1 < \ep_2 $ consider the approximating
functions $$ \p_{\ep_i}(r) \, :=  \, r^{-\tau} \, \left( \log \, r
\right)^{ -  \frac{\tau}{2} \left(1 + \ep_i \right) } \hspace{1cm}
(i=1,2)  \; \; . $$ It is easily verified that for any $\ep_i >
0$, $$ m(W(\p_{\ep_i})) = 0 \, , \quad \dim W(\p_{\ep_i}) = 2/\tau
\quad {\rm and}\quad \cH^{2/\tau} (W(\p_{\ep_i})) =0 \ . $$
However, consider the dimension  function $f$ given by $ f(r) =
r^{2/\tau}( \log r^{-1/\tau})^{\ep_1}$.  Then $ \sum_{r=1}^\infty
r \, f\left(\p_{\ep_i} (r)\right) \asymp \sum_{r=1}^\infty \; (r
\, (\log r)^{1+ \ep_i - \ep_{1} } )^{-1} $, where as usual the
symbol $\asymp $ denotes comparability (the quotient of the
associated quantities is bounded from above and below by positive,
finite constants). Hence, Jarn\'{\i}k's zero--infinity law implies
that
\[
\hf\left(W(\p_{ \ep_1 } )\right) \ = \ \infty \hspace{12mm}  {\rm
whilst \ }  \hspace{12mm}    \hf\left(W(\p_{ \ep_2 } )\right) \ =
\ 0 \; .
\]
Thus the Hausdorff measure $ \hf $ does make a distinction between
the `sizes' of the sets  under consideration; unlike
$s$--dimensional Hausdorff measure.

Within this classical setup, it is apparent that Khintchine's
theorem together with Jarn\'{\i}k's zero--infinity law provide a
complete measure theoretic description of $W(\p)$. In short, our
central  aim is to establish  analogues of the divergence parts of
these classical results within a general framework. Recall, that
the divergence parts constitute the main substance of the
classical statements.

\subsection{The general setup and fundamental problems} Let
$(\Omega,d)$ be a compact metric space equipped with a non-atomic,
probability measure $m$. Let $\cR=\{\ra \subset \Omega : \alpha
\in J \}$  be a family of subsets $\ra$ of $\Omega$ indexed by an
infinite, countable set $J$. The sets $\ra$ will be referred to as
{\em resonant sets} for reasons which will become apparent later.
Next, let $\beta: J \to \R^+ : \alpha \mapsto \ma $ be a positive
function on $J$. Thus, the function $\beta$ attaches a `weight'
$\ma$  to the  resonant set $\ra$. To avoid pathological
situations within our framework, we shall assume that the number
of $\alpha$ in $J$ with $\ma$ bounded above is always finite. For
a set $A \subset \Omega$, let $$\De(A,\de) := \{ x \in \Omega :
\dist(x,A) < \de \} $$ where $\dist(x,A) := \inf \{d(x,a) : a \in
A \} $. Thus, $\De(A,\de)$ is simply the $\de$-neighborhood of
$A$. Given a decreasing function $\psi : \R^+ \to \R^+ $ let
$$\La(\psi)=\{x\in\Omega:x\in\De(\ra,\p(\ma))\ \mbox{for\
infinitely\ many\ }\alpha\in J\} \ . $$ The set $\La(\psi)$ is a
`$\limsup$' set; it consists of points in $\Omega$ which lie in
infinitely many of the `thickenings' $\De(\ra,\p(\ma))$. Clearly,
even in this abstract setup  it is natural to refer to the
function $\p$ as the {\em approximating function}. It governs the
`rate' at which points in $\Omega$ must be approximated by
resonant sets in order to lie in $\La(\p)$. Notice, that in the
case the resonant sets are points, the thickenings
$\De(\ra,\p(\ma))$ are simply balls $B(\ra,\p(\ma))$ centred at
resonant points.

Before continuing our discussion, we rewrite $\La(\p)$ in a
fashion which brings  its `$\limsup$' nature to the forefront. For
$n \in \N$, let $$
 \De(\p,n) :=  \!\!\!\!\!  \bigcup_{\alpha\in J \ : \ { k^{n-1} < \ma \leq k^n}}
\!\!\!\!\!\!\!\!\!\! \De(\ra,\p(\ma)) \ \ \  \mbox{\rm where $k>1$
is fixed.} $$ By assumption the number of $ \alpha$ in $J$ with $
k^{n-1} < \ma \leq k^n $ is finite regardless of the value of $k$.
Thus, $\La(\psi)$ is precisely the set of points in $\Omega$ which
lie in infinitely many $\De(\p,n)$; that is $$ \La(\psi) =
\limsup_{n \to \infty} \De(\p,n) := \bigcap_{m=1}^{\infty}
\bigcup_{n=m}^{\infty} \De(\p,n) \ . $$

The main line of our investigation is motivated by the following
pair of fundamental  problems regarding the measure theoretic
structure of $\La(\p)$. In turn the fundamental problems are
motivated by the classical theory described in the  previous
section. It is reasonably straightforward to determine conditions
under which $m(\La(\p)) = 0$. In fact, this is implied by the
convergence part of the Borel--Cantelli lemma from probability
theory whenever
\begin{equation}\label{*1}
\sum_{n=1}^\infty m(\De(\p,n)) < \infty  \ .
\end{equation}
In view of this it is natural to consider:

\begin{problem}
Under what conditions is $m(\La(\p))$ strictly positive ?
\end{problem}

Under a `global ubiquity' hypothesis and a divergent sum
condition, together with  mild conditions on the measure, our
first theorem provides a complete solution to this problem.
Moreover, if we replace the `global ubiquity' hypothesis by a
`local ubiquity' hypothesis then $\La(\p)$ has full $m$--measure
and this statement can be viewed as the analogue of Khintchine's
theorem or more generally as the analogue of the classical linear
forms theorem of Khintchine--Groshev.

Reiterating the above measure zero statement, if the approximating
function $\p$  decreases sufficiently quickly so that (\ref{*1})
is satisfied, the corresponding $\limsup $ set $\La(\p)$ is of
zero $m$--measure.  As with the classical setup of \S\ref{beg}, in
this case we cannot obtain any further information regarding the
`size' of $\La(\p)$ in terms of the ambient measure $m$  --- it is
always zero.  In short, we require a more delicate notion of
`size' than simply the given $m$-measure. In keeping with the
classical development, we investigate the `size' of $\La(\p)$ with
respect to  the Hausdorff measures $ \mathcal{ H}^{f} $ where $f$
is a dimension function.  Again, provided a certain `$f$-volume'
sum converges, it is reasonably simple to determine conditions
under which $\mathcal{ H}^{f}(\La(\p)) = 0$. Naturally, we
consider:

\begin{problem}
Under what conditions  is $\mathcal{ H}^{f}(\La(\p))$ strictly
positive ?
\end{problem}

 This problem turns out to be far more subtle than
the previous one regarding $m$-measure. To make any substantial
progress, we impose the condition that the $m$-measure of any ball
centred at a point in $\Omega$ is  comparable to some fixed  power
of its radius. Then, under a `local ubiquity' hypothesis and an
`$f$-volume' divergent sum condition, together with mild
conditions on the dimension function, our second theorem shows
that $\mathcal{ H}^{f}(\La(\p)) = \infty$. Thus, $\mathcal{
H}^{f}(\La(\p))$ satisfies an  elegant `zero--infinity' law
whenever the convergence of the `$f$-volume' sum implies
$\mathcal{ H}^{f}(\La(\p)) = 0$ as is often the case. In
particular, this latter statement is true for the standard
$s$-dimensional Hausdorff measure $\mathcal{ H}^{s}$. Thus, in the
language of geometric measure theory the sets $\La(\p)$ are not
$s$-sets. Furthermore, from such zero--infinity laws it is easy to
deduce the Hausdorff dimension of $\La(\p)$.

Examples of $\limsup$ sets which fall into the above abstract
framework include the classical sets of well approximable
numbers/vectors in the  theory of Diophantine approximation as
well as the `shrinking target' sets associated with the phase
space of a given dynamical system.

In order to illustrate and clarify the above setup and our line of
investigation, we return to the  basic $\limsup$  set  of
\S\ref{beg}. The classical set $W(\p)$  of $\p$--well approximable
numbers in the theory of one dimensional Diophantine approximation
can clearly  be expressed in the form $\La(\p)$ with
\begin{eqnarray*}
& &\Omega:= [0,1] \  , \ \  J:= \{ (p,q) \in \N \times \N : 0\le p
\leq q \} \ , \ \ \alpha := (p,q) \in J \ , \
\\[1.5ex]
& &\ma := q \ , \ \ \ra:= p/q  \ \  \ {\rm and \ } \ \
\De(\ra,\p(\ma)) := B(p/q, \p(q)) \ .
\end{eqnarray*}
The metric $d$ is of course the standard Euclidean metric; $d(x,y)
:= |x-y| \, $.  Thus in this basic example, the resonant sets
$\ra$ are simply rational points $p/q$. Furthermore,
$$\De(\p,n) := 
\bigcup_{k^{n-1} < q \leq k^n} \ \bigcup_{p=0}^{q} B(p/q, \p(q)) \
$$ and   $ W (\p) = \limsup \De(\p,n)$ as ${n \to \infty}$.

For this basic example, the solution to our first fundamental
problem is given by   Khintchine's theorem and the solution to the
second  by Jarn\'{\i}k's theorem. Together, these theorems
provide a complete measure theoretic description of $W(\p)$. In
the case of the general framework, analogues of these results
should be regarded as the ultimate  pair of results describing the
metric structure of the $\limsup$ sets $\La(\p)$. Alternatively,
they provide extremely satisfactory solutions  to the fundamental
problems. Analogues  of the convergence parts of the classical
results usually follow by adapting the `natural cover'
$$\{\De(\p,n) \, : \,  n=m, m+1, \cdots  \}   \hspace{16mm} (m \in
\N) \ $$ of $\La(\p)$. Our key  aim is to establish analogues of
the divergence parts of the classical results for $\La(\p)$.

\section{Ubiquity and conditions on the general setup\label{congs}}

In order to make any reasonable progress with the fundamental
problems we impose various conditions on the metric measure space
$(\Omega,d,m)$. Moreover, we require the notion of a `global' and
`local' ubiquitous system which will underpin our line of
investigation. The general setup is independent of  the
approximating function $\p$.

Throughout, a ball centred at a point $x$ and radius $r$ is
defined to be the set $\{ y \in \Omega : d(x,y) <  r \}$ or $\{ y
\in \Omega : d(x,y) \leq  r \}$ depending on whether it is open or
closed. In general,  we do not specify whether a certain ball is
open or close since it will be  irrelevant. Notice, that by
definition any ball is automatically a subset of $\Omega$.

\subsection{Upper and lower sequences and the sets $J^u_l(n)$}
Let $l:=\{l_n\}$ and $u:=\{u_n\}$ be positive increasing sequences
such that $$ l_n < u_n \ \hspace{14mm} {\rm and } \hspace{14mm}
\lim_{n \to \infty} l_n = \infty   \ . $$ Thus,  $ \lim_{n \to
\infty} u_n = \infty $.  Now, define $$
 \De^u_l(\p,n) := \ \bigcup_{\alpha\in J^u_l(n)}
\De(\ra,\p(\ma)) $$ where $$ J^u_l(n):=\{\alpha\in J \, : \,
l_n<\ma\le u_n\}. $$

\noindent By assumption the cardinality of $J^u_l(n)$ is finite
regardless of $l$ and $u$.  In view of this and the fact that $
l_n \to \infty$ as $n \to \infty$, it follows that   $$ \La(\psi)
= \limsup_{n \to \infty} \De^u_l(\p,n) := \bigcap_{m=1}^{\infty}
\bigcup_{n=m}^{\infty} \De^u_l(\p,n) \ . $$ This statement is
irrespective of the choice of the sequences $l$ and $u$. Note that
without additional assumptions, the fact that $ l_n \to  \infty$
as $n \to \infty$ is crucial. For obvious reasons, the sequence
$l$ will be referred to as {\it the lower sequence}\/ and $u$ as
{\it the upper sequence}\/.

\subsection{The conditions on the measure  and the space}
The two central conditions on the measure $m$ are as follows and
will always be assumed throughout. Firstly, the $m$--measure of
any ball centred at a point of the space $\Omega$ is strictly
positive; i.e. $ m(B(x,r))  >  0 $ for $x$ in $\Omega$ and $r>0$.
Secondly, the measure  $m$ is {\em doubling}.  That is to say that
there exists a constant $C \geq 1 $ such that for any $x$ in
$\Omega$ $$ m(B(x,2r)) \, \leq \, C \ m(B(x,r))  \ . $$ The
doubling condition allows us to blow up a given ball by a constant
factor without drastically affecting its measure. Also note that
it implies that  $ m(B(x, t r)) \, \leq C(t) \, m(B(x,r)) $ for
any $t > 1 $ with $x$ in $\Omega$. The metric measure space
$(\Omega,d,m)$ is also said to be  doubling if $m$ is doubling
\cite{jh}.

Regarding `Problem 1', we shall impose the following reasonably
mild conditions on the measure. Essentially, it asserts that balls
of the same radius centred at points on resonant sets $\ra $ with
$\alpha \in J^u_l(n):=\{\alpha\in J:l_n<\ma\le u_n\}$ have roughly
the same measure for some choice of $l$ and $u$.

\medskip

\noindent {\bf (M1) \ } For $c\in\ra$, $c'\in R_{\alpha'}$ with
$\alpha,\alpha'\in J^u_l(n)$ and $r < r_o$
\begin{eqnarray}
a \ \leq \  \frac{ m(B(c,r)) }{   m(B(c',r))} \ \leq \  b     \ ,
\label{measurem1}
\end{eqnarray}
where the constants $a, b >0$ are  independent of $n$ and  the
balls under consideration, but may depend on the $l$ and $u$.

\medskip

Regarding `Problem 2', more is required. Namely, that the measure
of a ball centred at a point in $\Omega$ is comparable to some
fixed  power  of its radius.

\medskip

\noindent {\bf (M2) \ } There exist positive constants $\delta$
and $r_o$   such that for any $ x \in \Omega $ and $ r \leq r_o $,
\begin{equation}
 a \, r ^{\delta}  \ \leq  \  m(B(x,r))  \ \leq  \   b \, r
^{\delta}    \ . \label{measure}
\end{equation}

\medskip

The  constants $a$ and $b$ are independent of the ball and without
loss of generality we assume that $0<a<1<b$. Notice that if $m$
satisfies condition (M2), then (M1) is trivially satisfied for any
choice of $l$ and $u$, as are the central conditions.  Also, (M2)
implies that $\dim \Omega = \delta $ -- see \S\ref{hmd} for the
details.

\subsection{The intersection conditions}
In the case that the resonant sets are not points,  we will
require measure theoretic control on the intersection of certain
balls centred at points on resonant sets with certain
`thickenings' of the resonant sets. The radii of the balls and the
thickenings of the resonant sets are governed by a positive
function $\rho$ which is intimately tied up with the notion of
ubiquity -- see below.

\bigskip

\newtheorem{intcond}{The intersection conditions \!\!\!\!\!\!\! }
\renewcommand{\theintcond}{}

\begin{intcond}\rm
There exists a constant $\g$ with $ 0 \leq \g \leq \dim \Omega $,
such that for any  $\alpha\in J$ with $\beta_{\alpha} \leq u_n$,
$c\in\ra$ and   $0< \la \le \r(u_n)$ the following are satisfied
for $n$ sufficiently large:

\begin{enumerate}

\item[] ~~

\item[$(i)$]
~~~$\displaystyle m(B(c, {\mbox{\small
$\frac{1}{2}$}}\r(u_n))\cap\De(\ra,\la))\geq c_{\mbox{\tiny 1}} \,
m(B(c,\la))\left(\frac{\r(u_n)}{\la}\right)^{\g} \  $

\item[] ~~

\item[] ~~

\item[$(ii)$]
~~~$\displaystyle m(B\cap B(c,3\r(u_{n}))\cap\De(\ra,3\la))\leq
c_{\mbox{\tiny 2}} \,
m(B(c,\la))\left(\frac{r(B)}{\la}\right)^{\g} \  $
\end{enumerate}

\noindent  where $B$ is an arbitrary ball centred on a resonant
set with radius $r(B)\le 3 \, \r(u_{n})$. The   constants
$c_{\mbox{\tiny 1}}$ and $ c_{\mbox{\tiny 2}}$  are positive and
absolute. Without loss of generality we assume that
$0<c_{\mbox{\tiny 1}}<1<c_{\mbox{\tiny 2}}$.

\end{intcond}

\medskip

When the resonant sets are points so that $\De(\ra,\lambda):=
B(c,\lambda)$, the above conditions are trivially satisfied with
$\g=0$. In  applications, when the resonant sets are not points
they are usually subsets of smooth manifolds or simply planes, all
of the same dimension. In such cases, the intersection conditions
hold with $\g = \dim \ra $. In particular, it is readily verified
that when $(\Omega,d)$ is a subspace of $\R^n$ and the resonant
sets are $\gamma$--dimensional affine subspaces of $\Omega$ then
the intersection conditions are inevitably satisfied. In view of
this we refer to $\g$ as the {\em common dimension} of the
resonant sets in $\cR$.

\subsection{The ubiquitous systems\label{thesystems}}
The following `systems' contain the key measure theoretic
structure necessary for our attack on the fundamental problems.
Recall that $\cR$ denotes the family of resonant sets $\ra$ and
that the function $\beta$ attaches a `weight' $\beta_{\alpha}$ to
each resonant set $\ra \in \cR$.



Let $\r : \R^+ \to \R^+ $ be a  function with $\r(r) \to 0 $ as $r
\to \infty $ and let

$$ \De^u_l(\r,n) := \bigcup_{\alpha\in J^u_l(n)} \De(\ra,\r(u_n))
\ .$$

\medskip

\noindent{\bf Definition (Local $\, m$--ubiquity)}~ Let $B=B(x,r)$
be an arbitrary ball with centre $x$ in $\Omega$ and radius $r
\leq r_0$. Suppose there exists a function $\r$, sequences $l$ and
$u$ and an absolute constant $\ka>0$ such that
\begin{equation}\label{ub}
  m\left( B \cap \De^u_l(\r,n) \right) \ \ge \ \ka \ m(B) \qquad \mbox{for $n \geq n_o(B)$} .
\end{equation}
Furthermore, suppose the intersection conditions are satisfied.
Then the pair $(\cR,\beta)$ is said to be a {\em local
$m$-ubiquitous system relative to $(\r,l,u)$}.

\medskip

\noindent{\bf Definition (Global $\, m$--ubiquity)}~ Suppose there
exists a function $\r$, sequences $l$ and $u$  and an  absolute
constant $\ka>0$ such that for $n \geq n_o$, (\ref{ub}) is
satisfied for $B:=\Omega$. Furthermore, suppose the intersection
conditions are satisfied. Then the  pair $(\cR,\beta)$ is said to
be a {\em global $m$-ubiquitous system relative to $(\r,l,u)$}.

\medskip

The function $\r$, in either form of ubiquity will be referred to
as the {\em ubiquitous function}. It is  clear that for global
ubiquity (\ref{ub}) reduces to $ m\left( \De^u_l(\r,n) \right) \
\ge \ \ka$. Essentially, in the `global' case all that is required
is that the set  $\De^u_l(\r,n)$ approximates the underlying space
$\Omega$ in terms of the measure $m$. In the `local' case, this
approximating property is required to hold locally on balls
centred at points in $\Omega$.

 Clearly, local $m$--ubiquity implies global
$m$--ubiquity. Simply take  a  ball $B$ centred at a point of
$\Omega$ with radius $\leq r_o$, then for $n$ sufficiently large $
m\left( \De^u_l(\r,n) \right)  \geq  m\left( B \cap \De^u_l(\r,n)
\right)  \ge  \ka \, m(B)   :=  \ka_1  >  0  $. In other words
local ubiquity with a constant $\ka$ implies global ubiquity with
some constant $\ka_1$ where $0<\ka_1 \leq \ka$. In general the
converse is not true. However, if $$ m\left( \De^u_l(\r,n) \right)
\ \to \  1 \,  = \, m(\Omega) \hspace{8mm} {\rm as } \hspace{8mm}
n  \, \to \,  \infty  \ , $$ then it is easy to  show that global
$m$--ubiquity implies local $m$--ubiquity. To see this, let $B$ be
any ball and assume without loss of generality that $m(B) = \ep >
0 $. Then, for $n$ sufficiently large  $ m\left( \De^u_l(\r,n)
\right)
> m(\Omega) - \ep/2$. Hence,  $ m\left( B \cap \De^u_l(\r,n)
\right) \geq \ep/2$ as required.  This rather simple observation
can be extremely useful when trying to establish that a given
system is locally $m$-ubiquitous.

In attempting to establish the measure theoretic inequality in
either form of ubiquity, the presence of the lower sequence $l$ is
irrelevant. To see this, suppose we are able to show that for $n
\geq n_o(B)$, $$ \textstyle m\left( B \cap  \bigcup_{\alpha\in J:
\ma \leq u_n} \De(\ra,\r(u_n)) \right) \ \ge \ \ka \ m(B) \ . $$
Since $\r(r) \to 0 $ as $r \to \infty$, for any $t \in \N$ there
exists an integer $n_t$ such that for $n \geq n_t $ $$ \textstyle
m\left( B \cap \bigcup_{\alpha\in J: \ma \leq t}\ \De(\ra,\r(u_n))
\right) \ < \  \frac{1}{2} \,  \ka \ m(B) \  . $$ Without loss of
generality we can assume that $n_{t+1} \geq n_t +1 $. Now consider
the lower sequence $l$ given by $l_n:=t$ for  $n \in [n_t,
n_{t+1})$. Clearly, $l$ is an increasing sequence with $l_n \to
\infty $ as $n \to \infty$. Moreover, for $n$ sufficiently large
we have that $
 m\left( B \cap \De^u_l(\r,n) \right)  \, \ge \,
 \mbox{\footnotesize{$\frac{1}{2}$}} \,  \ka \, m(B)
$; i.e. the pair $(\cR,\beta)$ is  a local $m$-ubiquitous system
relative to $(\r,l,u)$.

The above discussion indicates that  the lower sequence $l$ is
irrelevant within the ubiquity framework.    Regarding the upper
sequence, notice that any subsequence $s$ of $u$  will also do;
i.e. the measure theoretic inequalities are satisfied for
$\De^s_l(\r,n)$.  To see that this is the case, simply observe
that for each $t \in \N$ we have that  $s_t= u_n$ for some $n \geq
t $. Then, $l_t \leq l_n$ and so $J^s_l(t) \supseteq J^u_l(n) $.
Thus $\De^s_l(\r,t) \supseteq \De^u_l(\r,n)$, and it follows that
$m\left( B \cap \De^s_l(\r,t) \right)  \, \ge \,  m\left( B \cap
\De^u_l(\r,n) \right) \, \ge \,   \ka \, m(B)  $; i.e. the pair
$(\cR,\beta)$ is  a local $m$-ubiquitous system relative to
$(\r,l,s)$.



 In practice, the global or local $m$--ubiquity of a
system can be established using standard arguments concerning the
distribution of the resonant sets in $\Omega$, from which the
function $\r$ arises naturally. To  illustrate this, we return to
our basic example.

\noindent  {\it The basic example again: \ }For the set $W(\p)$ of
$\p$--well approximable numbers the resonant sets are simply
rational points. Thus the intersection conditions are
automatically satisfied with $\g=0$. Of course, the measure  $m$
is one--dimensional Lebesgue measure and satisfies the measure
condition (M2) with $\delta =1 $.

\bigskip

\begin{lemma}\label{reg_be}
There is a constant $k>1$ such that the pair $(\cR,\beta)$ is a
local $m$-ubiquitous system relative to $(\r,l,u)$ where $l_{n+1}
\; = \; u_n \; := \;  k^n$ and $\r : r \mapsto  {\rm constant}
\times r^{-2}$.
\end{lemma}
\noindent{\em Proof.}~ \  Let $I=[a,b]\subset[0,1]$. By
Dirichlet's theorem, for any $x\in I$ there are coprime integers
$p,q$ with $1\le q\le k^n$ satisfying $|x-p/q|<(qk^n)^{-1}$.
Clearly, $aq-1\le p\le bq+1$. Thus, for a fixed $q$ there are at
most $m(I)q+3$ possible values of $p$. Trivially, for $n$ large
 $$
m\left( I \cap  \!\! \bigcup_{q \le k^{n-1} }    \bigcup_{p} \
\textstyle{B\left(\frac{p}{q}, \frac{1}{q k^n} \right) } \right) \
\le \ 2 \!\! \sum_{q \, \leq k^{n-1} } \!\!\! \mbox{\large
$\frac{1}{qk^{n}} $ } (m(I)q+3) \le \mbox{\large $\frac{3}{k} $ }
 m(I)  \, .$$
 \noindent It follows that for  $k\geq 6$,
$$m\left( I \cap  \!\! \bigcup_{k^{n-1} < q \le k^n } \bigcup_{p}
\ \textstyle{B\left(\frac{p}{q}, \frac{k}{k^{2n}} \right) }
\right) \ \geq \ m(I) - \mbox{\large $\frac{3}{k} $ }m(I)  \ \geq
\ \mbox{\large $\frac{1}{2} $ } m(I) \ . \vspace{-6mm} $$
 \QED

\medskip

It will be evident from our `ubiquity' theorems, that
Lemma~\ref{reg_be} is sufficient  for directly  establishing the
divergence part of both Khintchine's theorem and Jarn\'{\i}k's
zero--infinity law -- see \S\ref{basiceg}.

\subsection{A remark on related systems}
In the case that $\Omega$ is a bounded subset of $\R^n$ and $m$ is
$n$-dimensional Lebesgue measure, the notion of ubiquity was
originally formulated by Dodson, Rynne \& Vickers \cite{DRV}  to
obtain lower bounds for the Hausdorff dimension of the sets
$\Lambda(\p)$ -- see \S\ref{r3}. Their notion of ubiquity is
closely related to our notion of a `local $m$-ubiquitous' system.
In the case that the resonant sets are points the ubiquitous
systems of Dodson, Rynne \& Vickers  coincide with the `regular
systems' of Baker \& Schmidt \cite{BS}. Both these systems have
proved very useful in obtaining lower bounds for the Hausdorff
dimension of $\limsup$ sets. However, both \cite{BS} and
\cite{DRV} fail to shed any light on the fundamental problems
considered in this paper. For further details regarding  regular
systems and the original formulation of ubiquitous systems see
\cite{BD}.

Recently and independently, in \cite{yann2} the notion of an
optimal regular system introduced in \cite{Ber2000} has been
re-formulated to obtain divergent type Hausdorff measures results
for subsets of $\R^n$. This re-formulated notion is essentially
equivalent to our notion of local $m$-ubiquity in which  $m$ is
$n$-dimensional Lebesgue measure, the resonant sets are points
($\gamma=0$), the ubiquity function is comparable to $\rho: r \to
r^{\mbox{\tiny $-1/n$ }} $ and the sequences $l$ and $u$ are given
by $l_{n+1}  =  u_n := \; 2^n$. These highly restrictive
conditions, in particular the latter two which fix the function
$\rho$ and the sequences $l$ and $u$, excludes many of the
applications we have in mind even when $\Omega$ is a subset of
$\R^n$. Furthermore, even with the restrictions  our notion of
local $m$-ubiquity is not equivalent to that of an optimal regular
system since we make no assumption on the growth of $\# J^u_l(n)$.

\section{The statements of the main theorems}
First some useful notation.  Let $m$ be a  measure  satisfying
condition (M1) with respect to the  sequences $l$ and $u$.  Then
$B_n(r)$ will  denote a generic ball of radius $r$ centred at a
point of a resonant set $\ra$ with $\alpha$ in $J^u_l(n)$. Given
the conditions imposed on the measure $m$, we have that for any
ball $B(c,r)$ with $c \in \ra$ and $\alpha \in J^u_l(n)$ $$
m(B(c,r)) \ \asymp \ m(B_n(r)) \ . $$ This comparability is
obviously satisfied for any $c$ in $\Omega$ if the measure
satisfies (M2). With this in mind, we now state our main results.
Recall, that an approximating function $\psi$ is a real, positive
decreasing function and that a ubiquity function $\rho$ is a real,
positive function such that $\r(r) \to 0 $ as $r \to \infty $.

\bigskip

\begin{theorem} \label{THM1}
Let $(\Omega,d)$ be a compact metric space equipped with a measure
$m$ satisfying condition\/ {\rm(M1)}\/ with respect to sequences
$l$ and $u$. Suppose that $(\cR,\beta)$ is a global
$\,m$-ubiquitous system relative to $(\r,l,u)$ and that $\, \p$ is
an approximating function. Assume that
\begin{equation}
\limsup_{n \to \infty} \frac{\psi(u_n)}{\r(u_n)} \ > \ 0
\label{af1}
\end{equation}
or assume that
\begin{equation}\label{af2}
\sum_{n=1}^{\infty} \ \frac{m(B_n(\p(u_n)))}{m(B_n(\r(u_n)))} \ \
\left(\frac{\r(u_n)}{\p(u_n)}\right)^{\g} \ = \ \infty
\end{equation}
and  for $Q$ sufficiently large
\begin{equation}\label{af}
\begin{array}[b]{rl}
\displaystyle \sum_{s=1}^{Q-1} \
\frac{\r(u_s)^{\g}}{m(B_s(\r(u_s)))} &\displaystyle
\sum_{\substack{s+1 \, \leq \, t \leq \,  Q \, :\\ \p(u_s) \, < \,
\r(u_t)}}
 \frac{m(B_t(\p(u_t)))}{\p(u_t)^{\g}} \\[5ex]
 &\displaystyle \ll   \left( \sum_{n=1}^{Q} \ \frac{m(B_n(\p(u_n)))}{m(B_n(\r(u_n)))} \
 \left(\frac{\r(u_n)}{\p(u_n)}\right)^{\g} \right)^2  .
\end{array}
\end{equation}
 Then,
 $ m \left( \Lambda(\p) \right)  >  0 $. In addition, if any open subset
  of $\Omega$ is $m$--measurable and  $(\cR,\beta)$ is locally $\,m$-ubiquitous
relative to $(\r,l,u)$, then $ m \left( \Lambda(\p) \right) = 1 $.
\end{theorem}

\medskip

 Before stating the Hausdorff measure analogue of the above
theorems we introduce the following notion. Given a sequence $u$,
a function $h$ will be said to be {\bf ${\bf u}$-regular} if there
exists a strictly positive constant $\lambda < 1$ such that for
$n$ sufficiently large
\begin{equation} h(u_{n+1}) \leq \lambda \, h(u_n) \ .  \label{afmh}
\end{equation}
The constant $\lambda$ is independent of $n$ but may depend on
$u$. Clearly, if $h$ is $u$-regular then the function $h$ is
eventually, strictly decreasing along the sequence $u$. Thus, the
regularity condition imposes the condition that  $u$ is
eventually, strictly increasing.  Also, note that if $h$ is
$u$-regular then it is $s$--regular for any subsequence $s$ of
$u$.

\bigskip

\begin{theorem}\label{THM3}
 Let $(\Omega,d)$ be a compact metric space
equipped with a measure $m$ satisfying condition\/ {\rm(M2)}\/.
Suppose that  $(\cR,\beta)$ is a locally $\,m$-ubiquitous system
relative
 to $(\rho,l,u)$ and that $\psi$ is an approximation
function. Let $f$ be a dimension function such that $r^{-\delta}
\, f(r) \to \infty$ as $r \to 0$ and $r^{-\delta} \, f(r)$  is
decreasing. Furthermore, suppose that $r^{-\gamma} \, f(r) $ is
increasing. Let $g$ be the real, positive  function given by
\begin{equation}\label{my5}
g(r)  :=   f(\p(r))  \p(r)^{-\ga}  \r(r)^{\ga -\delta} \ \mbox{
and let } \ \  G \, := \, \limsup_{n \to \infty} \, g(u_n).
\end{equation}
\begin{enumerate}
\item[(i)]
Suppose that $G = 0 $  and that $\r$ is $u$-regular. Then,
\begin{equation}\label{my6}
 \mathcal{ H}^f \! \left( \Lambda(\p) \right)  =
 \infty \hspace{10mm} {\rm if } \hspace{10mm} \sum_{n=1}^{\infty}
g(u_n) =  \infty \; .
\end{equation}
\item[(ii)]
Suppose that $0< G \leq  \infty $. Then, $\mathcal{ H}^f \! \left(
\Lambda(\p) \right) \ = \ \infty$.
\end{enumerate}
\end{theorem}

\bigskip

\noindent{\em An important general observation: \ } In statements
such as Theorem \ref{THM3} in which the measure is of type (M2),
the lower sequence $l$ is actually redundant  from the  hypothesis
that $(\cR,\beta)$ is a locally $\,m$-ubiquitous system relative
to $(\rho,l,u)$.   The point is that the measure condition (M2) is
independent of the sequences $l$ and $u$.  Hence, in view of the
discussion in \S\ref{thesystems}, given an upper sequence $u$ for
which $$\textstyle
 m\left( B \cap  \bigcup_{\alpha\in J: \ma \leq u_n}
 \De(\ra,\r(u_n)) \right) \ \ge \ \ka \ m(B) \   ,
 $$
 a lower
sequence $l$ can always be constructed so that $(\cR,\beta)$ is a
locally $\,m$-ubiquitous system relative to $(\rho,l,u)$.

In the statement of Theorem \ref{THM1} the sequences $l$ and $u$
are determined by the measure condition (M1) as well as by
(\ref{af1}) -- (\ref{af}) and it is important that we have
ubiquity with respect to these sequences.

\section{Remarks  and corollaries to Theorem
\ref{THM1}\label{r12}}
Obviously the first conclusion of  Theorem \ref{THM1} is
significantly weaker than the other --- positive measure of
$\Lambda(\p)$ as opposed to full measure. However, in practice it
is much easier to establish `global ubiquity' than `local
ubiquity'.  Moreover, for certain applications it is possible to
establish the subsidiary result that $m(\Lambda(\p))$ is either
zero or one. For example, this is the case for the classical set
of $\p$-well approximable numbers -- see Theorem 2.7 of \cite{Har}
and indeed for the majority of applications considered in
\S\ref{app}.  Thus for such applications establishing
$m(\Lambda(\p))>0 $ is enough to deduce the full measure result
and in view of Theorem \ref{THM1} `global ubiquity' is all that is
required.

It will become evident during the course of establishing Theorem
\ref{THM1} that the intersection conditions associated with either
form of ubiquity are only required for $\alpha\in J$ with
$l_n<\ma\le u_n$ rather than for $\alpha \in J$ with
$\beta_{\alpha} \leq u_n$. Also, for the positive measure
statement of Theorem \ref{THM1} the doubling property of the
measure $m$ is only required for balls centred at resonant sets
rather than at arbitrary points of $\Omega$.

 It is easy to verify that the $\limsup$ condition
(\ref{af1}) implies the divergent sum condition (\ref{af2}). Thus,
whenever (\ref{af}) is satisfied  the $\limsup$ condition is
redundant. At first glance, (\ref{af}) may look like a horrendous
condition. Nevertheless, we shall see that it is both natural and
not particularly restrictive. For example, suppose throughout the
following discussion that the measure $m$ satisfies condition
(M2). Then the divergent sum condition (\ref{af2}) becomes
\begin{equation}  \label{dsm2}
  \sum_{n=1}^{\infty} \
 \left(\frac{\p(u_n)}{\r(u_n)}\right)^{\delta - \g} \ = \ \infty  \  ,
\end{equation}
and (\ref{af}) simplifies to
\begin{equation}  \label{afm1}
 \sum_{s=1}^{Q-1} \ \r(u_s)^{\g-\delta} \ \sum_{\substack{s+1 \, \leq \, t \leq \, Q
\, :\\ \p(u_s) \, < \, \r(u_t)}} \!\!\!\!\! \p(u_t)^{\delta-\g} \
\ll \   \left( \sum_{n=1}^{Q} \
 \left(\frac{\p(u_n)}{\r(u_n)}\right)^{\delta - \g} \right)^2  .
\end{equation}
Reiterating the earlier remark, trivially the  $\limsup$ condition
(\ref{af1}) implies (\ref{dsm2}). Thus, whenever (\ref{afm1}) is
satisfied (\ref{af1}) is redundant since (\ref{dsm2}) and
(\ref{afm1}) together  already imply the desired conclusions.

Consider for the moment  the special case that $\g=\delta $.
Trivially, both (\ref{dsm2}) and  (\ref{afm1}) are satsified.
Thus, Theorem \ref{THM1} reduces to:

\begin{corollary}
Let $(\Omega,d)$ be a compact metric space equipped with a measure
$m$ satisfying condition\/ {\rm(M2)}\/. Suppose that $(\cR,\beta)$
is a global $\,m$-ubiquitous system relative to $(\r,l,u)$ and
that $\, \p$ is an  approximating function. If $\g=\delta$ then $
m \left( \Lambda(\p) \right)  >  0$. In addition, if any open
 subset of $\Omega$ is $m$--measurable and  $(\cR,\beta)$ is
 locally $\,m$-ubiquitous
relative to $(\r,l,u)$, then $ m \left( \Lambda(\p) \right) = 1 $.
\label{cor1}
\end{corollary}

Next, suppose that the approximating function $\p$ is $u$-regular.
Then, for $t>s$ with $s$ sufficiently large  we have that $$
\p(u_t) \leq \lambda^{t-s}  \, \p(u_s)  \ \  $$ for some
$0<\lambda < 1$.  Without loss of generality, assume that $\delta
- \g > 0 $ -- the case $\delta = \g $ is covered by the above
corollary and since the measure $m$ satisfies condition (M2) we
have that $ \gamma \leq \dim \Omega = \delta$. Then for $Q$
sufficiently large, the  L.H.S.  of  (\ref{afm1}) is
 $$  \ll \
\sum_{s=1}^{Q-1} \
\left(\frac{\p(u_s)}{\r(u_s)}\right)^{\delta-\g} \ \sum_{s \, < \,
t \leq \, Q}   (\lambda^{\delta-\g})^{t-s} \ \ll \ \sum_{n=1}^{Q}
\ \left(\frac{\p(u_n)}{\r(u_n)}\right)^{\delta - \g}  \ .
 $$
This together with (\ref{dsm2}) implies that (\ref{afm1}) is
satisfied. We now consider the case that the ubiquity function
$\r$ is $u$-regular. It is easily verified that for $Q$
sufficiently large,
\begin{eqnarray*}
 {\rm L.H.S. \ of \  (\ref{afm1})} \ \, & \ll & \, \sum_{n=2}^{Q}
\ \p(u_n)^{\delta-\g} \ \sum_{m=1}^{n-1} \r(u_m)^{\g-\delta} \\ &
 &\hspace*{-10ex} \ \ll \   \sum_{n=2}^{Q}
\left(\frac{\p(u_n)}{\r(u_n)}\right)^{\delta-\g} \ \
\sum_{m=1}^{n-1}
 (\lambda^{\delta-\g})^{n-m}  \ \ll \
 \sum_{n=1}^{Q}
 \left(\frac{\p(u_n)}{\r(u_n)}\right)^{\delta - \g}  .
\end{eqnarray*}

\noindent This together with (\ref{dsm2}) implies that
(\ref{afm1}) is satisfied. On gathering together these
observations we have:

\begin{corollary}\label{cor2}
Let $(\Omega,d)$ be a compact metric space equipped with a measure
$m$ satisfying condition\/ {\rm(M2)}\/. Suppose that $(\cR,\beta)$
is a global $\,m$-ubiquitous system relative to $(\r,l,u)$ and
that $\, \p$ is an  approximating function. Furthermore, if
$\delta > \gamma$  suppose that either $\p$ or $\r$ is $u$-regular
and\/ {\rm(\ref{dsm2})} is satisfied. Then $ m \left( \Lambda(\p)
\right)  >  0 $. In addition, if any open  subset of $\Omega$ is
$m$--measurable and  $(\cR,\beta)$ is locally $\,m$-ubiquitous
relative to $(\r,l,u)$, then $ m \left( \Lambda(\p) \right) = 1 $.
\end{corollary}

In the numerous  applications considered in this paper, the
various ubiquitous functions $\r$ will always satisfy (\ref{afmh})
for the appropriate upper sequences $u$. Thus, if the measures are
also of type (M2) then Theorem \ref{THM1} simplifies to the above
corollaries.

\section{Remarks  and corollaries to Theorem
\ref{THM3}\label{r3}}
The case when $G$ is finite  constitutes the main substance  of
the theorem. When $G=0$, it would be desirable to remove the
condition that the ubiquity function $\r$ is $u$-regular -- see
below; in particular Corollary \ref{thm3*}. However, for the
numerous applications considered in this paper the hypotheses that
$\r$ is $u$-regular is always satisfied for the sequences $u$
under consideration. Clearly, the assumption that the function $0<
G \leq \infty $  in part (ii) implies the divergent sum condition
in part (i).

The case when the dimension function $f$ is $\delta$--dimensional
Hausdorff measure $\mathcal{ H}^{\delta}$ is excluded from the
statement of Theorem \ref{THM3} by the condition that $r^{-\delta}
\, f(r) \to \infty$ as $r \to 0$. This is natural since otherwise
Theorem \ref{THM1} implies that $m(\Lambda(\p))
> 0 $ which in turn implies that $\mathcal{ H}^{\delta} (\Lambda(\p))
$ is positive and finite -- see \S\ref{hmd}.  In other words
$\mathcal{ H}^{\delta} (\Lambda(\p)) $ is never infinite. To see that
Theorem \ref{THM1} is applicable, note that with $f(r)= r^{\d}$
the function $g$ of Theorem \ref{THM3} becomes $g(r) =
(\p(r)/\r(r))^{\delta - \gamma}$. Thus, if the sum in Theorem
\ref{THM3} diverges then so does the sum in Theorem \ref{THM1};
i.e. (\ref{af2}) is satisfied. Now, if $G =0 $ then (\ref{af}) is
satisfied since we assume that $\r$ is $u$-regular in part (i) of
Theorem \ref{THM3}. On the other hand, if $G > 0 $ then
(\ref{af1}) is satisfied. Thus, in either case Theorem \ref{THM1}
implies that $m(\Lambda(\p))
> 0 $.

Next  notice that if $\gamma=0$, then the hypothesis that
$r^{-\gamma} \, f(r) $ is increasing is redundant since, by
definition, any dimension function $f$ is increasing (cf.
\S\ref{hmd}). Also, notice that the conditions on $r^{-\gamma} \,
f(r) $ and $r^{-\d} \, f(r) $ exclude the possibility that $\g =
\d $. However, this is no great loss since if $\g = \d $ then
Corollary \ref{cor1} implies that $m(\Lambda(\p)) > 0 $. Here we
make use of the fact that `local' implies global' ubiquity.  Thus,
$\mathcal{ H}^{\delta} (\Lambda(\p))
> 0 $ and  is in fact finite. Now, $f$ is a dimension function such
that $r^{-\delta} \, f(r) \to \infty$ as $r \to 0$. It is
therefore a simple consequence of the elementary fact stated in
\S\ref{hmd} that $\mathcal{ H}^{f} (\Lambda(\p)) = \infty $. Thus
without loss of generality we can assume that $\gamma < \delta $
in Theorem \ref{THM3}. Finally,  note that in the case $\gamma
> 0$, if  both the functions  $\r$ and $g$ are  decreasing then one
always has that $r^{-\gamma} \, f(r) $ is   increasing.

As mentioned above, when $G=0$ it would be desirable to remove the
condition that the ubiquity function $\r$ is $u$-regular in
Theorem \ref{THM3}. At the expense of imposing growth conditions
on the functions $\p$ and $f$, the following result achieves
precisely this.

\begin{corollary}\label{thm3*}
 Let $(\Omega,d)$ be a compact metric space
equipped with a measure $m$ satisfying condition\/ {\rm(M2)}\/.
Suppose that  $(\cR,\beta)$ is a locally $\,m$-ubiquitous system
relative
 to $(\rho,l,u)$ and that $\psi$ is an approximation
function. Let $f$ be a dimension function such that $r^{-\delta}
\, f(r) \to \infty$ as $r \to 0$ and $r^{-\delta} \, f(r)$  is
decreasing. Furthermore, suppose that $r^{-\gamma} \, f(r) $ is
increasing. Let $g$ be the positive function given by\/
{\rm(\ref{my5})}.
\begin{enumerate}
\item[(i)]
Suppose that $G = 0 $  and that $\p$ is $u$-regular. Furthermore,
suppose there exist constants $r_0, \lambda_1,\lambda_2 \in (0,1)
$ so that for any $ r \in (0, r_0)$ one has $r^\g \, f(\lambda_1
r) \leq \lambda_2 \, f(r) (\lambda_1 r)^{\g}$. Then,
{\rm(\ref{my6})} is satisfied.
\medskip
\item[(ii)]
Suppose that $0< G \leq  \infty $. Then, $\mathcal{ H}^f \! \left(
\Lambda(\p) \right) \ = \ \infty$.
\end{enumerate}
\end{corollary}

\medskip

 \noindent{\it Proof of Corollary \ref{thm3*}.}~ Recall
that if the hypotheses of local ubiquity are satisfied for a
particular upper sequence $u$ then they are also satisfied for any
subsequence $s$ of $u$. The corollary follows from Theorem
\ref{THM3} by proving the existence of   an appropriate
subsequence $s$ of $u$ on which 
$\r$ is $s$-regular and $\sum g(s_n) = \infty $. To this end,
since $\p$ is $u$-regular there exists a constant $\lambda \in
(0,1) $ such that $ \psi(u_{n+1}) \leq \lambda \, \psi(u_n) $ for
all $n$ sufficiently large. Without loss of generality, we can
assume that $\lambda \leq \lambda_1$ -- see \S\ref{subseqwork}. In
view of the growth condition imposed on the dimension function $f$
and the fact that $r^{-\gamma}f(r)$ is increasing, we have that
for $n $ sufficiently large

$$ x_{n+1} = \frac{f(\psi(u_{n+1}))}{\psi(u_{n+1})^{\gamma}} \leq
\frac{f(\lambda\psi(u_{n}))}{\lambda\psi(u_{n})^{\gamma}} \le
\frac{f(\lambda_1\psi(u_{n}))}{\lambda_1\psi(u_{n})^{\gamma}} \leq
\lambda_2 \; \frac{f(\psi(u_{n}))}{\psi(u_{n})^{\gamma}} . $$

Hence, $ x_{n+1}\le\lambda_2 \, x_n$. Next, fix some sufficiently
large $n_1$ and for $k \geq 2 $ let $n_k$ be the least integer
strictly greater than $n_{k-1}$ such that $ \rho(u_{n_k}) \leq
\frac{1}{2} \rho(u_{n_{k-1}}) $. This is possible since $\r(r) \to
0$ as $ r \to \infty $. By construction, $\rho(u_{m}) \geq
\frac{1}{2} \rho(u_{n_{k-1}}) $ for any integer $ m \in [n_{k-1},
n_{k}-1 ] $. For shortness, let $\d':=\d-\g$. It follows that
\begin{eqnarray*}
\infty & = & \ \sum_{n=n_1}^{\infty} g(u_n) \  =  \
\sum_{n=n_1}^{\infty} x_n \; \rho(u_n)^{-\delta'}   =
\sum_{k=2}^{\infty} \ \ \sum_{n_{k-1} \leq m < n_k } \!\!\!\! x_m
\; \rho(u_m)^{-\delta'} \\ & \leq & \sum_{k=2}^{\infty}
\sum_{n_{k-1} \leq m < n_k } \!\!\!\! x_m \;
\rho(u_{n_{k-1}})^{-\delta'} \, 2^{\delta'}  =  2^{\delta'}
\sum_{k=2}^{\infty} \rho(u_{n_{k-1}})^{-\delta'} \!\!\!\!
\sum_{n_{k-1} \leq m < n_k } \!\!\!\!\!\!\!\! x_m \\ & \ll &
\sum_{k=2}^{\infty} \rho(u_{n_{k-1}})^{-\delta'}  x_{n_{k-1}} \
\sum_{i=0}^{\infty} \lambda_2^i
\ll \sum_{k=1}^{\infty} \rho(u_{n_{k}})^{-\delta'}  x_{n_{k
}}
 :=  \sum_{k=1}^{\infty} g(u_{n_{k}}) .
\end{eqnarray*}
Now set $s:= \{u_{n_k} \}$. By construction, $\rho$ is $s$-regular
and $\sum g(s_n)= \infty $. \QED

\medskip

Regarding Corollary \ref{thm3*},  the growth condition imposed  on
the function $f$ is not particularly restrictive. In particular,
when $f:r \to r^s$ and $s
> \gamma$ the growth condition is trivially satisfied. By
restricting our attention to $s$-dimensional Hausdorff measure
$\mathcal{ H}^s$, Theorem \ref{THM3} together with Corollary
\ref{thm3*} yield the following statement.

\medskip

\begin{corollary}\label{thm2cor}
 Let $(\Omega,d)$ be a compact metric space
equipped with a measure $m$ satisfying condition\/ {\rm(M2)}\/.
Suppose that  $(\cR,\beta)$ is a locally $\,m$-ubiquitous system
relative to $(\rho,l,u)$ and that $\psi$ is an approximation
function. Let $s \geq 0$ such that  $ \gamma < s< \delta$,
 $ g(r)  :=  \p(r)^{s-\gamma}
\r(r)^{\gamma - \delta }$  and let $G \, := \, \limsup_{n \to
\infty}  g(u_n)$.
\begin{enumerate}
\item[(i)]
Suppose that $G =0$  and that either $\p$ or $\r$ is $u$-regular.
Then,  $$
 \mathcal{ H}^s \! \left( \Lambda(\p) \right) \ =
\ \infty \hspace{10mm} {\rm if } \hspace{10mm} \sum_{n=1}^{\infty}
g(u_n)\ = \ \infty \ \ . $$
\item[(ii)] Suppose that $0< G \leq \infty$. Then, $\mathcal{ H}^s \!
\left( \Lambda(\p) \right) \ = \ \infty$.
\end{enumerate}
\end{corollary}

\medskip

  The following  lower  bound
statement for the dimension of $ \Lambda(\p)$ is essentially a
consequence of part (ii) of Corollary \ref{thm2cor}. The statement
is free of any regularity condition.

\medskip

\begin{corollary} \label{dimthm}
Let $(\Omega,d)$ be a compact metric space with a measure $m$
satisfying condition\/ {\rm(M2)}\/. Suppose  that $(\cR,\beta)$ is
a local $\,m$-ubiquitous system relative to $(\r,l,u)$ and let $
\p$ be an approximating function.
\begin{enumerate}
\item[(i)] If $ \delta \neq  \gamma$ and\/ $\lim_{n \to
\infty} \p(u_n)/\r(u_n) = 0$ then $$ \dim \Lambda(\p)\ge
d:=\gamma+\sigma(\delta-\gamma),\ \ {\rm where\  } \sigma  :=
\limsup_{n \to \infty} \frac{\log \rho(u_n)}{\log \p(u_n)}. $$
Moreover,  if\/ $ \liminf\limits_{n \to \infty} \rho(u_{n})/
\p(u_{n})^{\sigma} < \infty $, then $ \mathcal{ H}^d(\Lambda(\p))
= \infty $.
\item[]~~
\item[(ii)] If  either  $ \delta = \gamma$  or
$\limsup_{n \to \infty} \p(u_n)/\r(u_n)
> 0$ then  $ 0<  \mathcal{
H}^\delta(\Lambda(\p)) < \infty  \, $ and so  $ \dim \Lambda(\p) =
\delta $.
\end{enumerate}
\end{corollary}

It will be evident from the proof below that part (ii) of
Corollary \ref{dimthm} is a simple consequence of Theorem
\ref{THM1}. Part (i) contains the main substance  of Corollary
\ref{dimthm}. In order to establish part (i), all that is required
of Corollary \ref{thm2cor} is part (ii) and this explains why
there is no regularity condition on either $\p$ or $\rho$ in the
statement of Corollary \ref{dimthm}. Note that Corollary
\ref{thm2cor} implies both the dimension and measure statements of
part (i) and moreover provides a more general criteria for when
$\mathcal{ H}^d (\Lambda(\psi)) = \infty$. To illustrate this,
consider the classical  set $ W (\p) $ of $\psi$--well
approximable numbers with $\psi \!  : \! r \to r^{-\tau} (\log
r)^{-1} $ and $\tau > 2 $ -- see \S\ref{beg}.  The associated
ubiquitous system is given by Lemma~\ref{reg_be} in  \S\ref{congs}
-- in particular $\rho(r) := {\rm constant } \times r^{-2}$. Thus,
$\sigma = 2/\tau$ and  part (i) of  Corollary \ref{dimthm} implies
that $\dim W(\tau) \geq 2/\tau$ for $\tau > 2$. However, $ \lim_{n
\to \infty} \rho(2^{n})/ \p(2^{n})^{\sigma} = \infty$ and so we
obtain no information concerning $\mathcal{ H}^{2/\tau}
(\Lambda(\psi))$. On the other hand, since the function $\r$ is
$u$-regular Corollary \ref{thm2cor} implies the above  dimension
statement and shows that $\mathcal{ H}^{2/\tau} (\Lambda(\psi))=
\infty$. The dimension statement follows from the definition of
Hausdorff dimension -- see \S\ref{hmd}.

\noindent{\it Proof Corollary \ref{dimthm}.}~ To start with
consider the case $\delta = \gamma$, so $d := \gamma + \sigma
(\delta - \gamma) = \delta$. Since `local' implies `global'
ubiquity, Corollary \ref{cor1} of Theorem \ref{THM1} implies that
$m(\Lambda(\psi))
> 0 $. In turn,  Lemma \ref{covmdp} of  \S\ref{hmd} implies that $0 < \mathcal{
H}^\delta(\Lambda(\psi)) < \infty $  and $ \dim \Lambda(\psi) =
\delta $. This completes the proof in the case that $\delta =
\gamma$. 
 Next suppose,  $\limsup_{n \to \infty}
\p(u_n)/ \r(u_n) \, >  \, 0 $.   Then, Theorem \ref{THM1} implies
that $m(\Lambda(\psi)) > 0 $ and so $0 < \mathcal{
H}^\delta(\Lambda(\psi)) < \infty $  and $ \dim \Lambda(\psi) =
\delta $ for the same reasons as in the case $\delta = \gamma$.
The completes the proof of part (ii) of Corollary \ref{dimthm}.
Thus, without loss of generality we can assume that $ \lim_{n \to
\infty} \p(u_n)/ \r(u_n) \,= \, 0 $ and so $$ 0 \ \leq \
\limsup_{n \to \infty} \frac{\log \rho(u_n)}{\log \p(u_n)} \, :=
\,\sigma \, \leq 1 \ . $$ Regarding part (i), we first suppose
that there exists a strictly increasing sequence $\{n_i\}_{ i \in
\N }$ such that
\begin{equation} \lim_{i \to \infty} \frac{\rho(u_{n_i})}{
\p(u_{n_i})^{\sigma} } \ = L < \infty \ \ . \label{incseq1}
\end{equation}
Since $\lim_{n \to \infty} \p(u_n)/\r(u_n) = 0$  we have that
$\sigma < 1 $. This together with the fact that  $\delta
> \gamma$ implies that  $\g \leq d < \delta$. Now notice that
(\ref{incseq1}) implies that $$\lim_{i \to \infty}
\frac{\p(u_{n_i})^{d-\g}}{ \r(u_{n_i})^{\delta - \g}  } \ :=  \
\lim_{i \to \infty} g(u_{n_i})  \  =  \ L^{\gamma -\delta} \ > \ 0
 \ ,$$ and $g$ is precisely the function in Corollary
\ref{thm2cor} with $s=d$. Hence $G > 0  $  and part (ii) of
Corollary \ref{thm2cor} implies that $ \mathcal{ H}^d(\Lambda(\p)) =
\infty $, as required. Now suppose there is no sequence $\{n_i\}$
such that (\ref{incseq1}) is satisfied. Then $\sigma > 0 $ since
$\r(r) \to 0$ as $r \to \infty$. It follows from the definition of
$\sigma$ that for any $0 < \epsilon < \sigma$, there exists a
sequence $\{n_i\}_{i \in \N} $ such that (\ref{incseq1}) is
satisfied with $\sigma$ replaced by $\sigma_{\epsilon}:= \sigma -
\epsilon $. In fact, $L = 0 $ in (\ref{incseq1}). Thus, on
repeating the above argument with $\sigma$ replaced by
$\sigma_{\epsilon}$ and $d$ replaced by $d_{\epsilon} := \gamma +
\sigma_{\epsilon}(\delta - \gamma)$  we conclude that $\mathcal{
H}^{d_{\epsilon}} (\Lambda(\psi)) = \infty $. Hence $\dim
\Lambda(\psi) \geq d_{\epsilon} = d - \epsilon(\delta - \gamma) $.
On letting $\epsilon \to 0$, we obtain the  desired dimension
result. This completes the proof of Corollary \ref{dimthm}. \QED

\medskip

 The following corollary is a simple consequence of Corollary
\ref{dimthm}. In some sense it is no more than a slightly weaker,
alternative statement and is more in line with the original
`ubiquity result' of Dodson, Rynne \& Vickers \cite{DRV}.

\medskip

\begin{corollary}\label{dimcor}
Let $(\Omega,d)$ be a compact metric space equipped with a measure
$m$ satisfying condition\/ {\rm(M2)}\/. Suppose $(\cR,\beta)$ is a
local $\,m$-ubiquitous system relative to $(\r,l,u)$ and let $ \p$
be an approximating function. Then $\dim \Lambda(\p)  \geq  d : =
\gamma  +  \sigma  (\delta - \gamma)$, where $$
 \sigma  \ := \   \min \Big\{ 1, \Big|
\limsup\limits_{n \to \infty} \frac{\log \rho(u_n)}{\log \p(u_n)}
\Big|
  \Big\}  .
$$

\noindent Furthermore,  if $d < \delta$ and $ \liminf_{n \to
\infty} \rho(u_{n})/ \p(u_{n})^{\sigma} < \infty $, then $
\mathcal{ H}^d(\Lambda(\p)) = \infty $.
\end{corollary}

In the case that $(\Omega,d)$ is a bounded subset of $\R^n$, $m$\/
is the Lebes\-gue measure in $\R^n$ and $ \liminf_{n \to \infty}
\rho(u_{n})/ \p(u_{n})^{\sigma}  =   0 $, the above corollary is
essentially the `ubiquity result' of Dodson, Rynne \& Vickers
\cite{DRV}\footnote{The statement of the Theorem 1  in \cite{DRV}
is not  correct when $\delta = \gamma$ and the proof assumes that
$\lim_{r\to \infty} \p(r) = 0$. Also, a weaker statement of
Corollary \ref{dimcor}  appears as Theorem 5.6 in \cite{BD}.
However, the proof contains a flaw. The claim that a certain set
$T_{\infty}$ is a subset of $\Lambda(\p)$ is not necessarily true.
}. Also, in the case that $(\Omega,d)$ is a bounded subset of
$\R^n$, the resonant sets are points ($\gamma =0$) and $
\liminf_{n \to \infty} \rho(u_{n})/ \p(u_{n})^{\sigma} =   0 $,
the above corollary is essentially equivalent to Theorem 1 of
\cite{gang4}.

\medskip

We end this section  with a comment regarding Theorem \ref{THM3}
and the measure condition (M2).  Given the central conditions on
$m$ (namely that the $m$--measure of a ball $B(x,r)$ with $r>0$
and $x \in \Omega$ is strictly positive and that $m$ is doubling)
the main  property of (M2) that is utilized during the proof of
Theorem \ref{THM3}  is that $m(B(x,r))$ is comparable to a
function of $r$ alone and is  independent of $x$. In view of this,
consider the following measure condition.

\medskip

\noindent {\bf (M$2^{\prime}$) \ } There exists a positive
constant $r_o$   such that for any $ x \in \Omega $ and $ r \leq
r_o $, $$
 a \, {\bf m}(r)  \ \leq  \  m(B(x,r))  \ \leq  \   b \, {\bf m}(r)
   \  \label{m2'}
$$
 where ${\bf m}: \R^+ \to \R^+$ is an increasing, continuous
 function with ${\bf m}(r) \to 0 $ as $r \to 0$.

Clearly, with   $ {\bf m}(r) := r^{\delta} $  we obtain the
measure condition (M2).  By adapting  the proof of Theorem
\ref{THM3} in the obvious manner  we  obtain the more general
result:

\medskip

\begin{thdel}
Let $(\Omega,d)$ be a compact metric space equipped with a measure
$m$ satisfying condition \emph{(M{\em $2^{\prime}$})}.   Suppose
that $(\cR,\beta)$ is a locally $\,m$-ubiquitous system relative
 to $(\rho,l,u)$ and that $\psi$ is an approximation
function. Let $f$ be a dimension function such that $ f(r)/{\bf
m}(r) \to \infty$ as $r \to 0$ and $f(r)/{\bf m}(r)$  is
decreasing. Furthermore, suppose that  $r^{-\gamma} {\bf m}(r)$ is
increasing and that $r^{-\gamma} \, f(r) $ is increasing. Let $g$
be the function given by $$ g(r) := f(\p(r)) \; \p(r)^{-\g}
\r(r)^{\g}  \,  {\bf m}(\r(r))^{-1} \  \mbox{ and let } \   G  :=
\limsup_{r \to \infty} \, g(u_n)  . $$

\noindent (i)  Suppose that $G =0$   and that $\r$ is $u$-regular.
Then,  $$
 \mathcal{ H}^f \! \left( \Lambda(\p) \right) \ =
\ \infty \hspace{10mm} {\rm if } \hspace{10mm} \sum_{n=1}^{\infty}
g(u_n)\ = \ \infty \ \ . $$ (ii) Suppose that $0 < G \leq \infty$.
Then, $\mathcal{ H}^f \! \left( \Lambda(\p) \right) \ = \ \infty$.
\end{thdel}

\medskip

We have opted to work with condition (M2) rather than
(M$2^{\prime}$) simply for the sake of clarity and the  ease of
discussion. Furthermore, for the various applications considered
in \S\ref{app} if $m$ satisfies  (M$2^{\prime}$) then it always
satisfies (M2). Note that the extra condition that  $r^{-\gamma}
{\bf m}(r)$ is increasing in Theorem $\ref{THM3}^{\prime}$ plays
the role of the fact that $ \g \leq \delta$ in the case that $
{\bf m}(r) := r^{\delta} $.

\section{The classical results\label{basiceg}}
For the classical set $W(\p)$ of $\p$--well approximable numbers,
 Lemma~\ref{reg_be} in \S\ref{congs}
establishes local $m$-ubiquity. Clearly, the ubiquity function
$\r$ satisfies (\ref{afmh}) (i.e. $\r$  is $u$-regular) and so
Corollary \ref{cor2} establishes the divergent part of
Khintchine's Theorem. In fact  global $m$-ubiquity   would suffice
since to go from positive measure to full measure simply involves
making use of the subsidiary  result that $m(\Lambda(\p))$ is
either zero or one - see  Theorem 2.7 of \cite{Har}.  On the other
hand, Theorem \ref{THM3} establishes the divergent part of
Jarn\'{\i}k's Theorem. By making use of the `natural cover' of
$W(\p)$, the convergent parts of these classical results are
easily established.


\section{Hausdorff measures and dimension \label{hmd} }

A {\em dimension function} $f  :  \R^+ \to \R^+ $ is an
increasing, continuous  function such that $f(r)\to 0$ as $r\to 0
\, $. The Hausdorff $f$--measure with respect to the  dimension
function $f$ will be denoted throughout by $\mathcal{ H}^{f}$ and
is defined as follows. Suppose $F$ is  a non--empty subset  of
$(\Omega,d)$. For $\rho > 0$, a countable collection $
\left\{B_{i} \right\} $ of balls in $\Omega$ with radii $r_i  \leq
\rho $ for each $i$ such that $F \subset \bigcup_{i} B_{i} $ is
called a {\em $ \rho $-cover}\/ for $F$. Clearly such a cover
always exists for totally bounded metric spaces. For a dimension
function $f$ define $ \mathcal{ H}^{f}_{\rho} (F) \, = \, \inf
\left\{ \sum_{i} f(r_i) \ : \{ B_{i} \}  {\rm \ is\ a\ } \rho {\rm
-cover\ of\ } F \right\} \,
 $, where the infimum is over all $\rho$-covers.  The {\it
Hausdorff $f$--measure} $ \mathcal{ H}^{f} (F)$ of $F$ with respect to
the dimension function $f$ is defined by   $$ \mathcal{ H}^{f} (F) :=
\lim_{ \rho \rightarrow 0} \mathcal{ H}^{f}_{\rho} (F)  \; = \;
\sup_{\rho > 0 } \mathcal{ H}^{f}_{\rho}  (F) \; . $$

A simple consequence of the definition of $ \mathcal{ H}^f $ is the
following useful

\medskip

\begin{lemma}
If $ \, f$ and $g$ are two dimension functions such that the ratio
$f(r)/g(r) \to 0 $ as $ r \to 0 $, then $\mathcal{ H}^{f} (F) =0 $
whenever $\mathcal{ H}^{g} (F) < \infty $.
\end{lemma}

In the case that  $f(r) = r^s$ ($s \geq 0$), the measure $ \hf $
is the usual {\em $s$--dimensional Hausdorff measure}\/ $\hs $ and
the Hausdorff dimension $\dim F$ of a set $F$ is defined by $$
\dim \, F \, := \, \inf \left\{ s : \mathcal{ H}^{s} (F) =0 \right\} =
\sup \left\{ s : \mathcal{ H}^{s} (F) = \infty \right\} . $$ In
particular when $s$ is an integer $\hs$ is comparable to
$s$--dimensional Lebesgue measure. For further details see
\cite{det,falc}. A general and classical  method for obtaining a
lower bound for the Hausdorff $f$-measure of an arbitrary set $F$
is the following mass distribution principle.

\medskip

\begin{lemma}[Mass Distribution Principle]
Let $\mu$ be a probability measure supported on a subset $F$ of
$(\Omega,d)$. Suppose there are  positive constants $c$ and $r_o$
such that for any ball $B$ with radius $r \leq r_o$ $$ \mu(B) \le
c \;  f(r)\; . $$ If $X$ is a subset of $F$ with $\mu(X) = \lambda
> 0$ then $ \mathcal{ H}^{f} (X) \geq \lambda/c \; $.
\end{lemma}

\noindent{\em Proof.}~ If $ \left\{B_{i} \right\} $ is a
$\rho$--cover of $X$ with $\rho \leq r_o$ then $$ \lambda = \mu(X)
= \mu \left( \textstyle{ \bigcup\limits_{i} } B_i \right) \leq
\sum_i \mu \left( B_i \right) \leq c \sum_i f(r_i) \; . $$ It
follows that $ \mathcal{ H}^{f}_{\rho } (X) \geq \lambda/c $ for any
$\rho \leq r_o$. On letting $\rho \to 0 \, $, the quantity $\mathcal{
H}^{f}_{ \rho } (X)$ increases and so we obtain the required
result. \QED


The following rather simple covering result will be used at
various stages during  the proof of our theorems.

\medskip

\begin{lemma}[Covering lemma]
Let $(\Omega,d)$ be a  metric space and $\B$ be a finite
collection of balls with common radius $r>0$. Then there exists a
disjoint sub-collection $\{B_i\}$ such that $$\textstyle
\bigcup\limits_{B\in\B} \, B \ \subset \ \bigcup\limits_i \
{3B}_i \ .
$$
\end{lemma}

\noindent{\em Proof.}~ Let $S$ denote the set of centres of the
balls in $\B$.  Choose $c_1 \in S$ and for $k \geq 1 $,
$$\textstyle c_{k+1} \ \in \ S \ \backslash \  \bigcup\limits_{i=1}^{k}
B(c_i, 2r) $$ as long as $S  \ \backslash \ \bigcup_{i=1}^{k}
B(c_i, 2r) \neq \emptyset $. Since $\#S$ is finite, there exists
$k_1 \leq \#S $ such that $
 S \ \subset \  \bigcup_{i=1}^{k_1} B(c_i, 2r) \ .
 $ By
construction, any ball $B(c,r)$ in the original collection $\B$ is
contained in some ball $B(c_i,3r)$ and since $d(c_i,c_j) > 2r$
 the chosen balls $ B(c_i,r) $ are clearly  disjoint.
\QED

\medskip

We end this section by making use of the mass distribution
principle and the covering lemma to establish the following claim
mentioned in \S\ref{congs}.

\medskip

\begin{lemma}
Let $(\Omega,d)$ be a totally bounded metric space equipped with a
probability measure $m$ satisfying condition\/ {\rm(M2)}\/.   Then
for any $X \subseteq \Omega $ with $m(X)
> 0 $
$$ 0 \ < \ \mathcal{ H}^{\delta} (X) \ < \ \infty \hspace{9mm} {\rm
and \ } \hspace{9mm} \dim X = \delta \ .    $$ \label{covmdp}
\end{lemma}

\noindent{\em Proof.}~ Given the measure statement, the dimension
statement follows directly from the definition of Hausdorff
dimension. The fact that $\mathcal{ H}^{\delta}  (X)$ is strictly
positive is a simple consequence of the mass distribution
principle with $f(r)= r^\delta$. Thus the lemma follows on showing
that $\mathcal{ H}^{\delta}  (\Omega)$ is finite since $\mathcal{
H}^{\delta}  (X) \leq \mathcal{ H}^{\delta}  (\Omega)$. Since the
metric space is totally bounded, for any $\rho > 0 $ there exists
a finite collection $\B$ of balls $B(\rho)$ with centres in
$\Omega$ and common radius $\rho$ such that $$ \Omega \; \subset
\; \bigcup\limits_{B(\rho)\in\B} \ B(\rho) \ . $$ In other words,
$\B$ is a $\rho$--cover of $\Omega$. By the covering lemma, there
exists a sub-collection $\{B_i(\rho)\}$ such that $$
 \bigcup_{i}^{\circ} \,
B_i(\rho) \ \subset \ \bigcup_{B(\rho)\in\B} \, B(\rho) \ \subset
\ \bigcup_i \ B_i(3\rho) \ , $$ where the left hand union is
disjoint.  Thus the collection $\{ B_i(3\rho)\} $ is a
$3\rho$-cover of $\Omega$. Hence, $$
 \mathcal{ H}_{3\rho}^{\delta}  (\Omega)  \le   \sum_i(3\rho)^{\delta}
 \ll  \sum_i m(B_i(\rho))  =  m \left(
\bigcup_{i}^{\circ}B_i(\rho) \right)  \leq   m(\Omega)  =  1  . $$
On letting $\rho \to 0 \, $, we obtain that $ \mathcal{ H}^{\delta}
(\Omega) \, \ll \, 1 $ as required. \QED

\smallskip


\section{Positive and full $m$--measure sets}

 The aim of this section is to  determine  conditions under which
a subset $E$ of a metric measure space $(\Omega,d,m)$ has full
$m$-measure. We also state two important lemmas which enable us to
conclude that $m(E) $ is strictly positive in the case that $E$ is
a $\limsup$ set. It is worth mentioning that the lemmas are in
fact a generalization of results known for Lebesgue measure since
the 1920s --  see, for instance  \cite{knopp}.

\medskip

\begin{lemma}\label{lem1_0}
Let $(\Omega,d)$ be a metric space and let $m$ be a finite measure
on $\Omega$ such that any open set is $m$-measurable. Let $E$ be a
Borel subset of\/ $\Omega$ and $f:\Rp\to\Rp$ be increasing with
$f(x)\to0$ as $x\to0$. Assume that $$ m(E\cap U)\ge f(m(U))\
\mbox{ for any open set }U\subset\Omega. $$ Then $E$ has full
measure in $\Omega$; i.e.\ $m(\Omega\setminus E)=0$.
\end{lemma}

\noindent{\em Proof.}~ Assume  that   $m(\Omega\setminus E)>0$.
Thus there exists $\ve
>0$ such that $f(m(\Omega\setminus E)) > \ve$. Notice that
$m(E)\le m(\Omega)<\infty$. Then $E$ contains a closed set $C$ for
which  $m(E\setminus C)<\ve$ (see Theorem~2.2.2. \cite{Federer}).
Since $C$ is closed, the set $U=\Omega\setminus C$ is open and
thus $m(U\cap E)\ge f(m(U))$. As $C\subset E$, we have
$\Omega\setminus E\subset U$. Then $m(\Omega\setminus E)\le m(U)$
and by the monotonicity of $f$, $$ m(U\cap E)\ge f(m(U))\ge
f(m(\Omega\setminus E))>\ve. $$ However, $E\setminus C=U\cap E$.
So the previous set of inequalities contradicts the fact that
$m(U\cap E)=m(E\setminus C)<\ve$. \QED

\medskip

\begin{lemma}\label{lem1}
Let $(\Omega,d)$ be a metric space and let $m$ be a finite measure
on $\Omega$ such that any open set is $m$-measurable. Let $E$ be a
Borel subset of\/ $\Omega$. Assume that there are constants
$r_0,a>0$ such that for any ball $B$ with centre in $ \Omega$ of
radius $r(B) <r_0$ $$ m(E\cap B)\ge a\, m(B)  \  .  $$

\noindent Furthermore, assume that there is a constant $b>0$ such
that for every open set $U$ there is a  finite or countable
disjoint collection $G$ of  balls contained in $U$ with centres in
$\Omega$ and  radii less than $r_0$ such that
\begin{equation}\label{v1}
b\, m(U)\le \sum_{B\in G}m(B).
\end{equation}
Then $E$ has full measure in $\Omega$;  i.e.\ $m(\Omega\setminus
E)=0$.
\end{lemma}

\noindent{\em Proof.}~ Take any open set $U$. Then there is a
disjoint collection $G$ of balls in $U$ with centres in $\Omega$
and radii less than $r_0$ satisfying (\ref{v1}). It follows that $
m(U\cap E)\ge m\Big(\bigcup_{B\in G} B\cap E\Big)=\sum_{B\in
G}m(B\cap E) \ge\sum_{B\in G}a\,m(B) \geq a \, b \,m(U) . $
Applying Lemma \ref{lem1_0} with $f(x)=a b \, x$ \ completes the
proof. \QED

\medskip

In short, condition (\ref{v1}) implies that any open set $U$ can
be substantially packed with sufficiently small disjoint balls
centred in $\Omega$.

\bigskip

At this point we introduce various notions which can all  be found
in \S 2.8 of \cite{Federer}. Let $F$ be a family of closed subsets
of $\Omega$. We say that $F$ covers $A\subset\Omega$\/ {\it
finely}\/ if for any $a\in A$, $\ve>0$ there is a set $S$ in $F$
such that $a \in S$ with  $S$  contained in the open ball
$B(a,\ve)$. Next, $F$ is said to be $m$-{\it adequate}\/ for $A$
if for each open subset $V$ of $\Omega$ there is a countable
disjoint subfamily $G$ of $F$ with $$ \bigcup_{S\in G}S\subset
V\mbox{ and }m\Big((V\cap A)\setminus\bigcup_{S\in G}S\Big)=0. $$
Finally, we associate with each subset $S$ of $F$ its
$(\delta,\tau)$-{\it enlargement}\/ defined by
\begin{equation}\label{v2}
\hat S=\bigcup_{T\in F,\ T\cap S\not=\emptyset,\ \delta(T)\le\tau
\,  \delta(S)}T \, .
\end{equation}
Here   $\tau\in\Rp$ and $\delta$ is a non-negative bounded
function on $F$. The following lemma brings together the above
notions and appears as Theorem 2.8.7 in \cite{Federer}.

\medskip

\begin{adeq}
Assume that $m$ is a finite measure on $\Omega$ such that every
open set in $\Omega$ is $m$-measurable. If $F$ covers
$A\subset\Omega$ finely, $\delta$ is a non-negative bounded
function on $F$, $1<\tau,\lambda<\infty$ and
\begin{equation}\label{v3}
m(\hat S)\le\lambda\, m(S)
\end{equation}
whenever $S\in F$ and $\hat S$ is the $(\delta,\tau)$-enlargement
of $S$, then $F$ is $m$-adequate for $A$.
\end{adeq}

We now consider the case that $F$ is the family of closed balls
$B$ with centres in $\Omega$. With reference to Lemma \ref{lem1},
we assume that the radius $r(B)$ of any  ball in $F$ is less than
$r_o$.  It follows that $F$ covers finely every subset of $\Omega$
and in particular any open subset $U$. Next, consider the function
$ \delta  : B \to \delta(B)$ where $\delta(B)$ is the diameter of
$B$ in $F$. Then $ S =B \in F$ implies $\hat S \subset (1+2\tau)
\, B. $ Thus, if the measure $m$ satisfies the {\it diametric
regularity condition}\,:
\begin{equation}\label{v4}
m((1+2\tau)\, B ) \  \le \ \lambda  \, m(B)  \ ,
\end{equation}
then (\ref{v3}) is satisfied. Clearly, if  $m$ is doubling then
the diametric regularity condition (\ref{v4}) is trivially
satisfied. Recall, that $m$ is said to be doubling if there exists
a constant $C \geq 1 $ such that $\forall \ x \in \Omega$  \ \ \ $
m(B(x,2r)) \, \leq \, C \ m(B(x,r)) . $ Note that the doubling
condition is independent of whether the ball $B(x,r)$ is open or
closed. It is easy to see that if $m$ is doubling then the $m$
measure of a closed ball $B(x,r)$ with $x$ in $\Omega$ is
comparable to that of the open ball $B^{\circ}(x,r)\;\!$: $$
\mbox{\small $\frac{1}{C}$} \,  m(B) \leq m(\mbox{\small
$\frac{1}{2}$}B) \leq m(B^{\circ}) \leq m(B) \ . $$ The upshot of
this is that if $F$ and $\delta$ are as above and $m$ is doubling,
then $F$ is $m$-adequate for any open subset $U$ of $\Omega$. By
definition, there is a  finite or countable disjoint family
$G\subset F$ of closed balls with radii less than $r_o$ such that
$$ \bigcup_{B\in G}B\subset U\mbox{ and }
m\Big(U\setminus\bigcup_{B\in G}B\Big)=0. $$ Hence $$
m(U)=\sum_{B\in G}m(B). $$

Thus the family $G$ of balls satisfies (\ref{v1}) and we obtain
the following useful modification of Lemma \ref{lem1}.

\medskip

\begin{proposition}\label{lem1b}
Let $(\Omega,d)$ be a metric space and let $m$ be a finite,
doubling  measure on $\Omega$ such that any open set is
measurable. Let $E$ be a Borel subset of\/ $\Omega$. Assume that
there are constants $r_0,c>0$ such that for any ball $B$ of radius
$r(B)<r_0$ and centre in $\Omega$  we have that $$m(E\cap B)\ge c\
m(B) \ . $$ Then $E$ has full measure in $\Omega$, i.e.\
$m(\Omega\setminus E)=0$.
\end{proposition}

Proposition~\ref{lem1b} will be used in the proof of Theorem
\ref{THM1} to show the full measure result. The following two
propositions on the $m$--measure of $\limsup$ sets will be
required to establish that $m(\Lambda(\p))>0$.

\medskip

\begin{proposition}\label{lem2}
Let $(\Omega,A,m) $ be a probability  space and $E_n \in A$ be a
sequence of sets such that $\sum_{n=1}^\infty m(E_n)=\infty $.
Then
$$ m( \limsup_{n \to \infty} E_n ) \; \geq \;  \limsup_{Q \to
\infty} \frac{ \left( \sum_{s=1}^{Q} m(E_s) \right)^2 }{ \sum_{s,
t = 1}^{Q} m(E_s \cap E_t ) }  \  \  \ . $$
\end{proposition}

This result is a generalization of the divergent part of the
standard Borel--Cantelli lemma. For the proof see either
\cite[Lemma 2.3]{Har} or \cite[Lemma 5]{Sp}.

\medskip

\begin{proposition}\label{lem0}
Let $(\Omega,A,m) $ be a probability  space, $F\in A$ and $E_n \in
A$ a sequence of  sets. Suppose there exists a constant $c > 0$
such that $\limsup_{n \to \infty} m(F \, \cap \, E_n) \ge c \;
m(F)$. Then $$ m (F \, \cap  \, \limsup_{n\to\infty}E_n ) \ \ge \
c^2 \, m(F) \ \ . $$
\end{proposition}

\noindent{\em Proof.}~ Without loss of generality assume that
$m(F) > 0$. For any $0< \ve < c$, there is a subsequence $E_{n_i}$
with $n_i$ strictly increasing such that $m(F \cap E_{n_i})\ge
(c-\ve) \: m(F)$. Clearly $$\textstyle \Big(\sum\limits_{i=1}^N
m(F \cap E_{n_i})\,\Big)^2\ge \Big(\sum\limits_{i=1}^N (c-\ve) \;
m(F) \, \Big)^2= (c-\ve)^2N^2 \, m(F)^2 $$ and $$\textstyle{
\sum\limits_{n,m=1}^N m(F \cap E_n\cap E_m)\le
\sum\limits_{m,n=1}^N m(F)= m(F) \, N^2 } . $$ Also notice  that
$\sum_{i=1}^\infty m(F \cap E_{n_i}) \geq m(F) \,
\sum_{i=1}^\infty(c-\ve) =\infty$. Thus on applying Proposition
\ref{lem2} and observing that $ F \, \cap \,
\limsup_{n\to\infty}E_n   \supseteq     F \, \cap  \,
\limsup_{i\to\infty}E_{n_i} $ we have that $$ m\Big(F \cap
\limsup_{n\to\infty}E_n\Big)\ge
\limsup_{N\to\infty}\frac{(c-\ve)^2N^2 m(F)^2}{m(F)N^2}
=(c-\ve)^2\, m(F) \ . $$ As $\ve>0$ is arbitrary,  this completes
the proof of the proposition. \QED



\section{Proof of Theorem \ref{THM1}} \label{proof1}

Let $B$ be an arbitrary ball centred at a point in $\Omega$.  The
aim is to show that
\begin{equation}\label{my1}
m (\Lambda(\p)\cap B) \ \ge \  m(B)/C \ ,
\end{equation}
where  $C>0$  is a constant independent of $B$.

Under the global ubiquity hypothesis,  Theorem \ref{THM1} follows
on  establishing  (\ref{my1}) with $B:=\Omega$ -- the space
$\Omega$ can be regarded as a ball since it is compact. In the
case of local ubiquity, (\ref{my1}) will be established  for balls
$B$ with sufficiently small radii  so that the conditions of local
ubiquity and  Proposition~\ref{lem1b} are fulfilled. Then
(\ref{my1}) together with Proposition~\ref{lem1b} clearly implies
Theorem~\ref{THM1} for local ubiquity. Since in the `local' case
we appeal to Proposition~\ref{lem1b}, the extra  hypothesis that
any open subset of $\Omega$ is $m$-measurable  is necessary.

In view of the above discussion, let $B(x,r)$ be a ball for which
(\ref{ub})  is satisfied. In order to establish (\ref{my1}), we
begin by constructing a `good' subset $A(\p,B)$ of
$\Lambda(\p)\cap B$. Essentially, each thickening $\De(\ra,
\p(\beta_{\alpha}))$ of a resonant set $\ra$ will be replaced by
carefully chosen collections of balls contained in the set
$\De(\ra, \p(\beta_{\alpha}))$ with centers on $\ra$. In the case
that the resonant sets are points ($\gamma=0$), so that the
thickenings themselves are already balls, the argument is much
simplified but still crucial.

\subsection{The subset $A(\p,B)$ of $\Lambda(\p)\cap B
$\label{proof1.1}}
Since $(\Omega,d)$ is totally bounded, we can cover $\Omega$ by a
finite collection of balls $\tB$ with common radius $\r(u_n)$.
Suppose $\tB \cap \De(\ra,\r(u_n)) \neq \emptyset $. Then there
exists a point $c \in \ra$ such that $\tB \subset B(c,3\r(u_n))$.
Thus, for each $\alpha\in J^u_l(n)$, there is a finite cover of
$\De(\ra,\r(u_n))$ by balls $ B(c,3\r(u_n))$ centred at points $c$
on $\ra$. This statement is of course obvious in the case when the
resonant sets are points.  Denote by $G_{\mbox{\tiny
$\Omega$}}(n,\alpha)$ the collection of centers $c \in \ra$ of the
balls $B(c,3\r(u_n))$ and by $G^{*}_{\mbox{\tiny $\Omega$}}(n)$
the set of all such centers as $\alpha$ runs through $J^u_l(n)$;
that is $G^{*}_{\mbox{\tiny $\Omega$}}(n):=\{c \in G_{\mbox{\tiny
$\Omega$}}(n,\alpha) \, : \, \alpha\in J^u_l(n) \}$. In the case
that $c$ lies on more than one $\ra$ simply choose one of them.
Clearly, the  collection of balls $ \B^{*}_{\mbox{\tiny
$\Omega$}}(n):=\{B(c,3\r(u_n)):c\in G^{*}_{\mbox{\tiny
$\Omega$}}(n)\}$  is a cover for $\De^u_l(\rho,n) :=
\bigcup_{\alpha\in J^u_l(n)} \De(\ra,\r(u_n)) $. In view of the
covering lemma, there exists a disjoint sub-collection
$\tilde\B_{\mbox{\tiny $\Omega$}}(n)$ of $\B^{*}_{\mbox{\tiny
$\Omega$}}(n)$ with centers $c\in  \tilde G_{\mbox{\tiny
$\Omega$}}(n)$ such that
\begin{equation}\label{my2}
 \bigcup_{c \, \in
\tilde G_{\mbox{\tiny $\Omega$}}(n)}^\circ B(c,\r(u_n)) \ \subset
\ \De^u_l(\rho,n) \ \subset \ \bigcup_{c \, \in \tilde
G_{\mbox{\tiny $\Omega$}}(n)} B(c,9\r(u_n)) \ .
\end{equation}
The left hand side follows from the fact that the balls $B(c, 3
\r(u_n))$ with $c\in \tilde G_{\mbox{\tiny $\Omega$}}(n)$ are
disjoint and that $B(c,\r(u_n))\subseteq\De(\ra,\r(u_n))$ for any
point $c$ on $\ra$.

Choose $n$ sufficiently large so that $36\r(u_n) < r $ (by
definition,  $ \r(u_n) \to  0$ as $n \to \infty$) and let $$
\gb(n) \, := \, \left\{ c \in \tilde G_{\Omega}(n) \ : \ c \in
\mbox{\small $\frac{1}{2}$} \,  B \right\} \ \ . $$ Now by
definition and (\ref{my2}), $$ \bigcup_{c \, \in \gb(n)}^\circ
B(c,\r(u_n)) \ \subset \ \De^u_l(\rho,n) \ \cap \ B $$ and $$
\bigcup_{c \, \in \gb(n)} B(c,9\r(u_n))) \ \supset \
\De^u_l(\rho,n) \ \cap \ \mbox{\small $\frac{1}{4}$} \, B  \ \ .
$$

We now estimate the cardinality of $\gb(n)$. By (\ref{ub}) and the
fact that the measure $m$ is of type (M1), for $n$ sufficiently
large
\begin{eqnarray*}
\#\gb(n) & \, & \!\!\!\!\!\!\!\!\!\!\!\!\!\!  m(B_n(\r(u_n))) \
\gg  \ m \textstyle{\left(\bigcup_{c\in \gb(n)}
B(c,9\r(u_n))\right) } \ \geq \\ & \ge &
 \!\!\!  m\left( \De^u_l(\rho,n)  \cap  \mbox{\small $\frac{1}{4}$} \,
B \right) \ \
 \ge \ \ka \, m(\mbox{\small $\frac{1}{4}$} \, B) \ \gg \
\ka \, m(B)
 \ .
\end{eqnarray*}
where $B_n(r)$ is a generic ball centred at a point $c\in
G_{\mbox{\tiny $\Omega$}}(n)$ of radius $r$.  On the other hand
\begin{eqnarray*}
  m(B)
&\ge& m\left(\bigcup_{c\in \gb(n)}^\circ B(c,\r(u_n))\right)
\ \gg \ \#\gb(n) \  m(B_n(\r(u_n)))   \ ,
\end{eqnarray*}
where the implied constant is dependent only on the constant $a$
of (\ref{measurem1}). The upshot of this is that
 \begin{equation}\label{NGn2}
 \#\gb(n) \ \asymp \
 \frac{m(B)}{ m(B_n(\r(u_n)))} \ \ .
\end{equation}
In the case $B=\Omega$, (\ref{NGn2}) is satisfied with $m(B)$
replaced by $m(\Omega):=1$.

We are now already in the position to prove the theorem under the
 $\limsup$ hypothesis (\ref{af1}).  Suppose for some sufficiently
large $n \in \N$ we have that $\p(u_n) \geq k \, \r(u_n)$, where
$k>0$ is a constant. If $k \geq 1$, (\ref{ub})  implies that $$
m(\De^u_l(\p,n) \cap B ) \, \geq \, m(\De^u_l(\rho,n) \cap B ) \,
\geq \, \ka \, m(B)\ . $$ On the other hand if  $ \r(u_n) >
\psi(u_n) > k \, \r(u_n) $, then (\ref{NGn2}) together with the
fact that $m$ is doubling and that $k<1$ implies that
\begin{eqnarray*}
 m(\De^u_l(\p,n) \cap B )  & \geq & \, m\Big(\bigcup\limits^{\circ}_{c\in
\gb(n)} B(c,\p(u_n))\Big)   \\ & \gg &  \ \#\gb(n) \
m\left(B_n(\p(u_n))\right) \  \gg \  m(B) \ .
\end{eqnarray*}
Thus, if
$\p(u_n) \geq k \, \r(u_n)$ for infinitely many $n \in \N$,
Proposition \ref{lem0} with $F=\Omega$ implies (\ref{my1}) and
thereby completes the proof of Theorem~\ref{THM1}. It remains to
establish the theorem under the hypotheses (\ref{af2}) and
(\ref{af}). Moreover, given any constant $k > 0 $ we can assume
without loss of generality that for $n$ sufficiently large
\begin{equation}\label{my3}
\r(u_n) >  k \; \p(u_n)\; .
\end{equation}

By definition, for each $c\in \gb(n)$ there exists an $\alpha \in
J^u_l(n)$ such that $c\in \ra \cap \mbox{\footnotesize $
\frac{1}{2}$} B$ or simply that $c\in \ra $ when $B= \Omega$.
Assume for the moment that $\ga
> 0 $. Cover the set $$ B(c, \textstyle \frac{1}{2}
\r(u_n))\cap\De(\ra,\p(u_n)) \ \subset \ B $$ by balls $\tB$ of
common radius $\p(u_n)$. Suppose $\tB \cap B(c, \textstyle
\frac{1}{2} \r(u_n))\cap\De(\ra,\p(u_n)) \neq \emptyset $. Then
there exists some $c' \in \ra$ such that $\tB \subset B(c',
3\p(u_n))$. Let $\B^{*}_{\mbox{\tiny B}}(n,c)$ denote the
collection of balls $B(c', 3\p(u_n))$ arising in this way. Clearly
this collection of balls centred at points on $\ra$ is a cover for
the set  $ B(c, \mbox{\footnotesize $ \frac{1}{2} $}
\r(u_n))\cap\De(\ra,\p(u_n))$. Note that since we can assume that
$ \r(u_n)> 24 \, \p(u_n) $ for $n$ large enough, the collection
$\B^{*}_{\mbox{\tiny B}}(n,c)$ of balls is contained in $B(c,
\mbox{\footnotesize $ \frac{3}{4} $} \r(u_n))$.

By the covering lemma, there is a disjoint sub-collection
$\B_{\mbox{\tiny B}}(n,c)$ of $\B^{*}_{\mbox{\tiny B}}(n,c)$ with
centers $c'\in \gb(n,c)$ such that for $n$ sufficiently large
\begin{equation}\label{deltapsiB}
 B(c, \mbox{\footnotesize
$ \frac{1}{2} $} \r(u_n)) \, \cap \,   \De(\ra,\p(u_n))\ \subset \
\bigcup_{c'\in \gb(n,c)} B(c',9\p(u_n)) ,
\end{equation}
and
\begin{equation}\label{2deltapsiB}
\bigcup_{c'\in \gb(n,c)}^\circ B(c',\p(u_n)) \ \subset \
B(c,\mbox{\footnotesize $ \frac{3}{4} $}\r(u_n))\cap
\De(\ra,\p(u_n))\ .
\end{equation}
Now, (\ref{deltapsiB}) together with intersection condition (i),
implies that
\begin{eqnarray*} \#G_{\mbox{\tiny B}}(n,c) \
m\left(B_n(\p(u_n))\right) &\gg&
 m \left( \bigcup_{c'\in G_{\mbox{\tiny B}}(n,c)}
B(c',9\p(u_n))  \right)\\
&&\\
& \geq & m\left(
B(c,\mbox{\footnotesize $ \frac{1}{2} $}\r(u_n)) \cap
\De(\ra,\p(u_n)) \right)\\
&&\\
& \gg & m \left(B_n(\p(u_n))\right) \times  \left(
\frac{\r(u_n)}{\p(u_n)} \right)^\gamma  .
\end{eqnarray*}  Similarly,
(\ref{2deltapsiB}) together with intersection condition (ii)
implies that
\begin{eqnarray*}
 \#G_{\mbox{\tiny B}}(n,c) \ m(B_n(\p(u_n))) &\asymp&
m \left( \bigcup_{c'\in G_{\mbox{\tiny B}}(n,c)}^\circ
B(c',\p(u_n))) \right) \\
&&\\
&\le&   m\left( B(c,\textstyle
\frac{3}{4}\r(u_n)) \cap \De(\ra,\p(u_n)) \right)\\
&&\\
&\ll &  m
(B_n(\p(u_n))) \times \left( \frac{\r(u_n)}{\p(u_n)}
\right)^\gamma  . \end{eqnarray*} Hence,
\begin{equation}\label{NGncl}
 \# \gb(n,c) \, \asymp \, \left(
\frac{\r(u_n)}{\p(u_n)} \right)^\gamma \ .
\end{equation}

In the case $\ga = 0$, we define $\gb(n,c):= \{c\}$ and so $\#
\gb(n,c) = 1$. Thus  (\ref{NGncl}) is satisfied even when $\gamma
= 0$. Now let $$ A_n(\p,B) \, := \, \bigcup_{c\in
\gb(n)}\bigcup_{c'\in \gb(n,c)} B(c',\p(u_n)) \  \ .$$

\noindent It is easily verified  that the balls in the above
definition of $A_n(\p,B)$ are disjoint. Indeed, for any $c \in
G_{\mbox{\tiny B}}(n)$ the balls $B(c',\p(u_n)) $ with $c'\in
G_{\mbox{\tiny B}}(n,c)$ are disjoint. Also, for $c_1,c_2 \in
G_{\mbox{\tiny B}}(n)$ the balls $B(c_1',\p(u_n))$ and
$B(c_2',\p(u_n))$ with $c_i' \in G_{\mbox{\tiny B}}(n,c_i)$ are
disjoint since $B(c_i',\p(u_n)) \subset B(c_i, \mbox{\footnotesize
$ \frac{3}{4} $}\r(u_n))$ and $B(c_1,3\r(u_n))\cap
B(c_2,3\r(u_n))=\emptyset$. Therefore,
$$m(A_n(\p,B))  \asymp  m(B_n(\p(u_n))) \ \# \gb(n,c) \ \#\gb(n)$$
and in view of  (\ref{NGn2}) and (\ref{NGncl}),
\begin{equation}\label{mAncomp2}
m(A_n(\p,B))  \asymp
 m(B) \ \times \  \frac{m(B_n(\p(u_n)))}{m(B_n(\r(u_n)))}
\left(\frac{\r(u_n)}{\p(u_n)}\right)^\gamma \ .
\end{equation}
Finally, let $$ A(\p,B) \, := \, \limsup_{n \to \infty} \,
A_n(\p,B) \, := \, \bigcap_{m=1}^{\infty} \ \bigcup_{n=m}^{\infty}
A_n(\p,B ) \, . $$ By construction and the fact that $\psi$ is
decreasing, we have $ A_n(\p,B)  \subset  \De^u_l(\p,n) \cap B  $
and so  $A(\p,B)$ is a subset of $\La(\p) \cap B$. Now in view of
(\ref{my1}) the proof of Theorem~\ref{THM1} will be completed on
showing that
\begin{equation}\label{my1+}
m (A(\p,B)\cap B) \ \ge \  m(B)/C.
\end{equation}

Notice that estimate (\ref{mAncomp2}) on $m(A_n(\p,B))$ together
with the divergent sum hypothesis (\ref{af2}) of the theorem
implies that
\begin{equation}
\sum_{n = 1}^{\infty} m(A_n(\p,B)) = \infty \ . \label{div}
\end{equation}
This is a good sign as if the above sum was to  converge, then a
simple consequence of the Borel--Cantelli lemma is  that $ m
(A(\p,B)) = 0  \, $. However, the divergent sum alone is not
enough to ensure positive measure; independence of some sort is
also required. The following quasi-independence on average will be
sufficient.

\medskip

\begin{lemma}[Quasi--independence on average]
\label{In2} There exists a constant $C >1$ such that for $Q$
sufficiently large, $$
 \sum_{s, t = 1}^{Q} m (A_s(\p,B) \cap A_t(\p,B) ) \; \leq \; \frac{C}{m(B)} \,
\left( \sum_{s=1}^{Q} m (A_s(\p,B)) \right)^2 \; . $$
\end{lemma}

Clearly, Lemma \ref{In2} together with the divergent sum
(\ref{div}) and Proposition \ref{lem2} implies (\ref{my1+}). This
therefore completes the proof of Theorem \ref{THM1}, assuming of
course the quasi--independence on average result which we now
prove.

\subsection{Proof of Lemma  \ref{In2} : quasi--independence on
average \label{quasi} \!\! }
Throughout, we fix the ball $B$  and write $A_t(\p)$ for $
A_t(\p,B)$.  Also, let $s<t$ and  note that
\begin{eqnarray}
m(A_s(\p)\cap A_t(\p)) & = & m\left(\bigcup_{c\in G_{\mbox{\tiny
B}}(s)} \bigcup_{c'\in G_{\mbox{\tiny B}}(s,c)} B(c',\p(u_s)) \cap
A_t(\p)\right) \nonumber\\ &&\nonumber\\
 & = & \sum_{c\in G_{\mbox{\tiny B}}(s)}
 \sum_{c'\in G_{\mbox{\tiny B}}(s,c)} m\left(\, B(c',\p(u_s))\cap
 A_t(\p)\, \right)\nonumber\\ &&\nonumber\\
  & \asymp & \#G_{\mbox{\tiny B}}(s) \
  \#G_{\mbox{\tiny B}}(s,c) \  m\left(\, B_s(\p(u_s)) \, \cap \,  A_t(\p)\, \right).
  \label{asintat}
  \end{eqnarray}
We now obtain an upper bound for $m(B_s(\p(u_s))\cap A_t(\p))$
where $B_s(\p(u_s))$ is by definition any ball of $A_s(\p)$.
Trivially,
\begin{eqnarray}
m(B_s(\p(u_s))\cap A_t(\p))  & := &
m\Big(B_s(\p(u_s))\cap\bigcup\limits_{c\in G_{\mbox{\tiny
B}}(t)}\bigcup\limits^{\circ}_{c'\in G_{\mbox{\tiny B}}(t,c)} \
B(c',\p(u_t))\Big)  \nonumber \\  \nonumber \\ & \asymp &
   \sum\limits_{c\in G_{\mbox{\tiny B}}(t)}
 \sum\limits_{c'\in G_{\mbox{\tiny B}}(t,c)}
 m\Big(B_s(\p(u_s))\cap B(c',\p(u_t))\Big).
 \label{mbsintat}
\end{eqnarray}
We proceed by considering  two  cases depending on the size of
$\p(u_s)$ compared to $\r(u_t)$.

\noindent{\bf Case (i):} $t > s $ such that $2 \,
\p(u_s)<\r(u_t)$. \ \ \ Suppose that there are two elements $c_1$,
$c_2\in G_{\mbox{\tiny B}}(t)$ such that $$ B_s(\p(u_s))\cap
B(c_i,\r(u_t))\neq \emptyset \hspace{10mm} (i=\{1,2\}) \ . $$
Then, $\dist(c_1,c_2)\le 2\p(u_s)+2\r(u_t)<3\r(u_t)$. However, by
construction the balls $B(c_i,3\r(u_t))$  are disjoint, thus
$\dist(c_1,c_2)\ge 3\r(u_t)$. Hence, there is at most one ball
$B(c,\r(u_t))$ with $c\in G_{\mbox{\tiny B}}(t)$ that can possibly
intersect $B_s(\p(u_s))$. Now, (\ref{2deltapsiB}) together with
the upper bound intersection condition implies that
\begin{eqnarray*}
\sum_{c'\in G_{\mbox{\tiny B}}(t,c)}  m\Big(B_s(\p(u_s)) \cap
B(c',\p(u_t))\Big)  & \asymp &
  m\left(B_s(\p(u_s)) \ \ \cap \textstyle \bigcup\limits^\circ_{c'\in G_{\mbox{\tiny
  B}}(t,c)}\!\!\!\!
B(c',\p(u_t))\right) \\   & \hspace*{-10ex} \ll  & \hspace*{-6ex}
  m\left(B_s(\p(u_s))\cap B(c,\textstyle \frac{3}{4}\r(u_t))\cap
 \De(\ra,\p(u_t))\right)   \\ [1ex] &   \hspace*{-10ex} \ll  & \hspace*{-6ex}
  m(B_t(\p(u_t))) \left(\frac{\p(u_s)}{\p(u_t)}\right)^\gamma.
 \end{eqnarray*}

\noindent In view of  (\ref{mbsintat}) and the fact that at most
one ball $B(c,\r(u_t))$ with $c\in G_{\mbox{\tiny B}}(t)$  can
intersect $B_s(\p(u_s))$,
 we have that $$ m\left(B_s(\p(u_s))\cap
A_t(\p)\right) \ll m\left(B_t(\p(u_t))\right)
\left(\frac{\p(u_s)}{\p(u_t)}\right)^\gamma \ .$$ This together
with (\ref{NGn2}), (\ref{NGncl}) and (\ref{asintat}) implies that
\begin{eqnarray*}
m(A_s(\p)\cap A_t(\p)) \ \ll \ m(B) \ \times \
\frac{m\left(B_t(\p(u_t))\right)}{m\left(B_s(\r(u_s))\right)}
 \ \left(\frac{\r(u_s)}{\p(u_t)}\right)^\gamma .
\end{eqnarray*}

\noindent {\bf Case (ii):} $t > s $ such that $2 \,
\p(u_s)\ge\r(u_t)$. \ \ \ It follows from (\ref{mbsintat}) that $$
\begin{array}{lcl}
m\left(B_s(\p(u_s))\cap A_t(\p)\right) & \ll & \sum^*
\sum\limits_{c'\in G_{\mbox{\tiny
B}}(t,c)}m\left(B(c',\p(u_t))\right)
\\[3ex]
& \ll & m\left(B_t(\p(u_t))\right) \ \#G_{\mbox{\tiny B}}(t,c) \
N(t,s)\; ,
\end{array}
$$ where the sum $\sum^*$ is taken over $c\in G_{\mbox{\tiny
B}}(t)$ such that $B(c,\r(u_t))\cap B_s(\p(u_s))\neq\emptyset$ and
$N(t,s)$ denotes the number of such $c$. Clearly, in case (ii)
$B(c,\r(u_t))\cap B_s(\p(u_s))\neq \emptyset$ implies that
$B(c,\r(u_t))\subset B_s(5\p(u_s))$. Since the balls
$B(c,\r(u_t))$ with $c\in G_{\mbox{\tiny B}}(t)$ are disjoint we
obtain the following trivial estimate

$$ N(t,s) \ \le \
\frac{m\left(5B_s(\p(u_s))\right)}{m\left(B_t(\r(u_t))\right)} \
\asymp \
\frac{m\left(B_s(\p(u_s))\right)}{m\left(B_t(\r(u_t))\right)} \
.$$ Thus $$ m\left(B_s(\p(u_s))\cap A_t(\p)\right)\ll
m(B_t(\p(u_t))) \ \#G_{\mbox{\tiny B}}(t,c) \
\frac{m\left(B_s(\p(u_s))\right)}{m\left(B_t(\r(u_t))\right)} \
,$$

\noindent which together with (\ref{NGn2}), (\ref{NGncl}),
(\ref{mAncomp2}) and (\ref{asintat}) implies that $$ m(A_s(\p)\cap
A_t(\p))  \ll
 \frac{1}{m(B)} \ m(A_s(\p)) \ m(A_t(\p)).
$$

The upshot of these two cases, is that for $Q$ sufficiently large
$$
\begin{array}{l}
\hspace*{-4mm} \displaystyle
 \sum_{s, t = 1}^{Q} m(A_s(\p) \cap A_t(\p) ) \\[2ex] \hspace*{11ex} \displaystyle
 = \
\sum_{s=1}^{Q} m(A_s(\p)) \ + \ 2\,  \sum_{s=1}^{Q-1}
\sum_{\substack{s+1 \leq  t \leq Q  \\ {\rm case (i)}} } m(A_s(\p)
\cap A_t(\p) )\\[2ex] \hspace*{17ex} \displaystyle + \  2 \,
\sum_{s=1}^{Q-1} \sum_{\substack{s+1 \leq t \leq Q \\ {\rm case
(ii)}} } m(A_s(\p) \cap A_t(\p) )
\\[7ex] \hspace*{11ex} \displaystyle
 \ll  \ \sum_{s=1}^{Q}
m(A_s(\p)) \quad  + \ \frac{1}{m(B)}\Big( \sum\limits_{s=1}^{Q}
m(A_s(\p)) \Big)^2 \\[3ex] \hspace*{17ex} \displaystyle + \ \
 m(B)\ \sum\limits_{s=1}^{Q-1}  \sum_{\substack{s+1 \leq t \leq Q \\
\p(u_s)<\r(u_t) } }
\frac{m\left(B_t(\p(u_t))\right)}{m\left(B_s(\r(u_s))\right)}
   \left(\frac{\r(u_s)}{\p(u_t)}\right)^\gamma \!\! .
\end{array}
$$ By (\ref{mAncomp2}) and condition (\ref{af}) imposed in the
statement of the theorem, the latter double sum is $\ll \left(
m(B)^{-1} \sum_{s=1}^{Q} m(A_s(\p)) \right)^2 $. By (\ref{div}),
for $Q$ sufficiently large $\sum_{s=1}^{Q} m(A_s(\p))\le
m(B)^{-1}(\sum_{s=1}^{Q} m(A_s(\p)))^2$. The statement of Lemma
\ref{In2} now readily follows and thereby completes the proof of
Theorem \ref{THM1}. \QED

\section{Proof  of Theorem \ref{THM3}: $ 0\leq G < \infty$
\label{proof3}}

We begin by observing that the case  $\gamma=\delta$ is excluded
by the various hypotheses imposed on the dimension function $f$ --
see also \S\ref{r3}. {\em Thus, without loss of generality we can
assume in proving  Theorem \ref{THM3} that} $$0 \ \leq \ \gamma \
< \ \delta \, .$$

To prove Theorem \ref{THM3} we proceed as follows. For any fixed
$\eta \gg 1$ we construct a Cantor subset ${\bf K_{\eta}} $ of
$\Lambda(\psi)$ and a probability measure $\mu$ supported on ${\bf
K_{\eta}}$ satisfying the condition that for an arbitrary ball $A$
of sufficiently small radius $r(A) $
\begin{equation}
\mu(A) \; \ll \; \frac{f(r(A))}{\eta}  \;   , \label{task}
\end{equation}
where the implied constant is absolute. By the Mass Distribution
Principle, the above inequality implies that
$
\mathcal{ H}^f({\bf K_{\eta}}) \; \gg \; \eta \;.
$
Since $ {\bf K_{\eta}}\subset \Lambda(\psi)$, we obtain that
$\mathcal{ H}^f \left( \Lambda(\psi) \right)\gg \eta$. However, $\eta
\gg 1  $ is arbitrarily large whence $\mathcal{ H}^f   \left(
\Lambda(\psi) \right)=\infty$ and this proves  Theorem \ref{THM3}.

In view of the above outline, the whole strategy of our proof is
centred around the  construction of a `right type' of Cantor set
${\bf K_{\eta}}$ which  supports a measure $\mu$ with the desired
property. The actual nature of the construction of ${\bf
K_{\eta}}$ will depend heavily on whether  $G$ defined by
(\ref{my5}) is finite  or infinite. In this section we deal with
the case that $0\leq G < \infty$. The case that $G = \infty$ is
substantially easier -- see \S\ref{pfdimthm}.

\subsection{Preliminaries}
In this section we group together for clarity and convenience
various concepts and results which will be required in
constructing the Cantor set ${\bf K_{\eta}}$. We shall make use of
the various  hypotheses of Theorem  \ref{THM3} as required. In
particular, the measure $m$ is of type (M2).

\subsubsection{The sets $\gb(n)$ and
$\gb(n,c)$.}\label{pregb}
Let $B=B(x,r)$ be an arbitrary ball with centre $x \in \Omega$.
Assume that its radius $r$ is sufficiently small so that local
$m$-ubiquity and the measure estimate (\ref{measure}) are
fulfilled for $B$.  Relabel the sets $\gb(n)$ and $\tilde\go(n)$
constructed in \S\ref{proof1}  by $\gb'(n)$ and $\go'(n)$
respectively. By keeping track of constants, the estimate
(\ref{NGn2}) for $\#\gb'(n)$ is explicitly as follows: $$
 \frac{a \, \ka}{b \, (36)^\delta} \
\left(\frac{r}{\r(u_n)} \right)^\delta \ \leq \  \#\gb'(n) \ \leq
\ \frac{b }{a } \ \left(\frac{r}{\r(u_n)} \right)^\delta \ , $$
where $a, b $ and $\delta$ are as in (\ref{measure}) and $\ka$ is
as in (\ref{ub}). Since `local' implies `global' ubiquity,  the
corresponding estimates for $\#\go'(n)$ are explicitly as follows:
$$
 \frac{ \ka_1}{b \, 9^\delta} \
\left(\frac{1}{\r(u_n)} \right)^\delta \ \leq \ \#\go'(n) \ \leq \
\frac{1}{a} \ \left(\frac{1}{\r(u_n)} \right)^\delta \  $$ where
$0 < \ka_1 \leq \ka$ is the global ubiquity constant arising from
local ubiquity.

Now let $ 0 \, < \, c_3:= \min\{ \frac{a \, \ka}{b \,
(36)^\delta}, \frac{ \ka_{\mbox{\tiny 1}}     }{b \, 9^\delta} \}
\, < 1 $
 and  define $\gb(n)$ to be any
sub-collection of $ \gb'(n) $ such that $$ \#\gb(n) \ = \left[c_3
\left(\frac{r}{\r(u_n)} \right)^\delta \right] \ ,  $$ where $[x]$
denotes the integer part of a real number $x$.
 Thus, for $n$ sufficiently large
\begin{equation}
 \mbox{\small $ \frac{1}{2} $ }  \, c_3 \
\left(\frac{r}{\r(u_n)} \right)^\delta \ \leq \
 \#\gb(n) \ \leq  \ c_3 \left(\frac{r}{\r(u_n)} \right)^\delta
 \ ,  \label{NGnB}
\end{equation}

\noindent where we take $r=1$ when $B$ is replaced by $\Omega$.

\vspace{2mm}


In the following the arbitrary ball $B$ can be replaced by the
whole space $\Omega$ without any  loss of generality. Assume for
the moment that $\gamma > 0 $.  Associated with  any $ c\in
\gb(n)$  is the set $\gb(n,c) $ for which (\ref{deltapsiB}) and
(\ref{2deltapsiB}) are satisfied. Then the estimate (\ref{NGncl})
for $\#\gb(n,c)$ is explicitly given by
\begin{equation}
 c_4  \ \left( \frac{\r(u_n)}{\p(u_n)} \right)^\gamma
\, \leq \,  \# \gb(n,c) \,\leq \, 
\ c_5 \ \left( \frac{\r(u_n)}{\p(u_n)} \right)^\gamma \ .
\label{gbncl}
\end{equation}
\noindent where $0 \, < \,  c_4 := \frac{c_1 \, a}{2^\gamma
9^\delta b} \, < \, 1 \, < \, c_5 := \frac{c_2 \, b}{a}  $ and
$c_1$ and $c_2$ are the constants appearing in the intersection
conditions. In the case $\ga = 0$, we define $\gb(n,c):= \{c\}$
and so $\# \gb(n,c) = 1$.

\noindent Note that for any distinct $c',c'' \in \gb(n,c)$  we
have that $d(c',c'') \geq 3 \p(u_n)$. This follows from the fact
that by construction the respective balls of radius $3 \p(u_n)$
are disjoint. Hence, for any $x \in B(c',\p(u_n))$ and $y \in
B(c'',\p(u_n))$ we have that
\begin{equation}
 d(x,y) \ \geq \ 2 \, \p(u_n)   \ \label{dxy} .
\end{equation}
Also recall that any ball $B(c',\p(u_n))$ with $c' \in \gb(n,c)$
is contained in $B(c,\mbox{\footnotesize $ \frac{3}{4}$}\r(u_n))$
and in turn the ball $B(c,\r(u_n))$ is contained in $B$.

\smallskip

\noindent{\em Remark: \ } In the construction of  the set
$\gb(n,c)$ we make use of the fact that without loss of
generality, $\r(u_n)
> 24 \p(u_n)$ for $n$ sufficiently large.
This guarantees  that any ball $B(c',3\p(u_n))$ with $c' \in
\gb(n,c)$ is contained in $B(c,\mbox{\footnotesize $
\frac{3}{4}$}\r(u_n))$.  To see that the above  fact remains valid
under the hypotheses of Theorem \ref{THM3} we  observe that,
without loss of generality, we can assume  that $\rho(u_n)^{-1} \,
\psi(u_n) \to 0$ as $n \to \infty $. If this was not the case then
$\limsup \rho(u_n)^{-1} \, \psi(u_n) > 0$ as $n \to \infty $ and
since `local' implies `global' ubquity,  Theorem \ref{THM1}
implies that $m(\Lambda(\p))
> 0 $. In turn this implies that $\mathcal{ H}^{\delta} (\Lambda(\p))
$ is positive and finite. By the elementary fact stated in
\S\ref{hmd}, $\mathcal{ H}^f (\Lambda(\p)) = \infty $ for any
dimension function $f$ such that $r^{-\delta} \, f(r) \to \infty$
as $r \to 0$.

\subsubsection{Working on a  subsequence of
$u$.}\label{subseqwork}
To begin with recall the following simple  facts: (i) if the
hypothesis of `global' or `local' ubiquity are satisfied for a
particular upper sequence  $u$ then they are also satisfied for
any subsequence $s$ and (ii) if $\r$ is $u$-regular then it is
$s$-regular for any subsequence $s$. Also note that if $G$ is
finite, then $\limsup_{n \to \infty}g(s_n) < \infty $ for any
subsequence $s$ of $u$.

 Now notice that for any $m \in \N $, if   $\r$ is
$u$-regular with constant $\lambda < 1 $ then  we can find  a
subsequence $s$ of $u$ such that $$ \r(s_{t+1})  \ <  \ \lambda^m
\; \r(s_{t}) \hspace{7mm} \mbox{ and such that } \hspace{7mm} \sum
g(s_t) = \infty \ .  $$

\noindent Thus, without loss of generality  in establishing part
(i) of Theorem \ref{THM3} ($G=0$) we can assume that $\r$ is
$u$-regular with constant $\lambda$ as small as we please. The
existence of such a subsequence $s$ is easy to verify. Trivially,
for $p \geq m $ we have that $ \r(u_{n+p})  \ < \ \lambda^m \;
\r(u_{n})$. For $t \in \N$,  let $ g(u_{r_t}) := \max \{ g(u_r) :
m(t-1) < r \leq mt \} $. Then
\begin{eqnarray*}
\infty  \; = \; \sum_{r=1}^{\infty} g(u_r) \ & = & \
\sum_{t=1}^{\infty}  \ \ \ \ \sum_{m(t-1) < r \leq  mt} g(u_r) \
\leq \ m \, \sum_{t=1}^{\infty} g(u_{r_t} ) \\ & \ll & \ \,
\sum_{n=1}^{\infty} g(u_{r_{2n}}) \ + \ \  \sum_{n=1}^{\infty}
g(u_{r_{2n-1}}) \ \  .
\end{eqnarray*}
Thus on  both the sequences $\{u_{r_{2n}}\} $ and $ \{
u_{r_{2n-1}} \} $ the function $\rho$ satisfies the required
regularity condition and for one of them the divergent  sum
condition is satisfied.

Next notice that if $0 < G < \infty $, then there exists  a
strictly increasing  sequence $\{ n_i\} $ such that  $g(u_{n_i})
\geq G/2
> 0 $. Since $\lim_{r \to \infty}  \r(r) = 0$,  it follows
that for any $\lambda < 1 $ there exists a subsequence $s$  of $\{
u_{n_i} \} $ such that $\r(s_{t+1})   <   \lambda \, \r(s_{t}) $
and $  \sum g(s_t) = \infty $. Thus,  without loss of generality
in establishing Theorem \ref{THM3} for the case that $0 \leq G <
\infty $ we can assume that $\r$ is $u$ regular with constant as
small as we please.


\subsection{The Cantor set ${\bf K_{\eta}}$ }
We are  assuming  that $0 \leq G < \infty$. Let $G^*:= \max \{2,
\sup_{n \in \N} g(u_n) \} $. Then $$ g(u_n) \ < \ G^*
\hspace{8mm}  {\rm for \ all \ } n  \ \ . $$ Now fix a real number
$\eta$ such that $$\eta \
> \ G^* \ \ . $$
To avoid cumbersome expressions, let $\varpi$ denote the following
repeatedly occurring constant
\begin{equation} \varpi \ := \
\frac{c_3 c_4  a }{3^\delta  \, 32 \, \,   b^2 \, c_2  } \  \ < \
1 \ \  . \label{varpi}
\end{equation}
In view of \S\ref{subseqwork}, we can assume  that for $n$
sufficiently large $\r(u_{n+1}) \leq \lambda  \r(u_n) $ with
\begin{equation}
0 \ < \ \lambda \ <  \  \left( \frac{ a }{a + 3^\delta  \, 8 \, \,
b \, c_2  } \right)^{\frac{1}{\delta - \gamma}}  \ .
\label{lambda}
\end{equation}
Finally, unless stated otherwise  $B(r)$ will denote a generic
ball of radius $r$ centred at a point in $\Omega$.

\subsubsection{Constructing the first level ${\bf K(1)}$.}
Choose $t_1$ large enough so that
\begin{equation}
g(\kt) \  <  \ G^* \  <  \   \frac{a}{ 3^{\delta}\, 8 \, c_2 \, b}
\ \frac{\eta}{\varpi}  \ ,  \label{t1choiceb}
\end{equation}
\begin{equation} \frac{\fp\kt))}{\p(\kt)^\delta} \ > \ 3^{\delta -
\ga} \, \frac{\eta}{\varpi}  \ ,    \label{t1choiced}
\end{equation}
and so that the counting estimate (\ref{NGnB}) is valid for the
set $G_{\Omega}(t_1)$. Note that the first of these inequalities
is possible since   $g(u_n) < G^* $ for all n  and that
 $\eta > G^*$. The latter inequality is
possible  since $f(r)/r^\delta\to\infty$ as $r\to0$. Let $k_1 \geq
1 $ be the unique integer  such that
\begin{eqnarray}
\frac{ 3^{\delta}\, 2 \, c_2 \, b}{a} \ \frac{\varpi}{\eta} \
\sum_{i=0}^{k_1-1} g(u_{t_1+i}) &  \le &
\frac{1}{4}\label{k1choicea}\\ \nonumber \\  \frac{ 3^{\delta}\, 2
\, c_2 \, b}{a} \ \frac{\varpi}{\eta} \  \sum_{i=0}^{k_1}
g(u_{t_1+i})& > & \frac{1}{4}\label{k1choiceb}
\end{eqnarray}
Note, the fact that $k_1 \geq 1 $ is a consequence of
(\ref{t1choiceb}).

\bigskip

 The first level ${\bf K(1)}$ of the Cantor set ${\bf K_{\eta}}$ will
 now be constructed with the above $\eta$ in mind. This level will
 consist of sub-levels $K(t_1+ i)$  where $0 \, \leq \, i \, \leq k_1$.
\vspace{2ex}

\noindent $\bullet$ \ {\bf The sub-level $K(t_1)$ : \ } This
consists of balls of common radius $\p(u_{t_1})$ defined as
follows:- $$K(t_1)\ := \ \bigcup_{c\in
G_{\Omega}(t_1)}\bigcup_{c'\in G_{\Omega}(t_1,c)}^{\circ}
B(c',\p(\kt)) \ . $$

\vspace{2ex}

\noindent $\bullet$ \ {\bf The sub-level $K(t_1+1)$ : \ } The
second sub-level will consist of balls of common radius
$\p(u_{t_1+1})$  which substantially avoid balls from the previous
sub-level $K(t_1)$. Let $$ h(t_1) \ := \ \left(\frac{\varpi}{\eta}
\ \frac{\fp\kt))}{ \p(\kt)^\ga}\right)^{1/(\delta - \ga)} \ . $$

\noindent Consider some point $c\in G_{\Omega}(t_1)$. Thus, $c$
lies on a resonant set $\ra$ with $\alpha\in J^u_l(t_1)$.
Construct the `thickening' $$T_c(t_1) \ := \ \De\left(\ra,h(t_1)
\right)  \, \cap \, B(c, \r(\kt)) \ . $$ Note that in view of
(\ref{t1choiced}) we have $3 \p(u_{t_1})  <  h(t_1)$,  and so by
(\ref{2deltapsiB}) $$
 \bigcup_{c'\in G_{\Omega}(t_1,c)}^\circ B(c',\p(u_{t_1})) \
\subset \ T_c(t_1) \ . $$ Also, notice that in view of
(\ref{t1choiceb}) we have that $h(t_1) < \r(u_{t_1}) \, $. Now
define $$T(t_1) \ := \ \left\{ T_c(t_1) \, : \,  c\in
G_{\Omega}(t_1) \right\} \ , $$ thus  $ \#T(t_1) =
\#G_{\Omega}(t_1)$. Clearly, the `thickenings' in the collection $
T(t_1)$ are disjoint since the balls $ B(c, 3\r(\kt))$ with $ c\in
G_{\Omega}(t_1)$ are disjoint. We now introduce a collection of
balls from which the next sub-level is to be constructed. Consider
the set $G_{\Omega}(t_1+1)$ and for each $c \in G_{\Omega}(t_1+1)$
construct the ball $B(c,\r(\ktt))$. Clearly these balls are
disjoint. We disregard any of these balls which lie too close to
balls from the previous sub-level. To make this precise, we
introduce the sets
\be U_{\Omega}(t_1+1)&:=&\{c\in
G_{\Omega}(t_1+1):B(c,\r(\ktt))\cap T(t_1)\neq\emptyset\}\\[2ex]
V_{\Omega}(t_1+1)&:=&G_{\Omega}(t_1+1)\setminus U_{\Omega}(t_1+1)
\ . \ee We claim that $\#V_{\Omega}(t_1+1) \geq \frac{1}{2} \,
\#G_{\Omega}(t_1+1)  \ . $ This will obviously follow on
establishing the upper bound $$\#U_{\Omega}(t_1+1) < \frac{1}{2}
\, \#G_{\Omega}(t_1+1)  \ . $$


\noindent There are two cases to consider.

\vspace{2mm} \noindent{\rm \underline{Case} (i): \, $\r(\ktt) <
h(t_1)$. \, } Suppose that $B(c_1,\r(\ktt))$ with $c_1\in
G_{\Omega}(t_1+1)$ intersects  some $T_c(t_1)\in T(t_1)$. Then $
B(c_1,\r(\ktt)) \cap B(c,\r(\kt))\neq\emptyset $, and since
$\r(\ktt) \leq \r(\kt)$ we have the inclusion  $B(c_1,\r(\ktt))
\subset B(c,3\r(\kt)) $. Moreover, for $\r(\ktt) < h(t_1)$ $$
B(c_1,\r(\ktt)) \subset \De(\ra,3h(t_1))\cap  B(c,3\r(\kt)) \ . $$
 Let $N$ denote the number
of   balls $B(c_1,\r(\ktt))$ with $c_1\in G_{\Omega}(t_1+1)$ that
can possibly intersect some fixed $T_c(t_1)\in T(t_1)$. Then, $$
\begin{array}{rcl}
m\left(\De(\ra,3h(t_1))\cap  B(c,3\r(\kt))\right) & \geq & N \,
m(B(\r(\ktt))) \\[1ex]
 & \geq &  N \, a \ \r(\ktt)^{\delta} .
\end{array}
 $$
 For $h(t_1) < \r(\kt)$, the intersection condition (ii) implies
that $$ m\left(\De(\ra,3h(t_1))\cap  B(c,3\r(\kt))\right) \, \leq
\, c_2 \, b \, 3^{\delta} \, h(t_1)^{\delta-\ga} \, \r(\kt)^{\ga}
\ . $$ Thus, $$ N \ \leq \ \frac{c_2 \, b \, 3^{\delta}}{a} \
\frac{\varpi}{\eta} \ \frac{f(\p(\kt))}{\p(\kt)^{\ga}} \
\frac{\r(\kt)^{\ga}}{\r(\ktt)^{\delta}} \ \ . $$

\noindent It follows, by (\ref{NGnB}) and (\ref{t1choiceb}) or
equivalently (\ref{k1choicea}), that
\begin{eqnarray*}
\#U_{\Omega}(t_1+1)& \le &  N \ \#T(t_1) \ \leq \ \frac{c_2 \, b
\, 3^{\delta}}{a} \ \frac{\varpi}{\eta} \
\frac{f(\p(\kt))}{\p(\kt)^{\ga}} \ \frac{ c_3 \, \r(\kt)^{\ga -
\delta }}{\r(\ktt)^{\delta}} \\ & & \\
 & \le & \frac{2 \, c_2 \, b \, 3^{\delta}}{a} \
\frac{\varpi}{\eta} \ \ g(\kt) \ \#G_{\Omega}(t_1+1) \ < \
 \frac{1}{4}\#G_{\Omega}(t_1+1) \ .
\end{eqnarray*}

\vspace{2ex} \noindent{\rm \underline{Case} (ii): \, $\r(\ktt)
\geq h(t_1)$. \, } A similar argument to that given above implies
that if $B(c_1,\r(\ktt))$ with $c_1\in G_{\Omega}(t_1+1)$
intersects some $T_c(t_1)\in T(t_1)$ then
$$ B(c_1,\r(\ktt)) \subset
\De(\ra,3\r(\ktt))\cap B(c,3\r(\kt)) \ . $$

\noindent As before,  let $N$ denote the number of   balls
$B(c_1,\r(\ktt))$ with $c_1\in G_{\Omega}(t_1+1)$ that can
possibly intersect some fixed $T_c(t_1)\in T(t_1)$. Then, $$
\begin{array}{rcl}
m\left(\De(\ra,3\r(\ktt))\cap B(c,3\r(\kt))\right) & \geq & N \,
m(B(\r(\ktt))) \\[1ex]
 & \geq &  N \, a \ \r(\ktt)^{\delta}  .
\end{array}
$$ Since $\r(\ktt) <  \r(\kt)$, the intersection condition (ii) implies
that $$ m\left(\De(\ra,3\r(\ktt))\cap B(c,3\r(\kt))\right) \, \leq
\, c_2 \, b \, 3^{\delta} \, \r(\ktt)^{\delta-\ga} \,
\r(\kt)^{\ga} \ . $$ Thus, $$ N \ \leq \ \frac{c_2 \, b \,
3^{\delta}}{a} \ \left(\frac{\r(\kt)}{\r(\ktt)}\right)^{\ga} \ \ .
$$


\noindent It follows, by (\ref{NGnB}) and (\ref{lambda}), that
\begin{eqnarray*}
\#U_{\Omega}(t_1+1)& \le &  N \ \#T(t_1) \ \leq \ \frac{c_2 \, b
\, 3^{\delta}}{a} \ \left(\frac{\r(\kt)}{\r(\ktt)}\right)^{\ga} \
c_3 \ \left( \frac{1}{\r(\kt)}\right)^\delta
\\ & & \\
 & \le & \frac{2 \, c_2 \, b \, 3^{\delta}}{a} \
\left(\frac{\r(\ktt)}{\r(\kt)}\right)^{\delta - \ga} \
\#G_{\Omega}(t_1+1)
\\ & & \\
 & \le & \frac{2 \, c_2 \, b \, 3^{\delta}}{a} \
\lambda^{\delta - \ga} \ \#G_{\Omega}(t_1+1)
 \ < \
 \frac{1}{4}\#G_{\Omega}(t_1+1) \ .
\end{eqnarray*}
On combining the two cases, we have $ \#U_{\Omega}(t_1+1) <
\frac{1}{2} \ \#G_{\Omega}(t_1+1) $. Hence $$ \#V_{\Omega}(t_1+1)
\ \geq \  \frac{1}{2} \ \#G_{\Omega}(t_1+1) \ \ . $$ The second
sub-level is now defined to be: $$K(t_1+1)\ := \ \bigcup_{c\in
V_{\Omega}(t_1+1)}\bigcup_{c'\in G_{\Omega}(t_1+1,c)}^{\circ}
 B(c',\p(\ktt)) \ .
 $$
Note, by construction $K(t_1) \ \cap \ K(t_1+1) \ = \ \emptyset .
$

\medskip

\noindent $\bullet$ \ {\bf The sub-level $K(t_1+i)$  : \ } Suppose
$k_1 \geq 2$. Fix   $2 \leq i \leq k_1$ and for $1 \leq j \leq
i-1$ suppose we have constructed the sub-levels $$K(t_1+j) \
=\bigcup_{c\in V_{\Omega}(t_1+j)} \ \ \bigcup_{c'\in
G_{\Omega}(t_1+j, \, c)}^{\circ}
 B(c',\p(u_{t_1+j})) \ . $$
We proceed to construct $K(t_1+i)$. For a point $ c \in
V_{\Omega}(t_1+(i-1)) $ there exists a resonant set $\ra$ with
$\alpha \in J^u_l(t_1 + (i-1)) $ such that $c \in \ra$. Let $$
h(t_1+(i-1)) \ := \ \left(\frac{\varpi}{\eta}  \ \frac{
f(\p(u_{t_1+(i-1)}))    }{
\p(u_{t_1+(i-1)})^\ga}\right)^{1/(\delta - \ga)} \  , $$ and
construct the `thickening' $$T_c(t_1+(i-1)) \ := \
\De\left(\ra,h(t_1+(i-1)) \right)  \, \cap \, B(c,
\r(u_{t_1+(i-1)})) \ . $$ Note that in view of (\ref{t1choiced})
and the fact that $f(r)/r^{\delta}$ is decreasing $$3
\p(u_{t_1+(i-1)}) \ < \ h(t_1+(i-1)) \ , $$ and so by
(\ref{2deltapsiB}) $$ \bigcup_{c'\in G_{\Omega}(t_1+(i-1), \,
c)}^\circ B(c',\p(u_{t_1+(i-1)})) \ \subset \ T_c(t_1+(i-1)) \ .
$$ Also, notice that in view of (\ref{t1choiceb}) and the fact
that $g(u_n) < G^*$ for all $n$  we have that $h(t_1+(i-1)) <
\r(u_{t_1+(i-1)}) \, $. Define $$T(t_1+(i-1)) \ := \ \left\{
T_c(t_1+(i-1)) \, : \,  c\in V_{\Omega}(t_1+(i-1)) \right\} \ . $$
Thus,  $ \#T(t_1+(i-1))  =   \# V_{\Omega}(t_1+(i-1))  \leq
\#G_{\Omega}(t_1+(i-1)) $. Clearly  the above `thickenings' in $
T(t_1+(i-1))$ are disjoint since the balls $ B(c,
3\r(u_{t_1+(i-1)}))$ with $ c\in G_{\Omega}(t_1+(i-1))$ are
disjoint.  Now introduce  the set $G_{\Omega}(t_1+i)$ and for each
$c \in G_{\Omega}(t_1+i)$ construct the ball $B(c,\r(u_{t_1+i}))$.
Obviously  these balls are disjoint and we proceed by disregarding
any of  those  which lie too close to balls from  any of the
previous sub-levels $K(t_1+j)$. To make this precise,  introduce
$$
\begin{array}{l}
\displaystyle
U_{\Omega}(t_1+i):=\{c\in
G_{\Omega}(t_1+i):B(c,\r(\ktk))\cap \bigcup_{j=0}^{i-1}T(t_1+j)
\neq\emptyset\} \\ \displaystyle
V_{\Omega}(t_1+i):=G_{\Omega}(t_1+i)\setminus U_{\Omega}(t_1+i) \
.
\end{array}
$$
\noindent We claim that $\#V_{\Omega}(t_1+i) \geq \frac{1}{2} \,
G_{\Omega}(t_1+i)  \ . $ This will obviously follow on
establishing the upper bound $\#U_{\Omega}(t_1+i) < \frac{1}{2} \,
G_{\Omega}(t_1+i) $. As before there are two cases to consider.

\vspace{3mm} \noindent{\rm \underline{Case} (i):  $~ \;  \ 0 \leq
j \leq i-1$ such that $ \r(u_{t_1+i}) < h(t_1+j)$. \, }

\vspace{2mm} \noindent{\rm \underline{Case} (ii): \, $ 0 \leq j
\leq i-1$ such that $\r(u_{t_1+i}) \geq h(t_1+j)$. \, }

\vspace{3mm}

\noindent On following the arguments as in the $i=1$ case, we
obtain that
\begin{eqnarray*}
\#U_{\Omega}(t_1+i) & \le &  \sum_{\rm{case\ (i)\ }} \frac{c_2 \,
b \, 3^{\delta} }{a} \ \frac{\varpi}{\eta} \
\frac{f(\p(u_{t_1+j})) }{ \p(u_{t_1+j})^{\ga}  } \
\frac{\r(u_{t_1+j})^{\ga} }{ \r(u_{t_1+i})^{\delta} } \ \#T(t_1+j)
\\ & & + \sum_{\rm{case\ (ii)\ }} \frac{c_2 \, b \,
3^{\delta}}{a} \
\left(\frac{\r(u_{t_1+j})}{\r(u_{t_1+i})}\right)^{\ga}  \,
\#T(t_1+j).
\end{eqnarray*}
The contribution from  case (i) is:
\begin{eqnarray*}
& \leq  & \sum_{\rm{case\ (i)\ }} \frac{c_2 \, b \, 3^{\delta}}{a}
\ \frac{\varpi}{\eta} \
\frac{f(\p(u_{t_1+j}))}{\p(u_{t_1+j})^{\ga} } \
\frac{\r(u_{t_1+j})^{\ga} }{ \r(u_{t_1+i})^{\delta} } \
\#G_{\Omega}(t_1+j)  \\ \\  & \leq & \sum_{j=0}^{k_1-1} \frac{2 \,
c_2 \, b \, 3^{\delta}}{a} \ \frac{\varpi}{\eta} \ g(u_{t_1+j}) \
\#G_{\Omega}(t_1+i) \ \le \ \frac{1}{4} \ \#G_{\Omega}(t_1+i) \ ,
\end{eqnarray*}  by the choice of $k_1$ -- see (\ref{k1choicea}).
The contribution from case (ii) is:
\begin{eqnarray*}
& \leq  &  \sum_{\rm{case\ (ii)\ }} \frac{c_2 \, b \,
3^{\delta}}{a} \
\left(\frac{\r(u_{t_1+j})}{\r(u_{t_1+i})}\right)^{\ga}  \,
\#G_{\Omega}(t_1+j)
\\ \\  & \leq &
\frac{2 \, c_2 \, b \, 3^{\delta}}{a} \  \sum_{j=0}^{i-1}
\left(\frac{\r(u_{t_1+i})}{\r(u_{t_1+j})}\right)^{\delta - \ga} \
\#G_{\Omega}(t_1+i)
\\ \\  & \leq &
\frac{2 \, c_2 \, b \, 3^{\delta}}{a} \  \sum_{j=0}^{i-1}
\lambda^{(i-j)(\delta - \ga)} \ \#G_{\Omega}(t_1+i)
\\ \\  & \leq &
\frac{2 \, c_2 \, b \, 3^{\delta}}{a} \  \sum_{s=1}^{\infty}
(\lambda^{\delta - \ga})^s \ \#G_{\Omega}(t_1+i) \ < \ \frac{1}{4}
\  \#G_{\Omega}(t_1+i) \ ,
\end{eqnarray*}
by the choice of $\lambda$ -- see (\ref{lambda}).  On combining
the two cases, we obtain that $ \#U_{\Omega}(t_1+i) < \frac{1}{2}
\ \#G_{\Omega}(t_1+i) $ and so  $$ \#V_{\Omega}(t_1+i) \ \geq \
\frac{1}{2} \ \#G_{\Omega}(t_1+i) \ $$ as claimed.  The
 sub-level $K(t_1+i)$ is defined to be:
$$K(t_1+i) \ := \ \bigcup_{c\in V_{\Omega}(t_1+i)} \
\bigcup_{c'\in G_{\Omega}(t_1+1,c)}^{\circ}
 B(c',\p(\ktk)) \ . $$
Also, note that by construction for $ 0 \leq i \neq j \leq k_1 $
$$K(t_1+i) \ \cap \ K(t_1+j) \ = \ \emptyset \ . $$

The first level ${\bf K(1)}$ of the Cantor set is defined to be $$
{\bf K(1)} \ := \ \bigcup_{i=0}^{k_1}K(t_1+i)  \ . $$

\subsubsection{Constructing the second level {\bf K(2)}.}
The second level of the Cantor set is constructed by `looking'
locally at each ball from the previous level. Thus the second
level {\bf K(2)} will be defined in terms of local levels ${\rm
K(2,B)}$ associated with $B$ in {\bf K(1)}.

Choose $t_2 > t_1 $ sufficiently large so that for any ball $B$ in
{\bf K(1)} the counting estimate (\ref{NGnB}) is valid and so that
\begin{equation}
\frac{\fp\kd))}{\p(\kd)^\delta } \ > \ 3^{\delta - \ga} \,
\frac{1}{\varpi} \ \frac{f (r(B))}{m(B) }  \ .
 \label{t2choice}
\end{equation}
Here and throughout, $r(B)$ denotes the radius of the ball $B$.
Note that for  $B \in  {\bf K(1)} $, $$ \p(u_{t_1+k_1}) \leq r(B)
\leq \p(u_{t_1}) \ . $$ In view of (\ref{t1choiced}), the fact
that $g(u_n) < G^* $ for all $n$ and that $f(r)/r^{\delta}$ is
decreasing as $r$ increases, it is easily verified that
\begin{equation}
g(u_{t_2}) \ < \ G^* \ <  \ \frac{a}{ 3^{\delta}\, 8 \, c_2 \, b}
\ \frac{1}{\varpi} \ \frac{f (r(B))}{m(B) } \  \hspace{10mm}
\forall \ \ \ \ B \in {\bf K(1)} . \label{t2gchoice}
\end{equation}

\vspace{7mm}

\noindent $\bullet$ \ {\bf The local sub-level $K(t_2,B)$ : \ }
Fix a ball $B$  in  ${\bf K(1)}$. Thus $B=B(c',\p(u_{t_1+i})) $ is
a ball in the sub-level $K(t_1+i)$ for some $0\leq i \leq k_1$.
Consider the set $G_B(t_2)$. Each $c \in G_B(t_2)$ gives rise to
the set $G_B(t_2,c)$. Let $$ K(t_2,B) \ := \ \bigcup_{c\in
\gb(t_2)}\bigcup_{c'\in \gb(t_2, \, c)}^{\circ} B(c',\p(u_{t_2}))
\ . $$ By construction, $ K(t_2,B) \subset B $ and indeed $
K(t_2,B) $ defines the first local sub-level associated with $B$.
We now proceed  to construct further local sub-levels $
K(t_2+i,B)$ where $1 \leq i \leq k_2(B)$ and $k_2(B) $ is the
unique integer such that
\begin{eqnarray}
\frac{ 3^{\delta}\, 2 \, c_2 \, b \, \varpi }{a} \
\frac{m(B)}{f(r(B))}
 \ \sum_{i=0}^{k_2(B)-1} g(u_{t_2+i}) &  \le &
\frac{1}{4}\label{k2choicea}\\ \nonumber \\  \frac{ 3^{\delta}\, 2
\, c_2 \, b \, \varpi }{a} \ \frac{m(B)}{f(r(B))}
 \ \sum_{i=0}^{k_2(B)} g(u_{t_2+i})& > & \frac{1}{4}\label{k2choiceb}
\end{eqnarray}

\vspace{7mm}

\noindent $\bullet$ \ {\bf The local sub-level $K(t_2+1,B)$ : \ }
Consider some point $c \in \gb(t_2)$. Thus $c$ lies on the set
$\ra \cap \mbox{\small $\frac{1}{2}$} B $ for some $\alpha \in
J^u_l(t_2)$. Construct the thickening

$$T_c(t_2,B) \ := \ \De\left(\ra,h_B(t_2) \right)  \, \cap \, B(c,
\r(\kd)) \  , $$ where $$
 h_B(t_2) \ := \ \left(\frac{\varpi \, m(B)}{f(r(B))} \ \frac{\fp\kd))}{
\p(\kd)^\ga}\right)^{1/(\delta - \ga)} \ . $$

\noindent In view of (\ref{t2choice}), $3 \p(u_{t_2})  < h_B(t_2)
\ , $ and so by (\ref{2deltapsiB})

$$ \bigcup_{c'\in G_{B}(t_2,c)}^\circ B(c',\p(u_{t_2})) \ \subset
\ T_c(t_2,B) \ . $$

\noindent Also, notice that in view of (\ref{t2gchoice}) we have
that $h_B(t_2) < \r(u_{t_2}) \, $. Now define $$T(t_2,B) \ := \
\left\{ T_c(t_2,B) \, : \,  c\in G_{B}(t_2) \right\} \ , $$ thus $
\#T(t_2,B) = \#G_{B}(t_2)$. Moreover, $ K(t_2,B)  \subset T(t_2,B)
 $.  Consider the set $G_{B}(t_2+1)$ and for each $c \in
G_{B}(t_2+1)$ construct the ball $B(c,\r(\kdd))$. Clearly these
balls are disjoint. Introduce the sets \be
U_{B}(t_2+1)&:=&\left\{c\in G_{B}(t_2+1):B(c,\r(\kdd))\cap
T(t_2,B)\neq\emptyset \right\}\\[1ex]
V_{B}(t_2+1)&:=&G_{B}(t_2+1)\setminus U_{B}(t_2+1) \  \ \ . \ee

\noindent We show that  $\#U_{B}(t_2+1) < \frac{1}{2} \,
G_{B}(t_2+1)$ by considering the following  two cases:

\medskip

\noindent{\rm \underline{Case} (i):  $~ \ \, \r(\kdd) < h_B(t_2)$.
\, }

\vspace{2mm}

 \noindent{\rm \underline{Case} (ii): \, $ \r(\kdd)
\geq h_B(t_2)$. \, }

\medskip

\noindent As in the construction of the sub-level $K(t_1+1)$, we
find that for case (i):
\begin{eqnarray*}
\#U_{B}(t_2+1) & \le &   \frac{c_2 \, b \, 3^{\delta} \
h_B(t_2)^{\delta - \ga} }{a} \
 \ \frac{\r(u_{t_2})^{\ga} }{
\r(u_{t_2+1})^{\delta} } \ \  \#T(t_2,B)
\\ & = &
\frac{c_2 \, b \, 3^{\delta} \, \varpi }{a} \ \frac{m(B)}{f(r(B))}
\ \frac{f(\p(u_{t_2})) }{ \p(u_{t_2})^{\ga}  } \
\frac{\r(u_{t_2})^{\ga} }{ \r(u_{t_2+1})^{\delta} } \ \ \#G_B(t_2)
\ .
\end{eqnarray*}
By (\ref{NGnB}) and (\ref{t2gchoice}), it follows that $$
\begin{array}{rcl}
\#U_{B}(t_2+1) & \le & \frac{2 \, c_2 \, b \, 3^{\delta} \, \varpi
}{a} \ {\displaystyle\frac{m(B)}{f(r(B))}} \ g(u_{t_2}) \
\#G_B(t_2+1) \\[3ex] & <& \frac{1}{4} \ \#G_B(t_2+1) \ .
\end{array}
$$ For case (ii) we find that
\begin{eqnarray*}
\#U_{B}(t_2+1) & \le &  \frac{c_2 \, b \, 3^{\delta}}{a} \
\left(\frac{\r(u_{t_2})}{\r(u_{t_2+1})}\right)^{\ga}  \, \#T(t_2,
B)
\\ & = & \frac{c_2 \, b \, 3^{\delta}}{a} \
\left(\frac{\r(u_{t_2})}{\r(u_{t_2+1})}\right)^{\ga}  \ \#G_B(t_2)
\ .
\end{eqnarray*}
By (\ref{NGnB}) and (\ref{lambda}), it follows that
\begin{eqnarray*}
\#U_{B}(t_2+1) & \le &  \frac{2 \, c_2 \, b \, 3^{\delta}}{a} \
\left(\frac{\r(u_{t_2+1})}{\r(u_{t_2})}\right)^{\delta - \ga}  \,
\#G_B(t_2+1)
\\ &\leq  & \frac{2 \, c_2 \, b \, 3^{\delta}}{a} \
\lambda^{\delta-\ga}  \ \#G_B(t_2)  \ < \ \frac{1}{4} \
\#G_B(t_2+1) \ .
\end{eqnarray*}

\noindent The upshot of these estimates is that $$ \#V_{B}(t_2+1)
\ \geq \ \frac{1}{2} \ \#G_{B}(t_2+1) \  . $$

\noindent  The second local sub-level associated with $B$ is
defined to be $$ K(t_2+1,B) \ := \  \bigcup_{c\in V_B(t_2+1)} \ \
\bigcup_{c'\in \gb(t_2+1, \,  c)}^{\circ}  B(c',\p(\kdd)) \ . $$

Clearly, by construction $ K(t_2+1,B) \subset B $ and $K(t_2+1,B)
\cap K(t_2,B) = \emptyset $.

\vspace{7mm}

\noindent $\bullet$ \ {\bf The local level $K(2,B)$ : \ } For a
fixed ball $B$ in ${\bf K(1)}$ continue to construct the local
sub-levels $ K(t_2+i,B)$ associated with $B$  for $ 2 \leq i \leq
k_2(B)$, assuming of course that $k_2(B) \geq 2 $. Briefly, for $1
\leq j \leq i-1 $ suppose we have already constructed the local
sub-levels $$ K(t_2+j,B) \ := \ \bigcup_{c\in V_B(t_2+j)} \ \
\bigcup_{c'\in \gb(t_2+j, \, c)}^{\circ} B(c',\p(u_{t_2+j})) \ .
$$ For  $c \in V_B(t_2+(i-1))$ there exists a resonant set $\ra$
with $\alpha \in J^u_l(t_2+(i-1))$ such that   $c$ lies on the set
$\ra \cap \mbox{\small $\frac{1}{2}$} B $.  Construct the
thickening $$T_c(t_2+(i-1),B) \ := \ \De\left(\ra,h_B(t_2+(i-1))
\right)  \, \cap \, B(c, \r(u_{t_2+(i-1)})) \   $$ where $$
 h_B(t_2+(i-1)) \ := \ \left(\frac{\varpi \, m(B)}{f(r(B))} \
 \frac{f(\p(u_{t_2+(i-1)}))}{
\p(u_{t_2+(i-1)})^\ga}\right)^{1/(\delta - \ga)} \ , $$ and define
$$T(t_2+(i-1),B) \ := \ \left\{ T_c(t_2+(i-1),B) \, : \,  c \in
V_B(t_2+(i-1)) \right\} \ . $$ Introduce the sets \be
U_{B}(t_2+i)&:=& \{c\in G_{B}(t_2+i):B(c,\r(u_{t_2+i}))\cap
\bigcup_{j=0}^{i-1} T(t_2+j,B)\neq\emptyset \}\\
V_{B}(t_2+i)&:=&G_{B}(t_2+i)\setminus U_{B}(t_2+i) \  \ \ . \ee On
verifying that $ \#V_{B}(t_2+i) \, \geq \frac{1}{2} \ \#
G_{B}(t_2+i) $ we define $ K(t_2+i,B)$ in the obvious manner. By
construction for $0\leq i\neq j \leq k_2(B)$ we have that  $$
K(t_2+i,B) \ \cap \   K(t_2+j,B) \ = \ \emptyset \ . $$

\noindent The local level $ {\rm K(2,B)} $ associated with $B$ is
defined to be $$
 {\rm K(2,B)} \ := \ \bigcup_{i=0}^{k_2(B)}K(t_2+i,B)  \ .
 $$
In turn, the second level  of the Cantor set is defined to be $$
 {\bf K(2)} \ := \ \bigcup_{B \in {\bf K(1)} }^{\circ}  {\rm
K(2,B)}  \ . $$ Clearly, by construction ${\bf K(2)} \  \subset \
{\bf K(1)}.$

\subsubsection{Higher levels {\bf K(n)} and the Cantor set
${\bf K_{\eta}}$.}
Following the procedure for constructing {\bf K(1)} and  {\bf
K(2)}, for any integer  $n \geq 2  $ we define the $n$-th level
recursively  as follows: $$ {\bf K(n)} \ := \ \bigcup_{B \in {\bf
K(n-1)} }^{\circ}  {\rm K(n,B)}  \ , $$ where $${\rm K(n,B)} \ :=
\ \bigcup_{i=0}^{k_n(B)}K(t_n+i,B)  \ $$ is the  $n$-th local
level associated with the ball $B \in {\bf K(n-1)} $. Here $t_n >
t_{n-1} $ is chosen sufficiently large so that for any ball
$B\in{\bf K(n-1)}$ the counting estimate (\ref{NGnB}) is valid and
so that
\begin{equation}
\frac{\fp u_{t_n}))}{\p(u_{t_n})^\delta } \ > \ 3^{\delta - \ga}
\, \frac{1}{\varpi} \ \frac{f (r(B))}{m(B) }  \ .
 \label{tnchoice}
\end{equation}

\noindent Note that $  r(B) \leq \p(u_{t_{n-1}}) $ for $B \in {\bf
K(n-1)} $.  In view of  (\ref{t1choiced}), the fact that   $g(u_n)
< G^*$ for all $n$  and that $f(r)/r^{\delta}$ is decreasing as
$r$ increases,  we have that
\begin{equation}
g(u_{t_n}) \ < \ G^*  \ < \ \frac{a}{ 3^{\delta}\, 8 \, c_2 \, b}
\ \frac{1}{\varpi} \ \frac{f (r(B))}{m(B) } \  \hspace{10mm}
\forall \ \ \ \ B \in {\bf K(n-1)} . \label{tngchoice}
\end{equation}

\noindent Also in the above definition of ${\rm K(n,B)}$, the
quantity  $k_n(B)$ is the unique integer $\geq 1 $  such that
\begin{eqnarray}
\frac{ 3^{\delta}\, 2 \, c_2 \, b \, \varpi }{a} \
\frac{m(B)}{f(r(B))}
 \ \sum_{i=0}^{k_n(B)-1} g(u_{t_n+i}) &  \le &
\frac{1}{4}\label{knchoicea}\\ \nonumber \\  \frac{ 3^{\delta}\, 2
\, c_2 \, b \, \varpi }{a} \ \frac{m(B)}{f(r(B))}
 \ \sum_{i=0}^{k_n(B)} g(u_{t_n+i})& > & \frac{1}{4}\label{knchoiceb}
\end{eqnarray}

For completeness and to fix notation, we give a quick sketch of
the construction of local sub-levels $K(t_n+i,B)$ associated with
$B \in {\bf K(n-1)} $.  We define the first local sub-level
associated with $B$ in the usual manner: $$ K(t_n,B) \ := \
\bigcup_{c\in \gb(t_n)}\bigcup_{c'\in \gb(t_n, \, c)}^{\circ}
B(c',\p(u_{t_n}))  \ . $$ Any subsequent local sub-level $ K(t_n+
i ,B)$  for $1 \leq i \leq k_n(B)$,  is obtained by the following
recursive procedure. For $ 0 \leq j \leq i-1$ suppose we have
already constructed the local sub-levels $ K(t_n+ j ,B)$. Consider
some point $c \in V_B(t_n+(i-1)) \ \ (:= \gb(t_n) \ \ {\rm if \ }
\ \  i=1) $. Thus $c$ lies on the set $\ra \cap \mbox{\small
$\frac{1}{2}$} B $ for some $\alpha \in J^u_l(t_n+(i-1))$.
Construct the thickening $$T_c(t_n+(i-1),B) \ := \
\De\left(\ra,h_B(t_n+(i-1)) \right)  \, \cap \, B(c,
\r(u_{t_n+(i-1)})) \   $$ where $$
 h_B(t_n+(i-1)) \ := \ \left(\frac{\varpi \, m(B)}{f(r(B))} \
 \frac{f(\p(u_{t_n+(i-1)}))}{
\p(u_{t_n+(i-1)})^\ga}\right)^{1/(\delta - \ga)} \ , $$ and define
$$T(t_n+(i-1),B) \ := \ \left\{ T_c(t_n+(i-1),B) \, : \,  c \in
V_B(t_n+(i-1)) \right\} \ . $$ Then in view of (\ref{tnchoice})
and the fact that $f(r)/r^{\delta}$ is decreasing
\begin{equation}
3 \ \p(u_{t_n+(i-1)}) \ < \ h_B(t_n+(i-1)) \ , \label{p<h}
\end{equation}
and so by (\ref{2deltapsiB}) $$ K(t_n+ (i-1) ,B) \ \subset \
T(t_n+(i-1),B) \ . $$ Also, in view of (\ref{tngchoice}) and that
$g(u_n) < G^*$ for all $n$ we  have that   $ h_B(t_n+(i-1)) <
\r(u_{t_n+(i-1)}) $.  Next, introduce the sets \be
U_{B}(t_n\!+i)&:=& \{c\in G_{B}(t_n\!+i):B(c,\r(u_{t_n\!+i}))\cap
\textstyle\bigcup\limits_{j=0}^{i-1} T(t_n+j,B)\neq\emptyset \}\\
\\ V_{B}(t_n+i)&:=&G_{B}(t_n+i)\setminus U_{B}(t_n+i) \  \ \ . \ee
As in the previous stages of the construction, it can be verified
that  $ \#V_{B}(t_n+i) \, \geq \frac{1}{2} \ \# G_{B}(t_n+i) $.
Finally, define $$ K(t_n+i,B) \ := \  \bigcup_{c\in V_B(t_n+i)} \
\ \bigcup_{c'\in \gb(t_n+i, \, c)}^{\circ} B(c',\p(u_{t_n+i})) \ .
$$

\noindent Clearly, for $0\leq i\neq j \leq k_n(B)$ we have that $$
K(t_n+i,B) \ \cap \ K(t_n+j,B) \ = \ \emptyset \ . $$ Furthermore,
by construction the local level ${\rm K(n,B)}$ associated with $B
\in {\bf K(n-1)} $ is contained in $B$. Therefore $${\bf K(n)} \
\subset \  {\bf K(n-1)} \ . $$ The Cantor set ${\bf K_{\eta} }$ is
simply defined as $$ {\bf K_{\eta} } \ := \ \bigcap_{n=1}^{\infty}
{\bf K(n)} \ \ . $$ Trivially $$ {\bf K_{\eta} } \ \subset \
\Lambda(\psi) \ \ . $$

Before moving on to the   construction of a measure $\mu$
supported on ${\bf K_{\eta} }$, we prove an important lemma. We
adopt the notation that $V_{\Omega}(t_1) := G_{\Omega}(t_1) $ and
for $n \geq 2$ that $ V_{B}(t_n) := \gb(t_n) $.

\bigskip

\begin{lemma}\label{lem:sums}
\begin{enumerate}
\item[{\rm (i)}]  For $0\le i\le k_1$, $$\sum_{c\in
V_{\Omega}(t_1+i)} \ \sum_{c'\in G_{\Omega}(t_1+i,c)}\fp u_{t_1+i}
)) \ \ge \ \frac{c_3 \, c_4}{4} \ g(u_{t_1+i}) $$

\item[{\rm (ii)}] For $n\ge 2$, let $B$ be a ball in ${\bf K(n-1)}$.
Then, for
 $0\le i\le k_n(B)$
$$\sum_{c\in V_B(t_n+i)} \ \sum_{c'\in \gb(t_n+i,c)} \fp
u_{t_n+i})) \ \ge \ \frac{c_3 \, c_4}{4 \, b} \ g(u_{t_n+i})  \
m(B) \ \ . $$
\end{enumerate}
\end{lemma}

\smallskip

\noindent{\em Proof.}~ For either part we consider $i=0$ and $i
\geq 1$ separately.

\noindent (i)  For $i=0$, by (\ref{NGnB}) and (\ref{gbncl})
$$
\begin{array}{l}
\displaystyle \sum_{c\in V_{\Omega}(t_1)} \ \sum_{c'\in
G_{\Omega}(t_1,c)} f(\psi(u_{t_1}))  \ = \  \# V_{\Omega}(t_1)  \
\#G_{\Omega}(t_1,c) \ f(\psi(u_{t_1}))   \\[2ex] \hspace*{4ex}
\displaystyle \geq  \ f(\psi(u_{t_1})) \, \frac{1}{2} \, c_3 \,
\left(\frac{1}{\r(\kt)}\right)^{\delta} \, c_4 \,
\left(\frac{\r(\kt)}{\p(\kt)}\right)^\ga = \frac{c_3 \, c_4}{2} \
g(\kt)   .
\end{array}
$$

\smallskip
\noindent For $1 \leq i\leq k_1$, by (\ref{NGnB}) and
(\ref{gbncl}) and the fact that $ \# V_{\Omega}(t_1+i) \geq
\frac{1}{2} \,  \# G_{\Omega}(t_1+i)$
$$
\begin{array}{l}
\displaystyle \sum_{c\in V_{\Omega}(t_1+i)} \,
\sum_{c'\in
G_{\Omega}(t_1+i,c)}\fp u_{t_1+i} ))    \\[2ex] \hspace*{4ex} \displaystyle
\geq \ f(\psi(u_{t_1+i})) \, \frac{1}{4} \, c_3 \,
\left(\frac{1}{\r(u_{t_1+i})}\right)^{\delta} \!\! c_4
\left(\frac{\r(u_{t_1+i})}{\p(u_{t_1+i})}\right)^\ga \!\!\!
=\frac{c_3 \, c_4}{4}  \,  g(u_{t_1+i})  .
\end{array}
$$
\smallskip

\noindent (ii)  For $i=0$, by (\ref{NGnB}) and (\ref{gbncl}) $$
\begin{array}{l}
\displaystyle \sum_{c\in V_{B}(t_n)} \ \sum_{c'\in G_{B}(t_n,c)}
f(\psi(u_{t_n}))    \\[2ex] \hspace*{4ex} \displaystyle \geq \ f(\psi(u_{t_n})) \,
\frac{1}{2} \, c_3 \,
\left(\frac{r(B)}{\r(u_{t_n})}\right)^{\delta} \!\! c_4 \,
\left(\frac{\r(u_{t_n})}{\p(u_{t_n})}\right)^\ga \!\! \ge
\frac{c_3 \, c_4}{2 \, b } \, m(B) \,   g(u_{t_n})  .
\end{array}
$$

\noindent For $1 \leq i\leq k_n(B)$, we obtain exactly the same
estimate except for an extra `$\frac{1}{2}$' factor since $ \#
V_{B}(t_n+i) \geq \frac{1}{2} \,  \# G_{B}(t_n+i)$. \QED

\subsection{A measure on ${\bf K_{\eta} }$}
In this section, we construct a  probability measure $\mu$
supported on ${\bf K_{\eta} }$ satisfying (\ref{task}); that is,
$\mu(A)  \ll  f(r(A))/\eta $  for an arbitrary  ball $A$ of
sufficiently small radius $r(A) $.

Suppose $n \geq 2$ and $B \in {\bf K(n)}$. For $1\leq m \leq n-1$,
let $B_m$ denote the unique ball in ${\bf K(m)}$ containing the
ball $B$. With this notation in mind we now define a measure
$\mu$. For any $B \in {\bf K(n)}$, we attach a weight $\mu(B)$
defined recursively as follows:

\smallskip

For  $n=1$, $$ \mu(B)\ := \ \frac{f(r(B))}{\sum_{B'\in {\bf K(1)}}
f( r(B'))} \ $$ and for $n\ge 2$, $$ \mu(B) \ := \
\frac{f(r(B))}{\sum_{B'\in K(n,B_{n-1})} f(r(B'))} \ \times \
\mu(B_{n-1})  \ . $$

\vspace{1ex}

This procedure thus defines inductively a mass on any ball
appearing  in the construction of $\bf K$. In fact a lot more is
true --- $\mu$ can be further extended to all Borel subsets $F$ of
$\Omega$ to determine $\mu(F)$ so that $\mu$ constructed as above
actually defines a measure supported on ${\bf K_{\eta} }$; see
Proposition 1.7 \cite{falc}. We state this formally as a

\medskip

\noindent{\bf Fact.}\  The probability measure $\mu$ constructed
above is supported on ${\bf K_{\eta} }$ and for any Borel subset
$F$ of $\Omega$
\[
\mu(F):= \mu(F \cap {\bf K_{\eta} })  \; = \;
\inf\;\sum_{B\in\mathcal{ B}}\mu(B)  \ ,
\]
where the infimum is taken over all coverings $\mathcal B$ of $F \cap
{\bf K_{\eta} }$ by balls  $B\in \{{\bf K(n)} : n\in\N\}$.

\smallskip

It remains to prove the estimate (\ref{task}) on the measure $\mu
\, $.

\subsubsection{Measure of a ball in the Cantor
construction.}%
If $B\in{\bf K(n)}$ for some $n\in \N$, then by the definition of
$\mu$
\begin{eqnarray}
\mu(B) & := & \frac{f(r(B))}{\sum_{B'\in K(n,B_{n-1})}f(r(B'))}
\times  \; \mu\left(B_{n-1}\right) \nonumber \\ & & \nonumber \\
& = & \frac{f(r(B))}{\sum_{B'\in{\bf K(1)}}f(r(B'))} \;
\prod_{m=1}^{n-1} \frac{f(r(B_m))}{\sum_{B'\in
K(m+1,B_{m})}f(r(B'))} \label{measball} \;.
\end{eqnarray}

\noindent The above product term is taken to be one when $n=1$. To
proceed we require the following lemma which gives us a lower
bound on the terms in the denominator of the above expression.

\medskip

\begin{lemma}
$$  \sum\limits_{B\in{\bf K(1)}} f(r(B)) \ \ge \  \eta
\hspace{5mm} {  and   }   \hspace{3mm}
 \sum\limits_{B\in K(n, B_{n-1})} \!\!\!\!\!  f(r(B))  \ \ge \
f( r(B_{n-1}) )  \   \ \ \ ( n\ge2) . $$ \label{triv}
\end{lemma}

\noindent{\em Proof.}~ By Lemma \ref{lem:sums}, the choice of
$k_1$ (\ref{k1choiceb}) and $\varpi$  (\ref{varpi}) it follows
that
\begin{eqnarray*}
\sum_{B\in {\bf K(1)}}  f(r(B)) \ & = & \ \sum_{i=0}^{k_1}\ \
\sum_{c\in V_{\Omega}(t_1+i)} \ \sum_{c'\in G_{\Omega}(t_1+i,c)}
\fp u_{t_1+i}))\\ & \ge & \ \frac{c_3 \, c_4}{4} \
\sum_{i=0}^{k_1} g(u_{t_1+i}) \ > \ \frac{c_3 \, c_4 \, a}{32 \,
3^{\delta} \, c_2 \, b} \  \frac{\eta}{\varpi} \ \ge \ \eta   \ .
\end{eqnarray*}

\noindent For $n \geq 2$, by Lemma \ref{lem:sums}, the choice of
$k_n(B_{n-1})$ (\ref{knchoiceb}) and $\varpi$  (\ref{varpi}) it
follows that
\begin{eqnarray*}
\sum_{B\in K(n, B_{n-1})}  \!\!\!\!\!\!\!\! f( r(B)) \  &=& \
\sum_{i=0}^{k_n(B_{n-1})} \ \sum_{c\in \vbn(t_n+i)} \ \
\sum_{c'\in\gbn(t_n+i, c)} \hspace*{-4ex}\fp u_{t_n+i}))\\ & &
\\ & \ge& \frac{c_3 \, c_4 }{4 \, b}  \ m(B_{n-1}) \
\sum_{i=0}^{k_n(B_{n-1})}g(u_{t_n+i}) \\ & > &
 \frac{c_3 \, c_4 \, a}{32 \, 3^{\delta}  \,
c_2 \, b^2} \  \frac{f(r(B_{n-1}))}{\varpi}  \ \ge  \
f(r(B_{n-1}))  \ . \hspace*{0.9in} \QED
\end{eqnarray*}

\bigskip

 In view of the above lemma, it now follows from
(\ref{measball}) that for any ball $B\in{\bf K(n)}$

\begin{equation}
 \mu(B) \ \le \ \frac{f (r(B))}{\eta}  \ .
 \label{bink}
 \end{equation}

\subsubsection{Measure of an arbitrary ball
$A$.\label{secar}}
We now determine the $\mu$-measure of an arbitrary ball $A$ with
radius $r(A) \leq r_o $. The ball $A$ need not be centred at a
point of $\Omega$.  The aim is to show that: $$ \mu(A) \ \ll \
\frac{f(r(A))}{\eta} \ . $$

The measure $\mu$ is supported on ${\bf K_{\eta}}$. Thus, without
loss of generality we can assume that $ A \cap {\bf K_{\eta}} \neq
\emptyset $; otherwise $\mu(A) = 0$ and there is nothing to prove.

\vspace{4mm}

 We can also assume that for every $n$ large enough $A$ intersects at
 least two balls in  ${\bf K(n)}$; since if $B$ is the only ball in  ${\bf K(n)}$
which has non-empty intersection with  $A$, then in view of
(\ref{bink})
 $$
 \mu(A) \ \leq \ \mu(B) \ \leq  \  \frac{f(r(B))}{\eta} \ \to \ 0
 \hspace{8mm}  {\rm as } \hspace{5mm} n \to \infty  \
 $$
 ($r(B) \to 0$ as $n \to \infty$) and again there is nothing to
 prove. Thus we may assume that there exists an integer $n \geq 2$
 such that $A$ intersects only one ball $\widetilde{B}$ in ${\bf K(n-1)}$
 and at least two balls from ${\bf K(n)}$. The case that $A$
 intersects two or more balls from the first level can be excluded
 by choosing $r(A)$ sufficiently small. This follows from the fact
 that by construction balls in any one level are disjoint.
 Furthermore, we can assume that
 $$
 r(A) \ < \ r(\widetilde{B})  \ .
 $$
 Otherwise, since $f$ is increasing
 $
 \mu(A) \ \leq \ \mu(\widetilde{B}) \ \leq \ \frac{f(r(\widetilde{B}))}{\eta}
 \ \leq \ \frac{f(r(A))}{\eta}
 $
 and we are done.

\vspace{4mm}

Given that $A$ only intersects the ball $ \widetilde{B}$ in ${\bf
K(n-1)}$, the balls from level ${\bf K(n)}$  which  intersect $A$
must be contained in the local level $$ {\rm K(n,\widetilde{B})} \
:= \ \bigcup_{i=0}^{k_n(\widetilde{B})}K(t_n+i,\widetilde{B})  \ .
$$ By construction, any ball $B(c',\p(u_{t_n+i}))$ from ${\rm
K(n,\widetilde{B})}$ is contained in some thickening
$T_c(t_n+i,\widetilde{B})$. Thus $A$ intersects at least one
thickening $T_c(t_n+i,\widetilde{B}) \in T(t_n+i,\widetilde{B})$
for some $ 0 \leq i \leq k_n(\widetilde{B})$.

\vspace{6mm}

 Let $K(t_n+i',\widetilde{B})$ be the first local
sub-level associated with $\widetilde{B}$ such that $$
K(t_n+i',\widetilde{B}) \ \cap \ A \ \neq \ \emptyset \  . $$ Of
course  $ 0 \leq i \leq k_n(\widetilde{B})$ and by definition, for
any $i < i'$ $$ K(t_n+i,\widetilde{B}) \ \cap \ A \ = \ \emptyset
\  . $$ Thus, $A$ intersects at least one ball
$B(c',\p(u_{t_n+i'}))$ from $K(t_n+i',\widetilde{B})$ and such
balls are indeed the largest balls from the $n$-th level ${\bf
K(n)}$ that intersect $A$. Clearly, $A$ intersects at least one
thickening $$ T_*  \ := \ T_{c_*}(t_n+i',\widetilde{B})$$ in $
T(t_n+i',\widetilde{B})$ or equivalently with $c_* \in
V_{\widetilde{B}}(t_n+i') $. We now prove a trivial but crucial
geometric lemma.

\medskip

\begin{lemma}
\label{geo} For $i \geq i'$, if $A$ intersects
$B(c',\p(u_{t_n+i})) \subset T_c(t_n+i,\widetilde{B}) \neq T_*$
then $$ r(A) \  > \  \mbox{\footnotesize $ \frac{1}{8}$} \,
\r(u_{t_n+i}) \ . $$
\end{lemma}

\noindent{\em Proof.}~ We first establish the lemma for $i > i'$.
By definition, $ T_c(t_n+i,\widetilde{B}) \subset
B(c,\r(u_{t_n+i})) $ and by construction $ B(c,\r(u_{t_n+i})) \cap
T_* = \emptyset $. Also, in view of (\ref{2deltapsiB}) $$
B(c',\p(u_{t_n+i})) \subset B(c,\mbox{\footnotesize $
\frac{3}{4}$}\r(u_{t_n+i})) \cap \De(\ra,\p(u_{t_n+i}))  \subset
T_c(t_n+i,\widetilde{B})  . $$ Thus, there exists a point $x \in
B(c,\mbox{\footnotesize $ \frac{3}{4}$}\r(u_{t_n+i})) \cap A $.
Let $y \in T_* \cap A $. Since $ y \notin B(c,\r(u_{t_n+i}))$ we
have that $d(c,y) > \r(u_{t_n+i})$. It follows that $$ d(x,y) \
\geq \ d(y,c) \, - \, d(c,x) \ > \ \mbox{\footnotesize $
\frac{1}{4}$}\r(u_{t_n+i}) \ . $$ Since $x,y \in A$, $d(x,y) \leq
2\, r(A)$ which together with the above inequality implies that
$r(A) > \mbox{\footnotesize $ \frac{1}{8}$}\r(u_{t_n+i})$. Now
suppose $i = i'$. By construction, any thickening
$T_c(t_n+i',\widetilde{B}) \subset B(c,\r(u_{t_n+i'})) $ and the
balls $ B(c, 3 \, \r(u_{t_n+i'})) $ with $c  \in
V_{\widetilde{B}}(t_n+i') $ are disjoint. Now let $ x \in
T_c(t_n+i',\widetilde{B}) \cap A $ and $y \in  T_* \cap A $. It is
easily verified that $  d(x,y) > 2\, \r(u_{t_n+i'}) $ and so $r(A)
>\r(u_{t_n+i'})$.
\QED

In view of the definition of $i'$ and (\ref{bink}), we have that
\begin{eqnarray}
\mu(A) & \leq  & \sum_{i \, = \, i'}^{k_n(\tilde{B})} \sum_{c \,
\in \, V_{\tilde{B}}(t_n+i)} \ \sum_{\substack{ c' \, \in \,
G_{\tilde{B}} (t_n+i,c): \\ B(c',\p(u_{t_n+i})) \cap A \neq
\emptyset} } \mu(B(c',\p(u_{t_n+i}))) \nonumber \\ & & \nonumber
\\ & \leq & \frac{1}{\eta} \ \sum_{i \, = \, i'}^{k_n(\tilde{B})}
f(\p(u_{t_n+i})) \ \sum_{c \, \in \,  V_{\tilde{B}}(t_n+i)}  \
\sum_{\substack{c' \, \in \,  G_{\tilde{B}}(t_n+i,c): \\
B(c',\p(u_{t_n+i}))\cap A\neq\emptyset }}
\!\!\!\!\!\!\!\!\!\!\!\!\!\!\! 1     \label{ab} \  \  .
\end{eqnarray}

\noindent To proceed, two separate cases need to be considered:
\begin{enumerate}
\item[{\bf (i)}] when $A$ intersects at least two thickenings
in $T(t_n+i',\widetilde{B})$ \vspace{3mm}
\item[{\bf (ii)}]
when $A$ intersects only one thickening  in
$T(t_n+i',\widetilde{B})$; namely $T_*$.
\end{enumerate}

\noindent{\underline{Case} (i):}\ \ \ In view of Lemma \ref{geo},
if $A$ intersects some $T_c(t_n+i,\widetilde{B}) \in
T(t_n+i,\widetilde{B})$ then the ball $B(c,\r(u_{t_n+i}))$ which
contains $T_c(t_n+i,\widetilde{B})$ is itself contained in the
ball $17  A$. Let $N_i$ denote the number of balls
$B(c,\r(u_{t_n+i}))$ with $c \in   V_{\tilde{B}}(t_n+i) $  that
can possibly intersect $A$. By construction these balls are
disjoint  -- in fact the balls $B(c, 3\, \r(u_{t_n+i}))$ with $c
\in G_{\Omega}(t_n+i) $ are disjoint. Thus, $$ N_i \ m\left(B(3\,
\r(u_{t_n+i}))\right) \ \leq \ m(17 A)  \ . $$ Now $ A \cap {\bf
K_{\eta}} \neq \emptyset $, so there exists a point $ x \in \Omega
\cap A$. Hence, $ 17 A \subset B(x,34 \, r(A))$ and so  $m(17 A)
\leq b \, 34^{\delta} r(A)^{\delta} $. the upshot of this is that
$$ N_i \ \leq \ \frac{b \, 34^{\delta}}{a} \ \left(
\frac{r(A)}{\r(u_{t_n+i})}\right)^{\delta}   \ . $$

\noindent This implies, via (\ref{ab}) and  (\ref{gbncl}) that
\begin{eqnarray*}
\mu(A) & \leq  &  \frac{1}{\eta} \ \sum_{i \, = \,
i'}^{k_n(\tilde{B})} f(\p(u_{t_n+i})) \ \sum_{\substack{c \, \in
\, V_{\tilde{B}}(t_n+i) \\ T_c(t_n+i)\cap A\neq\emptyset } }  \
 \!\!\!\!\!\!\!\!\ \#G_{\tilde{B}}(t_n+i,c) \\ & &   \\
& \leq & \frac{c_5}{\eta} \ \sum_{i \, = \, i'}^{k_n(\tilde{B})}
f(\p(u_{t_n+i})) \   \left( \frac{\r(u_{t_n+i})}{\p(u_{t_n+i})}
\right)^{\gamma}\;N_i
\\ & &   \\ & \leq & \frac{c_5 \, b \, 34^{\delta}}{a \, \eta} \
r(A)^{\delta} \  \sum_{i \, = \, 0}^{k_n(\tilde{B})} g(u_{t_n+i})
\ .
\end{eqnarray*}
By (\ref{knchoicea}), $$  \sum_{i=0}^{k_n(B)-1} g(u_{t_n+i}) \
\leq  \   \frac{1}{ 3^{\delta}\, 8 \, c_2 \, b \, \varpi } \
\frac{f(r(\widetilde{B}))}{r(\widetilde{B})^{\delta}} \ , $$ and
by (\ref{tngchoice}) together with the fact that $g(u_n) < G^*$
for all $n$  $$ g(u_{t_n+k_n(\tilde{B})})  \ < \
\frac{f(r(\widetilde{B}))}{r(\widetilde{B})^{\delta} }  \
\frac{1}{\varpi}  \ \frac{1}{ 3^{\delta}\, 8 \, c_2 \, b}   \  .$$

\noindent Hence
\begin{equation} \mu(A) \ \ll \ \frac{1}{\eta} \
r(A)^{\delta} \
\frac{f(r(\widetilde{B}))}{r(\widetilde{B})^{\delta} } \ .
\label{case1ab}
\end{equation}

\noindent However, $r(A) < r(\widetilde{B}) $ and
$f(r)/r^{\delta}$ is decreasing. Thus $$
\frac{f(r(\widetilde{B}))}{r(\widetilde{B})^{\delta} } \leq
\frac{f(r(A))}{r(A)^{\delta} } \ , $$ which together with
(\ref{case1ab}) implies the desired inequality; namely
(\ref{task}).

\vspace{2ex}

\noindent{\underline{Case} (ii):}\ \ \ By assumption, $A$ only
intersects $T_*$ from the collection $T(t_n+i,\widetilde{B})$. On
rewriting (\ref{ab}) we have that
\begin{eqnarray}
\hspace*{-8ex} \mu(A) \ & \leq  & \  \mbox{$\sum_1$} \ := \
\frac{f(\p(u_{t_n+i'}))}{\eta} \ \sum_{\substack{c' \, \in \,
G_{\tilde{B}}(t_n+i',c_*): \\ B(c',\p(u_{t_n+i'}))\cap
A\neq\emptyset }} \!\!\!\!\!\!\!\!\!\!\!\!\!\!\! 1
 \nonumber \\ & & \nonumber
\\ \ & + & \
\mbox{$\sum_2$} \  := \ \frac{1}{\eta} \ \sum_{i \, = \,
i'}^{k_n(\tilde{B})} f(\p(u_{t_n+i})) \!\! \sum_{c \, \in \,
V_{\tilde{B}}(t_n+i)} \! \sum_{\substack{c' \, \in \,
G_{\tilde{B}}(t_n+i,c): \\ B(c',\p(u_{t_n+i}))\cap A\neq\emptyset
}} \!\!\!\!\!\!\!\!\!\!\!\!\!\!\! 1  \  \ . \label{abcase2}
\end{eqnarray}

\noindent In the case $A$ only intersects $T_*$  the  second term
$\sum_2$ on the right hand side of (\ref{abcase2}) is defined to
be zero.  In any case, we can estimate  $\sum_2$ in exactly the
same way as in case (i) to obtain that $$ \mbox{$\sum_2$} \  \ll \
\eta^{-1}f(r(A)) \  \ . $$ We now deal with the first term
$\sum_1$ on the right hand side of (\ref{abcase2}).  First note
that if $r(A) \gg \r(u_{t_n+i'}) $ then we are done; since by
(\ref{gbncl}) $$
\begin{array}{rcl}
\mbox{$\sum_1$} \ \leq \     \frac{f(\p(u_{t_n+i'}))}{\eta}  \; \#
G_{\tilde{B}}(t_n+i',c_*)  &  \ll & \frac{1}{\eta} \ g(u_{t_n+i'})
\;  \r(u_{t_n+i'})^{\delta} \\[2ex] & \ll & \frac{1}{\eta} \
g(u_{t_n+i'}) \; r(A)^{\delta} .
\end{array}
$$ However, by (\ref{tngchoice}) and the fact that $g(u_n) < G^*$
for all $n$
 $$ g(u_{t_n+i'})  \ \ll \
\frac{f(r(\widetilde{B}))}{r(\widetilde{B})^{\delta} }  \
 \  .$$
Hence $$ \mbox{$\sum_1$} \ \ll \ \frac{r(A)^{\delta}}{\eta} \
\frac{f(r(\widetilde{B}))}{r(\widetilde{B})^{\delta} }  \ \leq \
\frac{f(r(A))}{\eta}  \  \ , $$ since   $r(A) < r(\widetilde{B}) $
and $f(r)/r^{\delta}$ is decreasing. Thus, without loss of
generality we can assume that
\begin{equation}\label{rr}
 3 \, r(A) \ \leq \
\r(u_{t_n+i'})  \  .
\end{equation}
Now, $A$ must intersect at least two balls in
 $  {\rm K(n,\widetilde{B})}$. If at least two of them
 are contained in $T_*$,  say  $B(c',\p(u_{t_n+i'})) $ and
 $B(c'',\p(u_{t_n+i'})) $ with $c'\neq c'' \in
 G_{\tilde{B}}(t_n+i',c_*)$, then in view of (\ref{dxy})
 $$
  r(A) \ \geq \
\p(u_{t_n+i'})  \  . $$

\noindent On the other hand, if $A$ intersects only one ball
$B(c',\p(u_{t_n+i'})) \subset T_*$ then it must intersect some
other ball $B(c'',\p(u_{t_n+i})) \subset
T_{c}(t_n+i,\widetilde{B})$ with $i>i'$. By construction, $$
B(c'',\p(u_{t_n+i})) \ \cap \  T_* \ = \ \emptyset $$ and $$
B(c',\p(u_{t_n+i'})) \,  \subset \,   B(c',h_{\tilde{B}}(t_n+i') )
\, \subset \, T_*$$ where $h_{\tilde{B}}(t_n+i')$ is the
`thickening factor' associated  with $T_*$. Recall that $ 3 \,
\p(u_{t_n+i')})  <  h_{\tilde{B}}(t_n+i') $  -- this is
(\ref{p<h}). Now let $x \in B(c'',\p(u_{t_n+i})) \cap A $ and $ y
\in B(c',\p(u_{t_n+i'})) \cap A $. For $x \notin
B(c',h_{\tilde{B}}(t_n+i') )$ we have that $d(x,c') \geq
h_{\tilde{B}}(t_n+i') $. Thus $$ d(x,y) \ \geq \ d(x,c') \, - \,
 d(c',y) \ > \ 2\, \p(u_{t_n+i'})  \ \ , $$
 and so $r(A) > \p(u_{t_n+i'})$.   The upshot of this is that, without
 loss of generality, we can assume that
\begin{equation}
 r(A) \ \geq \
\p(u_{t_n+i'})  \  . \label{rp}
\end{equation}

In view of (\ref{rp}), it is easily verified that any ball
$B(c',\p(u_{t_n+i'}))$ which intersects $A$ is in fact contained
in the ball $3A$. In particular, in view of (\ref{2deltapsiB}) any
such ball is contained  in
 $$ 3A \, \cap \,
B(c_*,\mbox{\footnotesize $ \frac{3}{4}$}\r(u_{t_n+i'})) \, \cap
\, \De(\ra,\p(u_{t_n+i'}))  \ . $$

 \noindent Let $N$ denote the number of balls
 $B(c',\p(u_{t_n+i'}))$ with  $G_{\tilde{B}}(t_n+i',c_*)$ that
 intersect $A$. Then
\begin{eqnarray}
 m \left(  3A  \cap  B(c_*, \r(u_{t_n+i'}))  \cap
\De(\ra,\p(u_{t_n+i'}))  \right)   &\geq&    N \, m \left(
B(\p(u_{t_n+i'})) \right)  \nonumber 
\\ &\geq & a \, N \,
\p(u_{t_n+i'})^{\delta}  \ . \label{lbdab}
\end{eqnarray}

\noindent We would now like to apply the upper bound intersection
condition  to the quantity on the left hand side of (\ref{lbdab}).
However, this requires the ball $A$ to be centred at a point on a
resonant set. The following lemma is required.

\medskip

\begin{lemma}
\label{A*} Suppose $A \cap {\bf K_{\eta}} \neq \emptyset$. Then
there exists a ball $A^* \supseteq A$ with $r(A^*) \leq 3 \, r(A)$
and centre $a^*$ on a resonant set.
\end{lemma}

\noindent{\em Proof.}~ Without loss of generality assume that
 the centre $a$ of $A$ is not on a resonant set. By construction
 $ {\bf K_{\eta}}$ consists of points arbitrarily close to resonant
 sets. So if $x \in A \cap {\bf K_{\eta}}$, then there exists some
resonant  set $\ra$ such that $d(x, \ra) <  r(A)$. Thus,  $d(a,
\ra) \leq d(a, x) + d(x, \ra) < 2 \, r(A)$ which implies the
existence of some point $a^* \in \ra$ such that $d(a,a^*) < 2 \,
r(A)$. Now if $y$ is any point of $A$, then $d(y,a^*) <  3 \,
r(A)$ and this completes the proof of the lemma. \QED

In view of the lemma and (\ref{rr}), the upper bound intersection
condition implies that
\begin{eqnarray*}
{\rm L.H.S. \ of \  (\ref{lbdab})} \ &\leq&  \   m \left(  3A^* \,
\cap \, B(c_*, \r(u_{t_n+i'})) \, \cap \, \De(\ra,\p(u_{t_n+i'}))
\right) \\  \ &\leq&  \ c_2 \, b \, 9^{\delta} \
\p(u_{t_n+i'})^{\delta - \ga}   \  \ r(A)^{\ga}  \ .
\end{eqnarray*}

\noindent This together with (\ref{lbdab}) implies that $$ N \
\leq \  \frac{ \ c_2 \, b \, 9^{\delta}}{a} \ \left(
\frac{r(A)}{\p(u_{t_n+i'})} \right)^{\ga}   \ . $$

\noindent Hence,  by (\ref{rp}) and the fact that $f(r)/r^{\ga}$
is decreasing  as $ r \to 0 $
\begin{eqnarray*}
\mbox{$\sum_1$} \ := \ \frac{f(\p(u_{t_n+i'}))}{\eta} \ N   \ \ll
\ \frac{r(A)^{\ga}}{\eta}   \
\frac{f(\p(u_{t_n+i'}))}{\p(u_{t_n+i'})^{\ga}} \ \ll \
\frac{f(r(A))}{\eta} \ . \end{eqnarray*}

\noindent Thus for case (ii), $$ \mu(A) \ \leq \ \mbox{$\sum_1$} \
+ \ \mbox{$\sum_2$} \ \ll  \ \frac{f(r(A))}{\eta} \ . $$

\noindent On combining the two cases, we have shown that $  \mu(A)
\ll f(r(A))/\eta$ for an arbitrary ball $A$. This completes the
proof of  Theorem \ref{THM3} in the case that $G$ is finite. \QED

\section{Proof of  Theorem \ref{THM3}: $G= \infty$
\label{pfdimthm}}

The proof of  Theorem \ref{THM3} in the case that  $G$ is infinite
follows the same strategy as the proof when $G$ is finite. That is
to say, we construct a Cantor subset of $\Lambda(\psi)$ which
supports a certain probability measure and then apply the Mass
Distribution Theorem. However, to execute this strategy in the
case that  $G$ is infinite  is far simpler than in the finite
case. During the proof of the infinite case   we shall omit many
of the  details. After  the  proof of the finite case in the
previous section, the details should pose no real difficulties to
the reader.

To start with observe that we can assume,
 without loss of generality that $0  \leq
\gamma  <   \delta$. Also we can assume, without loss of
generality  that  $ \lim_{n \to \infty}  \p(u_n)/ \r(u_n) \,= \, 0
$.  In particular, we can assume that $\rho(u_n) > 24 \,
\psi(u_n)$ for $n$ large enough. Hence,  for an arbitrary ball
$B=B(x,r)$ with $r$ sufficiently small and $x \in \Omega $ or with
$B= \Omega $ we are able to construct the sets $\gb(n) $ and
$\gb(n,c)$ as in \S\ref{pregb}.

By definition, for each $c \in \gb(n)$ there exists an $\alpha \in
J^u_l(n) := \{ \alpha \in J : l_n < \beta_{\alpha} \leq u_n \} $
such that $c \in \ra$. In particular, by construction  for $c' \in
\gb(n,c) $ the ball $B(c', \p(u_n)) $ is contained in the ball
$B(c, \mbox{\footnotesize $ \frac{3}{4} $} \r(u_n))$ and $c' \in
\ra$.

\subsection{The Cantor set {\bf K} and the measure $\mu$ }
We start by defining a Cantor subset ${\bf K}$ of $\Lambda(\psi)$
which is dependent   on  a certain, strictly increasing sequence
of natural numbers $\{t_i : i \in \N  \}$.

\smallskip

\noindent{\em The Cantor set {\bf K}.  \  } Choose $t_1$
sufficiently large so that the counting estimate  (\ref{NGnB}) is
valid for the set $G_{\Omega}(t_1)$ and define the first level
${\bf K(1)}$ of the Cantor set ${\bf K}$  as follows: $${\bf
K(1)}\ := \ \bigcup_{c\in G_{\Omega}(t_1)} \ \bigcup_{c'\in
G_{\Omega}(t_1,c)}^{\circ} B(c', \p(\kt)) \ . $$ For $n \geq 2  $
we define the $n$-th level ${\bf K(n)}$ recursively as follows: $$
{\bf K(n)} \ := \ \bigcup_{{\rm B} \in {\bf K(n-1)} }^{\circ} {\rm
K(n,B)}  \ , $$ where $${\rm K(n,B)} \ :=  \ \bigcup_{c\in G_{{\rm
B}}(t_n)} \ \bigcup_{c'\in G_{{\rm B}}(t_n,c)}^{\circ}
B(c',\p(u_{t_n}))  $$ is the  $n$-th local level associated with
the ball ${\rm B}:= B(c', \p(u_{t_{n-1}})) \in {\bf K(n-1)} $.
Here $t_n
> t_{n-1} $ is chosen sufficiently large so that for any ball
${\rm B}$ in ${\bf K(n-1)}$ the counting estimate (\ref{NGnB}) is
valid.

The Cantor set ${\bf K}$ is simply given by $$\textstyle {\bf K }
\ := \ \bigcap\limits_{n=1}^{\infty} {\bf K(n)} \ \ . $$
Trivially,
$$
{\bf K } \ \subset \ \Lambda(\psi) \ \ .
$$

\noindent{\em The measure $\mu$.  \  } Suppose $n \geq 2$ and
${\rm B} \in {\bf K(n)}$. For $1\leq m \leq n-1$, let ${\rm B}_m$
denote the unique ball in ${\bf K(m)}$ containing the ball ${\rm
B}$. For any ${\rm B} \in {\bf K(n)}$, we attach a weight
$\mu({\rm B})$ defined recursively as follows:

For  $n=1$, $$\mu({\rm B})\ := \ \frac{1}{ \#G_{\Omega}(t_1) } \,
\frac{1}{ \#G_{\Omega}(t_1,c)} \ $$

and for $n\ge 2$, $$\mu({\rm B}) \ := \ \frac{1}{ \#G_{ {\rm
B}_{n-1} } (t_n) } \, \frac{1}{ \#G_{ {\rm B}_{n-1} } (t_n,c)}  \
\times \ \mu({\rm B}_{n-1}) \ . $$

\noindent By the definition of $\mu$ and the counting estimates
(\ref{NGnB}) and (\ref{gbncl}), it follows that

\begin{eqnarray}
\mu({\rm B}) & \leq  &   \r(u_{t_n})^{\delta -\gamma} \
\p(u_{t_n})^{\gamma} \ \left(\frac{2}{c_3 c_4}\right)^n
  \ \times \
\prod_{m=1}^{n-1} \left(\frac{\r(u_{t_m})}{\p(u_{t_m})}
\right)^{\delta -\gamma} \label{measballubexact} \\  \nonumber
\\ & \leq & f(r({\rm B})) \ c_6^n \
\frac{1}{g(u_{t_n})} \ \times \ \prod_{m=1}^{n-1}
\left(\frac{\r(u_{t_m})}{\p(u_{t_m})} \right)^{\delta -\gamma}
\label{measballub} \; ,
\end{eqnarray}
where $c_6 := 2/(c_3  c_4)  > 1 $ is a constant and  $r({\rm B})
:= \p(u_{t_n})$. The above product term is taken to be equal to
one when $n=1$.

\subsection{Completion of the proof }
Fix $\eta \geq 1 $. Since $\limsup_{n \to \infty}g(u_n) :=G =
\infty $, the sequence $\{t_i\}$ associated with the construction
of the Cantor set ${\bf K }$ can clearly be chosen so that
\begin{equation}
\eta \ \times \  c_6^i \  \times \ \prod_{j=1}^{i-1}
\left(\frac{\r(u_{t_j})}{\p(u_{t_j})} \right)^{\delta -\gamma} \
\leq \ g(u_{t_i})   \label{productchoice} \ \ . \end{equation} The
product term is one when $i=1$. It now immediately follows from
(\ref{measballub}) that for  any ${\rm B} \in {\bf K(n)}$,
\begin{equation}
\mu({\rm B}) \ \leq \ \frac{ f(r({\rm B})) }{\eta}   \ .
\label{bdone}
\end{equation}


We now show that $\mu(A) \ll  f(r(A)) / \eta $ where $A$ is an
arbitrary ball of radius $r(A) \leq r_o$. The same reasoning as in
\S\ref{secar}, enables us to assume that $A\cap {\bf K } \neq
\emptyset $, $A$ is centred on a resonant set, and that
 there exists an integer $n \geq 2$
 such that $A$ intersects only one ball $\widetilde{{\rm B}}$ in ${\bf K(n-1)}$
 and at least two balls from ${\bf K(n)}$. Thus,
without loss of generality we can assume that
\begin{equation}
 \psi(u_{t_n})  \
\leq \ r(A) \ \leq \  r(\widetilde{{\rm B}}) \, :=  \,
\p(u_{t_{n-1}}) \ \ . \label{rineq}
\end{equation} The left
hand side of the above inequality makes use of the fact that the
balls ${\rm 3 B}:= B(c', 3\, \psi(u_{t_n})) $ with ${\rm B} \in
{\bf K(n)}$  are disjoint. Consider the following two cases.

\medskip

\noindent{\underline{Case} (i):  $r(A) \leq \r(u_{t_n})$. } \ \ \
The balls  $ B(c, 3\, \r(u_{t_n})) $ with $c \in
G_{\widetilde{{\rm B}}}(t_n) $    are disjoint. Hence the ball $A$
 intersects only one  ball $B:= B(c, \r(u_{t_n}))$ with $c \in
G_{\widetilde{{\rm B}}}(t_n) $.  Let $N$ denote the number of
balls
 $  B(c', \psi(u_{t_n})) \subset B$ that can possibly intersect $A$.
The upper bound intersection condition implies that $$ m( A \cap
3B \cap \Delta(\ra,\p(u_{t_n})) \, ) \ \leq \ c_2
\p(u_{t_n})^{\delta - \gamma}  \, r(A)^{\gamma}  \ . $$ Thus, $$ N
\ \leq \ \frac{c_2}{a } \ \left( \frac{r(A)}{\p(u_{t_n})}
\right)^{\gamma} \ . $$ In view of (\ref{bdone}), (\ref{rineq})
and the fact that $f(r)/r^{\g}$ is increasing, it follows that $$
\mu(A) \ \leq \ N \, \mu({\rm B}) \ \ll \ \frac{1}{\eta} \
r(A)^{\gamma} \; \frac{f(\p(u_{t_n}))}{\p(u_{t_n})^\gamma} \ \leq
\ \frac{f(r(A))}{\eta} \ . $$

\medskip

\noindent{\underline{Case} (ii):  $r(A) > \r(u_{t_n})$. } \ \ \ If
$A$ only intersects one ball  $ B(c, \r(u_{t_n})) $ with $c \in
G_{\widetilde{{\rm B}}}(t_n) $, then by (\ref{gbncl}),
(\ref{bdone}) and the fact that $f(r)/r^{\g}$ is increasing, we
have that
\begin{eqnarray*}
\mu(A) & \leq &  \mu( B(c, \r(u_{t_n})) \, ) \  := \  \sum_{c' \in
G_{\widetilde{{\rm B}}}(t_n,c)} \mu ( B(c', \psi(u_{t_n})) \, )
\\ & \leq &  c_5 \,  \mu({\rm B}) \  \left(
\frac{\r(u_{t_n})}{\p(u_{t_n})} \right)^{\gamma} \ll \
\frac{1}{\eta} \ r(A)^{\gamma} \;
\frac{f(\p(u_{t_n}))}{\p(u_{t_n})^\gamma} \ \leq \
\frac{f(r(A))}{\eta} \ .
\end{eqnarray*}

\noindent Without loss of generality, assume that $A$ intersects
at least two  balls $ B(c, \r(u_{t_n})) $ with $c \in
G_{\widetilde{{\rm B}}}(t_n) $ and let $N$ be the number of such
balls that can possibly intersect $A$.  A simple geometric
argument, making use of the fact that the balls $B(c,
\r(u_{t_n}))$ are disjoint yields  that $$ N \ \leq \
\frac{3^{\delta} \, b }{a}  \; \left( \frac{r(A)}{\r(u_{t_n})}
\right)^{\delta}  \ . $$  In view of (\ref{gbncl}),
(\ref{measballubexact}),   (\ref{productchoice}), (\ref{rineq})
and the fact that $f(r)/r^{\delta}$  is decreasing,  we obtain
that
\begin{eqnarray*}
\mu(A) & \leq &  N \, \mu( B(c, \r(u_{t_n})) \, ) \  = \ N \
 \# G_{\widetilde{{\rm B}}}(t_n,c) \  \mu({\rm B})
 \\  & \ll &    r(A)^{\delta} \ \left(\frac{2}{c_3 c_4}\right)^n
  \ \times \
\prod_{m=1}^{n-1} \left(\frac{\r(u_{t_m})}{\p(u_{t_m})}
\right)^{\delta -\gamma} \\ & \ll &  f(r(A)) \ \ c_6^{n} \ \frac{
\p(u_{t_{n-1}})^{\delta}}{f(\p(u_{t_{n-1}})) }  \ \times \
\prod_{m=1}^{n-1} \left(\frac{\r(u_{t_m})}{\p(u_{t_m})}
\right)^{\delta -\gamma} \\ & \ll &  f(r(A))  \ \ c_6^{n-1} \
\frac{ 1}{g(u_{t_{n-1}}) }  \ \times \ \prod_{m=1}^{n-2}
\left(\frac{\r(u_{t_m})}{\p(u_{t_m})} \right)^{\delta -\gamma} \
\leq \  \frac{f(r(A))}{\eta} \ .
\end{eqnarray*}

The upshot of these cases is that $\mu(A)  \ll f(r(A))/\eta $ for
an arbitrary ball $A$. By the Mass Distribution Principle,
 $\mathcal{ H}^f(\Lambda(\psi)) \geq \mathcal{ H}^f({\bf
K} ) \gg \eta $.  However, $\eta \geq 1 $ is arbitrary whence
$\mathcal{ H}^f(\Lambda(\psi)) = \infty $. This thereby completes the
proof of  Theorem  \ref{THM3} in the case that  $G$ is infinite.
\QED


\section{Applications \label{app}}

Unless stated otherwise, in all the following statements the
convergent parts are easy to establish -- just use the `natural
cover' given by the $\limsup$ set under consideration. Also  we
will make use of the following simple fact. Suppose that $h: \R^+
\to \R^+ $ is a real, positive monotonic  function, $\alpha \in
\R$ and $k>1$.  Then the divergence and convergence properties of
the sums $$ \sum_{n=1}^{\infty} k^{n \, \alpha} \; h(k^n)
\hspace{1cm} {\rm and }  \hspace{1cm} \sum_{r=1}^{\infty} r^{
\alpha - 1 } \; h(r) \hspace{1cm} {\rm coincide. }   $$

The various  applications have been chosen to illustrate the
versatility of our general framework. There are many  other
applications, such as to inhomogeneous Diophantine approximation,
Markov maps and iterated function schemes, which we have decided
not to include -- mainly to avoid repetition and to keep the
length of the paper manageable. Throughout, `i.m.' is short for
`infinitely many'.

\subsection{Linear Forms \label{lfapp}}
Let $\psi$ be an approximating  function. An $m \times n $ matrix
$X = (x_{ij}) \in {\bf I}^{mn} := [0,1]^{mn}$ is said to be
$\psi$--well approximable if the system of inequalities $$ |q_1 \,
x_{1j} \, + \, q_2 \, x_{2j} \, + \dots +
 \, q_m \, x_{mj} \; - \;
p_j | \ <\ \psi \left( |\bq| \right) \ |\bq|
 \ \ \ \ \ \ (1 \leq j \leq n)
$$ is satisfied for infinitely many vectors $ \bq\in \Z^m \, $,
$\bp\in\Z^n \;$. Here $ | \bq | $ denotes the supremum norm of the
vector $\bq \,$; i.e. $ |\bq| = \max \{  |q_1| \, , \, \dots \, ,
\, |q_m| \} \,$.
 The
system $$ q_1 \, x_{1j} \, + \, q_2 \, x_{2j} \, + \dots +  \, q_m
\, x_{mj} \ \ \ \ \ \ (1 \leq j \leq n) $$ of $n$ real linear
forms in $m$ variables $ q_1 \, , \, \dots \, , \, q_m $ will be
written more concisely as $ \bq \, X   $, where the matrix $X$ is
regarded as a point in ${\bf I}^{mn} \, $. In view of this
notation, the set of $\psi$--well approximable points will be
denoted by $$ W \, \left(m,n;\psi\right) :=\{X\in {\bf
I}^{mn}:|\bq X-\bp| < \psi\left(|\bq|\right) |\bq|\ {\rm for \
i.m.}\ (\bp , \bq )\in\Z^n \times \Z^m \} \ . $$ By definition,
 $ |\bq X-\bp| \,
= \, \max_{1\leq j \leq n} |\bq . X^{(j)} \, - \, p_j| $ where $
X^{(j)} $ is the $j$'th column vector of $X$.  Note that  when $m
= n = 1$, the set $W(1,1;\psi)$ corresponds to the classical set
$W(\psi)$ in the one dimensional theory.

With reference to our general framework, let  $\Omega:= {\bf
I}^{mn}$, $J:= \{ ({\bf p},{\bf q}) \in \Z^n \times \Z^m
\backslash \{ {\bf 0}\}: |\bp| \leq |\bq| \}$, $\alpha :=
(\bp,\bq) \in J$, $\ma := |\bq|$  and
  $\ra:=  \{ X \in  {\bf I}^{mn}:  \bq X = \bp \}$.
Thus, the  family $\mathcal{ R}$ of resonant sets $\ra$ consists of
$(m-1)n$--dimensional, rational hyperplanes.   Furthermore,
$\De(\ra,\p(\ma)) := \{ X \in {\bf I}^{mn}:  {\rm dist} (X, \ra) <
\psi\left(|\bq|\right) \} \  $ and $$ \De^u_l(\p, n) :=
\bigcup_{2^{n-1} < |\bq| \leq 2^n} \ \bigcup_{|\bp| \leq |\bq| }
\De(\ra,\p(\ma)) \ . $$ Then  $$W \left(m,n;\psi\right)=
\Lambda(\p) := \limsup_{n \to \infty} \De^u_l(\p, n) \ \ . $$

\noindent Now let the measure  $m$  be $mn$--dimensional Lebesgue
measure, $ \delta = mn$ and  $\gamma= (m-1)n $.  Then a
probabilistic argument involving mean and variance techniques
yields the following statement  \cite[\S4.3]{det}.

\medskip

\begin{proposition} The pair $(\mathcal{ R},\beta) $ is
a local $m$--ubiquitous system relative to $(\r,l,u)$, where
$l_{t+1} = u_t := 2^t \ (t \in \N) $
 and $  \r : r \to \r(r) \:=  {\rm constant} \, \times \,
r^{-(m+n)/n}  \omega(r) $. Here $\omega$ is any real, positive
increasing function such that $1/\omega(r) \to 0 $ as $r \to
\infty $ and such that  for any  $C>1$ and  $r$ sufficiently large
$\omega(2\,r) < C \, \omega(r) $. \label{proplf}
\end{proposition}

In view of the proposition and the fact that the measure $m$ is of
type (M2) and that $\rho$ is $u$-regular, Corollary \ref{cor2} of
Theorem \ref{THM1} and Theorem \ref{THM3} yield the divergent
parts of the following statements.

\medskip

\begin{thgros}
Let $\psi $ be a real, positive decreasing function. Then $$
m(W(m,n;\psi) ) = \left\{
\begin{array}{ll}
 0 & {\rm if} \;\;\; \sum_{r=1}^\infty \;  \psi\left(r\right)^{n}
     \;\; r^{m+n-1} <\infty\; ,\\[2ex]
 1 & {\rm if} \;\;\; \sum_{r=1}^\infty \; \psi\left(r\right)^{n} \;\;
     r^{m+n-1} =\infty \; .
\end{array}
\right. $$
\end{thgros}

\medskip

\begin{thdet}
Let $f$ be a dimension function such that $r^{-mn} \, f(r)\to
\infty$ as $r\to 0 \, $ and  $r^{-mn} \, f(r) $ is decreasing.
Furthermore, suppose that $ r^{-(m-1)n} \, f(r)$ is increasing.
Let $\p $ be  a real, positive  decreasing function. Then $$
\hf\left(W(m,n;\p)\right)= \left\{
\begin{array}{ccl}
 0 && {\rm if} \ \sum   f\left(\p(r)\right)
     \p\left(r\right)^{-\left(m-1\right)n}  r^{m+n-1} <\infty
     ,\\[2ex]
 \infty && {\rm if} \ \sum  f\left(\p(r)\right)
          \p\left(r\right)^{-\left(m-1\right)n}
         r^{m+n-1} =\infty .
\end{array}
\right. $$
\end{thdet}

\smallskip

Notice that the function $\omega$ associated with $\rho$ does not
appear in the above statements. This is a consequence of  choosing
$\omega$ appropriately.  With $\r$ as in the proposition,
Corollary \ref{cor2}  implies that
$$m(W(m,n;\psi) )  \, = \, 1  \  \  {\rm if}   \ \
\mbox{$\sum_{\omega}$} \, := \, \sum \; \psi\left(r\right)^{n}
\;\; r^{m+n-1} \ \ \omega(r)^{-n} \ = \ \infty \ . $$   However,
to obtain the precise statement of the Khintchine-Groshev theorem
we need to remove the $\omega$ factor. To do this we choose
$\omega$ in such a way that the divergence/convergence properties
of $\sum_{\omega}$ and $\sum := \sum \psi(r)^{n} \ r^{m+n-1}$ are
the same. It is always possible to find such a  function. Clearly,
if $\sum_{\omega}= \infty$ then $\sum = \infty$. On the other
hand, if $\sum = \infty$, then we can find a strictly increasing
sequence of positive integers $\{r_i\}_{i\in\N}$ such that $$
\sum_{r_{i-1}<r \leq \ r_i} \psi\left(r\right)^{n} \ r^{m+n-1}  \
> \ 1 \ ,
 $$ and $ r_i > 2r_{i-1}$.   Now simply  define  $\omega $ be the
 step function
 $\omega(r) := i^{\frac{1}{n}}$ for $r_{i-1}<r\leq r_i$ and
this satisfies the required properties. A similar argument allows
us to conclude Theorem DV from Theorem \ref{THM3}  without the
presence of the $\omega$ factor.

\smallskip

\noindent{\em Remark.  } The above theorems remain valid if the
 set $W \left(m,n;\psi\right)$ is replaced by its
`inhomogeneous' analogue. Briefly, for a vector ${\bf b} \in {\bf
I}^{n} := [0,1]^{n}$ consider the set $ W_{\bf b}
\left(m,n;\psi\right) :=\{X\in {\bf I}^{mn}:|\bq X-\bp - {\bf b} |
< \psi\left(|\bq|\right) |\bq|\ {\rm for \ i.m.}\ (\bp , \bq
)\in\Z^n \times \Z^m \} . $ Obviously, the `homogeneous' set $W
\left(m,n;\psi\right)$ corresponds to the case when ${\bf b}$ is
the zero vector. Now define $\Lambda(\p)$ as above with the only
modification being that the family $\mathcal{ R}$ of  resonant sets
$\ra$  now consists of  $(m-1)n$--dimensional hyperplanes of the
form $ \{ X \in {\bf I}^{mn}: \bq X = \bp + {\bf b} \} $. Then, it
is possible to show that Proposition \ref{proplf} remains valid
for the  pair $(\mathcal{ R}, \beta) $ and so Corollary \ref{cor2} of
Theorem \ref{THM1} and Theorem \ref{THM3} yield the divergent
parts of the Khintchine-Groshev theorem and Theorem DV for the set
$W_{\bf b} \left(m,n;\psi\right)$.  The proof of the proposition
in the inhomogeneous setup follows from Theorem 2 in \cite{yann2}
if $m=1$ and from the mean-variance argument carried out in \S3.2
of \cite{jasan} if  $m \geq 2$.

\subsection{Algebraic Numbers}
Let $H(a)$ denote the height of an algebraic number $a$, that is
the maximum of the absolute values of the relatively prime integer
coefficients in its minimal defining polynomial. For $d \in \N$,
denote by  $A(d)$  the set of algebraic numbers $a$ with degree at
most $n$.   Given an approximating function $\psi$, let $$
K_d(\psi) := \{\xi \in [0,1] :|\xi - a| < \p \left( H(a) \right)
\; {\rm for\ i.m.\ } \; a \in A(d) \} \ . $$ The set $K_d(\psi)$
is a generalization of $W(\psi)$ since the rationals are algebraic
with degree one. In the case that $\psi(r) = r^{-(d+1)\tau} $ let
us write $K_d(\tau)$ for $K_d(\psi)$.  A. Baker and W.M. Schmidt
\cite{BS} have obtained the following analogue of the classical
Jarn\'{\i}k--Besicovitch Theorem.

\medskip

\begin{thbs}
\ \ For  $\tau \geq 1$, $ \dim K_d(\tau) =1/\tau  $.
\end{thbs}

As an application of our  main theorems we are able to give a
complete measure theoretic description of $K_d(\psi)$ which not
only implies the Baker-Schmidt Theorem but also shows that  $
\mathcal{ H}^{1/\tau}  \! \left( K_d(\tau) \right) = \infty $.


Let  $\Omega:= [0,1] \, ,  \,  J:= \{ a: a \in A(d)  \} \, ,  \,
\alpha := a \in J \, , \,  \ma \!\! := H(a) \ {\rm and \ }
  \ra:=  a. $ Thus, the  family $\mathcal{ R}$ of resonant sets $\ra$ consists
  of points corresponding to
algebraic numbers $a \in A(d)$. Furthermore, $\De(\ra,\p(\ma)) :=
B(a,\p( H(a)) ) $ and
 $$\De^u_l(\p, n) :=   \bigcup_{a \in J^u_l(n) }
 \!\!\!\!\!  B\left( a,\p( H(a)) \, \right) \ ,
 $$
 where $J^u_l(n) := \{ a \in A(d):  k^{n-1} < H(a) \leq k^n \} $.
 Here $k > 1 $ is a constant.
Then  $$K_d(\p)  \, = \, \Lambda(\p) \, := \, \limsup_{n \to
\infty} \De^u_l(\p, n) \ \ . $$

\noindent Now let $m$  be one--dimensional Lebesgue measure, $
\delta = 1$ and  $\gamma=0$.

\medskip

\begin{proposition}\label{prp1}
 The pair $(\mathcal{ R}, \beta) $ is a local
$m$--ubiquitous system relative to $(\rho,l,u)$ where for $k >
k_0$ -- a positive absolute constant, $l_{t+1} = u_t := k^t$ $(t
\in \N)$ and $\r : r \to \r(r)  :=  {\rm constant} \, \times \,
r^{-(d+1)}$.
\end{proposition}

Baker and Schmidt \cite{BS} established  Proposition~\ref{prp1}
with $\r(r) = r^{-(d+1)} \times (\log r)^{3d(d+1)}$. This is
sufficient only to determine the dimension result  (simply apply
Corollary \ref{dimthm}). The presence of the $\log$ term in their
ubiquity function $\rho$  rules out the possibility of obtaining
the more desirable measure theoretic laws for $K_d(\p)$. However,
a more  subtle analysis enables one to remove the $\log$ term
\cite{vb1}.

 In view of the proposition and the fact that the measure $m$ is of
type (M2) and that $\r$ is  $u$-regular, Corollary \ref{cor2} of
Theorem \ref{THM1}  and Theorem \ref{THM3} yield the divergent
parts of the following statements.

\medskip

\begin{theorem}
Let $\psi $ be  a real, positive decreasing function. Then $$ m
\left(K_d(\psi)\right)= \left\{
\begin{array}{ll}
  0 & {\rm if} \;\;\; \sum_{r=1}^\infty \; \psi\left(r\right) \;\; r^d
  <\infty\; ,\\[2ex]
  1 & {\rm if} \;\;\; \sum_{r=1}^\infty \; \psi\left(r\right)
   \;\; r^d =\infty \; .
\end{array}\right.
$$ \label{A1}
\end{theorem}

\medskip

\begin{theorem}
 Let $f$ be a dimension function such that
$r^{-1} \, f(r)\to \infty$ as $r\to 0 \, $ and  $r^{-1} \, f(r) $
is decreasing. Let $\p $ be  a real, positive  decreasing
function. Then $$ \hf\left(K_d(\psi))\right)= \left\{
\begin{array}{ll}
  0 & {\rm if} \;\;\; \sum_{r=1}^\infty \;  f\left(\p(r)\right) \;\; r^{d}
      <\infty\; ,\\[2ex]
  \infty & {\rm if} \;\;\; \sum_{r=1}^\infty \; f\left(\p(r)\right) \;\;  r^{d}
  =\infty\; .
\end{array}
\right. $$ \label{A2}
\end{theorem}

Theorem \ref{A1}   was first established in \cite{vb1}. As
mentioned above, Theorem \ref{A2}  not only implies the
Baker-Schmidt Theorem but also shows that  $ \mathcal{ H}^{1/\tau }
\left( K_d(\tau) \right) = \infty $. A weaker form of Theorem
\ref{A2}, 
has been recently established in
\cite{yann}.

\subsection{Kleinian Groups \label{appkg}}
The classical results of Diophantine approximation,\hspace{2pt} in
particular those from the one dimensional theory, have natural
counterparts and extensions in the hyperbolic space setting. In
this setting, instead of approximating real numbers by rationals,
 one approximates limit points of a fixed   Kleinian group $G$ by
 points in the orbit (under the group) of a certain distinguished
 limit point $y$.
Beardon and Maskit  have shown that the geometry of the group is
reflected in the approximation properties of points in the limit
set.  The elements of $G$ are orientation preserving
 M\"obius transformations of the $(n+1)$--dimensional
 unit ball $B^{n+1}$. Let $\Lambda $ denote the limit set of $G$ and
 let $\delta$ denote the Hausdorff dimension of $\Lambda$.
 For any element $g$ in $G$ we shall use the notation
 $L_g := |g^\prime(0)|^{-1}$, where $|g^\prime(0)|$ is the (Euclidean)
 conformal dilation of $g$ at the origin.

Let $\psi $ be an approximating  function and let $$ W_{y}(\psi)
:=\{ \xi \in \Lambda: | \xi - g(y) | < \psi(L_g) \, \mbox{for i.m.
$g$ in $G$}\}  . $$

\noindent  This is the set of points in the limit set $\Lambda$
which are `very  close' to   infinitely  many images  of  a
`distinguished' point $y$. The `closeness' is of course governed
by the approximating function $\psi$. The limit point $y$ is taken
to be a parabolic fixed point if the group has parabolic elements
and a hyperbolic fixed point otherwise.

\vspace{2mm}

\noindent{\bf Geometrically finite groups: \ } Let us assume that
the geometrically finite group has parabolic elements so it is not
convex co-compact. Thus our distinguished limit point $y$ is a
parabolic fixed point, say $p$. Associated with $p$ is a
geometrically motivated set $\mathcal{ T}_p$ of coset representatives
of $G_p \backslash G := \{g G_p : g \in G \} $; so chosen that if
$ g \in \mathcal{ T}_p$ then the orbit point $g(0)$ of the origin lies
within a bounded hyperbolic distance from the top of the standard
horoball $H_{g(p)}$. The latter, is an $(n+1)$--dimensional
Euclidean ball contained in $B^{n+1}$ such that its boundary
touches the unit ball $S^n$ at the point $g(p)$.   Let $R_g$
denote the Euclidean radius of  $H_{g(p)}$.  As a consequence of
the definition of $\mathcal{ T}_p$, it follows that $$ \frac{1}{C \,
L_g} \leq R_g \leq \frac{C}{L_g} \  $$ where $C>1$ is an absolute
constant. Also, it is worth mentioning that the balls in the
standard set of horoballs $\{H_{g(p)} \, : \,  g \in \mathcal{ T}_p \}
$ corresponding to the parabolic fixed point $p$ are pairwise
disjoint.   For further details and references  regarding the
above notions and statements see any of the papers
\cite{jbgfg,melian,Strap}. With reference to our general
framework, let
 $\Omega:= \Lambda  \  , \  J:= \{ g: g \in \mathcal{ T}_p  \} \ , \ \alpha :=
g \in J \ , \ \ma := C \, R_g^{-1} \ {\rm and \ }
  \ra:=  g(p). $ Thus, the  family $\mathcal{ R}$ of resonant sets $\ra$ consists
  of orbit points $g(p)$
with $g \in \mathcal{ T}_p$. Furthermore, $\De(\ra,\p(\ma)) :=
B(g(p),\p( C \, {R_g}^{-1} ) ) $ and
 $$\De^u_l(\p, n) :=   \bigcup_{g \in J^u_l(n)}
 \!\!  B\left(g(p),\p( C \, R_g^{-1} ) \right) \ , $$
 where $J^u_l(n) := \{ g \in \mathcal{ T}_p: k^{n-1} < \,  C \, R_g^{-1} \leq
 k^n \} $. Here $k >1 $ is a constant.
Then  $$   W_{p}(\psi) \, \supset \, \Lambda(\p) \, := \,
\limsup_{n \to \infty} \De^u_l(\p, n) \ \ . $$

\noindent Now, let $m$  be  Patterson measure,  $ \delta = \dim
\Lambda$ and $\gamma=0$. Thus $m$ is a non-atomic,
$\delta$--conformal probability measure supported on $\Lambda$.
Furthermore, $m$ is of type (M1) with respect to the sequences
$l:= \{k^{t-1}\} $ and $u:=\{k^t\} $ for any $k > 1$ -- see below.
In fact, the condition that $m(B(c,2r)) \ll m(B(c,r)) $ for balls
centred at resonant points is valid for any $c \in \Lambda$. We
have the following statement concerning  local ubiquity.

\medskip

\begin{proposition} The pair $(\mathcal{ R}, \beta) $ is a local
$m$--ubiquitous system relative to $(\r,l,u)$ where  for $k \geq
k_o $ -- a positive group constant, $$ l_{t+1} = u_t := k^t \ \ (t
\in \N)  \ \ \ \ { and } \ \ \ \ \r : r \to \r(r) \ := \ {\rm
constant} \ \times \ r^{-1} \ . $$  \label{kgps}
\end{proposition}

 The proposition follows from the following three facts which can
be found in  \cite{jbgfg,melian}.

\smallskip

\noindent {\bf $\bullet$} {\em Local Horoball Counting Result: \ }
Let $B$ be an arbitrary Euclidean ball in $S^n$ centred at a limit
point. For $\lambda \in (0,1)$ and $r \in \R^+$ define $$
A_{\lambda}(B,R) \, := \, \{ g \in \mathcal{ T}_p: g(p) \in B {\rm \
and \ } \lambda R \leq R_g < R \}  \ . $$ There exists a positive
group constant $\lambda_o$ such that if $\lambda  \leq \lambda_o $
and $R < R_o(B)$, then

$$ k_1^{-1} \,  R^{-\delta} \,  m(B) \ \leq \ \# A_{\lambda}(B,R)
\ \leq \  k_1 \,  R^{-\delta} \,  m(B)  \ , $$ where $k_1$ is a
positive constant independent of $B$ and $R_o(B)$ is a
sufficiently small positive  constant  which does depend on $B$.

\noindent {\bf $\bullet$} {\em Disjointness Lemma: \ } For
distinct elements  $g,h \in \mathcal{ T}_p $  with $ \lambda < R_g/R_h
< \lambda^{-1}$, one has $ B(g(p),\lambda R_g ) \, \cap \,
B(h(p),\lambda R_h ) \ = \ \emptyset  \ . $

\smallskip

\noindent {\bf $\bullet$} {\em Measure of balls centred at
parabolic points: \ } For $g \in \mathcal{ T}_p $ and $r \leq R_g$ $$
 k_2^{-1} \, r^{2 \delta - {\rm
rk}(p) } \, R_g^{ \,  {\rm rk}(p) - \delta }\ \leq \ m\left(
B(g(p),r) \, \right) \ \leq  \  k_2 \, r^{2 \delta - {\rm rk}(p) }
\, R_g^{ \, {\rm rk}(p) - \delta }  \  , $$ where ${\rm rk}(p) $
denotes the rank of the parabolic fixed point $p$ and $k_2 > 1 $
is a positive constant independent of $g$ and $r$. Clearly, this
implies that  $m$ satisfies condition (M1) with respect to the
sequences $l$ and $u$.

\vspace{2mm}

To prove the proposition, let $\rho(r) := C (k\, r)^{-1}$ where
$k:= 1/\lambda
> 1/\lambda_o$ and $B$ be an arbitrary ball centred at a limit
point. Then for $n$ sufficiently large
\begin{eqnarray*}
m (\  B   \cap \!\!\!\!\!\!\!\!\!\! \bigcup_{\substack{ g \in
\mathcal{ T}_p: \\ k^{n-1} <  C \, R_g^{-1} \leq k^n} }^{\circ}
 \!\!\!\!\!\!\!\!\!\! B\left(g(p),\r( k^n ) \right)   \ )   &\ge &
 m  (   \bigcup_{\substack{g \in \mathcal{ T}_p: \, g(p) \in \frac{1}{2} B \\ k^{n-1} < \,
 C \, R_g^{-1} \leq k^n} }^{\circ}
 \!\!\!\!\!  B\left(g(p),\r( k^n ) \right)   \ )  \\ & &
\\
&\gg &  k^{-n \, \delta}  \ \#  A_{\frac{1}{k}}(\mbox{\small
$\frac{1}{2}$} B \, , C \, k^{-(n-1)} ) \\
&&\\
& \gg & m (\mbox{\small
$\frac{1}{2}$} B) \ \gg \ m(B) \  .
\end{eqnarray*}

\noindent  Now let $\p$ be an approximating function and assume
without loss of generality that $\p( k^n) \leq \r( k^n) $  for $n$
sufficiently large.  If this were not the case then the $\limsup$
condition of Theorem \ref{THM1} can be invoked to imply the
desired result below.  Since  $\p(k^n) < R_g $ for  $g \in
J^u_l(n)$, the above measure fact for balls centred at parabolic
points implies that for any $g \in J^u_l(n)$ $$ m\left(
B\left(g(p),\p( k^n ) \right) \right) \asymp \p(k^n)^{ 2 \delta -
{\rm rk}(p) } \ k^{-n ({\rm rk}(p) - \delta) } \ . $$

\noindent Also notice that  $
 m\left(
B\left(g(p),\r( k^n ) \right) \right) \asymp k^{-n \delta} \ $. It
therefore follows that
\begin{eqnarray*} {\rm L.H.S. \ of \ (\ref{af}) } \  &
\leq & \sum_{s=1}^{Q-1} \ k^{s \delta}  \ \sum_{s+1 \, \leq \, t
\leq \, Q } \p(k^t)^{ 2 \delta - {\rm rk}(p) } \  \  k^{-t ({\rm
rk}(p) - \delta)  }   \\
 & =  &  \sum_{m=2}^{Q} \ \p(k^m)^{ 2 \delta - {\rm rk}(p) } \  k^{-m ({\rm
rk}(p) - \delta)  }  \ \sum_{r=1}^{m-1} k^{r \delta}  \\ & \ll &
 \sum_{m=2}^{Q} \ \p(k^m)^{ 2 \delta - {\rm rk}(p) } \  k^{m (2 \delta - {\rm rk}(p))  } \
 \ll \  {\rm R.H.S. \ of \ (\ref{af}) } \ \ .
 \end{eqnarray*}

\noindent  Thus, in view of Proposition \ref{kgps} and the fact
that the measure $m$ is of type (M1) and that $\r$ is $u$-regular,
Theorem \ref{THM1}  yields the divergent part of the following
statement.
\medskip

\begin{theorem}
Let $G$ be a geometrically finite Kleinian group with
pa\-ra\-bo\-lic elements and let  ${\rm rk}(p)$ denote  the  rank
of the parabolic fixed point $p$.  Let $\psi $ be  a real,
positive decreasing function. Then $$ m( W_{p}(\psi) ) = \left\{
\begin{array}{ll}
  0 & {\rm \ if} \;\;\;
  \sum_{r=1}^\infty \;  \psi\left(r\right)^{2 \delta - {\rm rk}(p)}
  \;\; r^{2\delta-{\rm rk}(p)- 1} <\infty\; ,\\[2ex]
  1 & {\rm \ if} \;\;\; \sum_{r=1}^\infty \;
 \psi\left(r\right)^{2 \delta - {\rm rk}(p)} \;\; r^{2\delta-{\rm
 rk}(p)- 1} =\infty \; .
\end{array}
\right. $$ \label{GF1}
\end{theorem}

\vspace{-3mm}

Theorem \ref{GF1} is not new. However, in previous statements
 of the theorem  a certain  regularity
condition on  $\psi$ is  assumed \cite{Paddyrs,Strap,SDS}. In
Theorem \ref{GF1}, the regularity is removed and replaced by the
`natural' condition that $  \psi \ $ is decreasing. Thus the above
Khintchine type theorem is the perfect analogue of the classical
statement.

 In general, for
geometrically finite Kleinian groups with parabolic elements,
Patterson measure $m$ is not of type (M2). Thus, Theorem
\ref{THM3}  is not applicable even though we have local
$m$-ubiquity (Proposition \ref{kgps}). In fact, in general $m$ is
not even of type (M$2^{\prime}$) so Theorem $\ref{THM3}^{\prime}$
of  \S\ref{r3} is  not applicable either. However, if the group is
of the first kind (so $\Lambda = S^n$) then $m$ is normalized
$n$-dimensional Lebesgue measure on the unit sphere $S^n$ and so
is certainly of type (M2). Also, for groups of the first kind
$\delta = n = {\rm rk}(p)$. Thus, for such groups Theorem
\ref{THM3} yields the divergent part of the following statement.

\medskip

\begin{theorem}\label{GF2}
Let $G$ be a geometrically finite Kleinian group of the first kind
with parabolic elements and $p$ be a parabolic fixed point.
 Let $f$ be a dimension function such that
$r^{-n} \, f(r)\to \infty$ as $r\to 0 \, $ and  $r^{-n} \, f(r) $
is decreasing.  Let $\psi$  be a real,  positive decreasing
function. Then $$ \hf\left( W_{p}(\psi)\right)= \left\{
\begin{array}{ll}
 0 & {\rm if}
 \;\;\; \sum_{r=1}^\infty \;  f\left(\p(r)\right) \;\;
 r^{n-1} <\infty\; ,\\[2ex]
 \infty & {\rm if} \;\;\; \sum_{r=1}^\infty \;
 f\left(\p(r)\right) \;\;  r^{n-1} =\infty \; .
\end{array}
\right. $$
\end{theorem}

 Regarding Theorem
\ref{GF2}, all that was previously known were dimension statements
for $ W_{p}(\psi)$ \cite{slv3}.  For example, in the case $\psi(r)
= r^{-\tau}$ let us write $W_{p}(\tau) $ for $ W_{p}(\psi)$. Then
$\dim W_{p}(\tau) = n/\tau $ $(\tau \geq 1)$. Clearly, Theorem
\ref{GF2} implies this statement and shows that the
$s$--dimensional Hausdorff measure of $W_{p}(\tau)$ at the
critical exponent $s=n/\tau$ is infinite. For completeness, we
mention  the following dimension result \cite{jbgfg}. Let $G$ be a
geometrically finite Kleinian group with parabolic elements and
let  ${\rm rk}(p)$ denote  the rank of the parabolic fixed point
$p$. Then for $\tau \geq 1 $ $$
 \dim W_{p} (\tau)
 \;=\;
 \min \left\{  \frac{\delta  + \rk(p)  \, (\tau - 1) }{2 \tau - 1},  \,
 \frac{\delta}{\tau}\right\}\;.
$$ So for groups of the second kind, although the  dimension of
$W_{p} (\tau)$ is known its Hausdorff measure at the critical
exponent is unknown. As already mentioned above,  our general
framework fails to shed any light on this,  since although for
groups of the second kind we are able to establish local
$m$-ubiquity (Proposition \ref{kgps}) the measure $m$  is not of
type (M2) or even of type (M$2^{\prime}$). Recently, it has been
shown that for sets closely related to $ W_{p}(\psi)$ the
Hausdorff measure is either zero or infinite \cite{Dodv}. However,
even for these related sets one is unable to establish the
analogue of Theorem \ref{GF2} for groups of the second kind.

When interpreted on the upper half plane model $\Half^2$ of
hyperbolic space and applied to the modular group ${\rm
SL}(2,\Z)$, the above theorems imply the classical results
associated with our basic example -- see \S\ref{beg}.  Next, let
$G_d$ denote the Bianchi group of 2x2 matrices of determinate one
with entries in the ring of integers $\vartheta =\vartheta(d) $ of
the imaginary quadratic field $\Q(\sqrt{-d})$. Here $d$ is a
positive integer which is not a perfect square.  For a real,
positive decreasing function  $\p $, let $W_{\vartheta}(\p) $
denote the set of complex numbers $z$ such that the inequality $$
|z \,  - \, p/q | \ < \ \p(|q|) $$ is satisfied for i.m. pairs
$p,q \in \vartheta \times \vartheta$ with ${\rm ideal}(p,q) =
\vartheta $. Following  \S7 of \cite{SDS}, it is easily verified
that when interpreted on the upper half space model $\Half^3$ of
hyperbolic space and applied to the Bianchi group $G_d$, the above
theorems imply the following statements.

\medskip

\begin{theorem}
Let $\psi$  be  a real, positive decreasing function and let  $m$
denote $2$--dimensional Lebesgue measure. Then $$ m
\left(W_{\vartheta}(\p)\right)=\left\{
\begin{array}{ll}
  0 & {\rm if} \;\;\; \sum_{r=1}^\infty \; \psi\left(r\right)^2 \;\; r^3
 <\infty \; ,\\[2ex]
  \infty & {\rm if} \;\;\; \sum_{r=1}^\infty \;
 \psi\left(r\right)^2 \;\; r^3 =\infty \; .
\end{array}
\right. $$ \label{c1}
\end{theorem}

\medskip

\begin{theorem}
 Let $f$ be a dimension function such that
$r^{-2} \, f(r)\to \infty$ as $r\to 0 \, $ and  $r^{-2} \, f(r) $
is decreasing. Let $\psi$  be a real,  positive decreasing
function. Then $$ \hf\left(W_{\vartheta}(\p)\right)=\left\{
\begin{array}{ll}
  0 & {\rm if} \;\;\; \sum_{r=1}^\infty \;  f\left(\p(r)\right) \;\;
 r^{3} <\infty \; ,\\[2ex]
 \infty & {\rm if} \;\;\; \sum_{r=1}^\infty \;
 f\left(\p(r)\right) \;\;  r^{3} =\infty \; .
\end{array}
\right. $$ \label{c2}
\end{theorem}

Theorem \ref{c1} is essentially due to Sullivan \cite{SDS}.
However, Sullivan assumed a certain regularity condition on $\p$.
This has been replaced  by the  more natural condition that $
\psi$ is decreasing. Thus, Theorem \ref{c1} is the precise
analogue of the classical statement of Khintchine. In the case
$\psi(r) = r^{-\tau} \!\!  \ $  write $W_{\vartheta}(\tau) $ for $
W_{\vartheta}(\psi)$. Then, Theorem \ref{c2} implies the following
`complex' analogue of the Jarnik--Besicovitch theorem.

\medskip

\begin{corollary} For $\tau \! \geq \!2$,
$ \dim  W_{\vartheta}(\tau)  \! = \! \frac{4}{\tau}$. Moreover,
$\mathcal{ H}^{4/\tau} \!\! \left(W_{\vartheta}(\tau)\right) \!\! =
\!\! \infty$.
\end{corollary}

\medskip

 \noindent{\bf Convex co-compact groups:}~~~ These are geometrically
finite Kleinian groups without parabolic elements. Thus, the
distinguished limit point $y$ is a hyperbolic fixed point. For
convex co-compact groups, Patterson measure $m$ is of type (M2)
and the situation becomes much more satisfactory.

 Let $L$ be the axis of the
conjugate pair of hyperbolic fixed points $y$ and $y'$, and let
$G_{yy'}$ denote the stabilizer of $y$ (or equivalently $y'$).
Then there is a geometrically motivated set $\mathcal{ T}_{yy'} $  of
coset representatives of $G_{yy'} \backslash G $; so chosen that
if $ g \in \mathcal{ T}_{yy'}$ then the orbit point $g(0)$ of the
origin lies within a bounded hyperbolic distance from the summit
$s_g$ of $g(L)$ -- the axis of the hyperbolic fixed pair $g(y)$
and $g(y')$.  The summit $s_g$ is simply the point on $g(L)$
`closest' to the origin. For $g \in \mathcal{ T}_{yy'}$, let
$H_{g(y)}$ be the horoball with base point at $g(y)$ and radius
$R_g := 1 - |s_g| $. Then the top of $H_{g(y)}$ lies within a
bounded hyperbolic distance of $g(0)$. Furthermore, as a
consequence of the definition of $\mathcal{ T}_{yy'}$, it follows that
$ C^{-1}  \leq R_g \, L_g \leq C $  where $C>1$ is an absolute
constant. We are now able to define the subset $\Lambda(\p)$ of
$W_{y}(\psi)$ in exactly the same way as in the parabolic  case
with $y$ replacing $p$ and $\mathcal{ T}_{yy'}$ replacing $\mathcal{
T}_{p}$.

Essentially the arguments given in \cite{melian}, can easily be
modified to obtain the analogue of the local  horoball counting
result stated above for the parabolic case.  We leave the details
to the reader. In turn, this enables  one to establish Proposition
\ref{kgps} for convex co-compact groups -- the statement remains
unchanged.  Since $m$ is of type (M2) and $\rho$ is $u$--regular
for any $k>1$, Corollary \ref{cor2} of Theorem \ref{THM1} and
Theorem \ref{THM3} yield the divergent parts of the following
statements.

\medskip

\begin{theorem}
Let $G$ be a convex co-compact Kleinian group and $y$ be a
hyperbolic fixed point. Let $\psi$  be a real,  positive
decreasing function. Then $$ m( W_{y}(\psi) ) = \left\{
\begin{array}{ll}
 0 & {\rm if} \;\;\; \sum_{r=1}^\infty \;  \psi\left(r\right)^{\delta} \;\;
r^{\delta-1} <\infty \; , \\[2ex]
 1 & {\rm if} \;\;\; \sum_{r=1}^\infty \; \psi\left(r\right)^{\delta} \;\;
r^{\delta-1} =\infty \;  .
\end{array}\right.
$$ \label{CC1}
\end{theorem}

\medskip

\begin{theorem}
Let $G$ be a convex co-compact Kleinian group and $y$ be a
hyperbolic fixed point. Let $f$ be a dimension function such that
$r^{-\delta} \, f(r)\to \infty$ as $r\to 0 \, $ and  $r^{-\delta}
\, f(r) $ is decreasing.  Let $\p $ be a real, positive decreasing
function. Then $$ \hf\left( W_{y}(\psi)\right)= \left\{
\begin{array}{ll}
  0 & {\rm if} \;\;\; \sum_{r=1}^\infty \;  f\left(\p(r)\right) \;\;
 r^{\delta-1} <\infty \; , \\[2ex]
 \infty & {\rm if} \;\;\; \sum_{r=1}^\infty \;
 f\left(\p(r)\right) \;\;  r^{\delta-1} =\infty \; .
\end{array}\right.
$$ \label{CC2}
\end{theorem}

Again, as in the parabolic case, the first of these theorems is
not new. However, in previous statements of the theorem  a certain
regularity condition on the function $\psi$ is  assumed.  Thus the
above Khintchine type theorem is the perfect analogue of the
classical statement. Regarding the second theorem, all that was
previously known were  dimension statements such as $\dim
W_{y}(\psi) = \delta/\tau $ $(\tau \geq 1)$ when $\psi(r) =
r^{-\tau}$. Clearly, Theorem \ref{CC2} implies this statement and
shows that the $s$--dimensional  Hausdorff measure of
$W_{y}(\psi)$ at the critical exponent $s= \delta/\tau$ is
infinite.

\subsection{Rational Maps}
We consider a special case of the general `shrinking target'
problem introduced in \cite{inv}. Let  $T$ be an expanding
rational map (degree $\geq 2$) of the Riemann sphere $
\overline{{\C}} = \C \cup \{ \infty \}$ and $J(T)$ be its Julia
set. For any $z_o \in J(T)$ and   $\psi$ an approximating
function, consider the set $$ \Wt
 = \{ z \in J(T) :   T^n(z) \in  B \, (  z_o , \psi(|\Tnd(z)| ))
 \hbox{ for i.m. } n\in\N \}  .
$$ In view of the bounded distortion property for expanding maps (
see Proposition 1, \cite{inv}), there exists a constant $C=C(T)>
1$ such that    the
 set of points $z$ in $J(T)$ which lie in the ball
 $$ B \, \left( y , \frac{\psi(C \, |\Tnd(y)| )}{C \, |\Tnd(y)|} \right) $$
 for i.m. pairs $(y,n) \in I :=  \{ (y,n) : n \in \N  {\rm \  with \ }
 T^n(y) =  z_o  \} $ is a subset of $\Wt$. On the other hand, if we replace
 $C$ by $C^{-1}$ in the above ball
 then  $\Wt$ is a subset of the  corresponding  set of points. It is now clear
 that $\Wt$ is a $\limsup $ set of the type considered within our framework.
 The backward orbit of the selected point $ z_o$ in $J(T)$
corresponds to the rationals in the classical theory. This set is
also the precise analogue of the set of well approximable limit
points associated with a Kleinian group.

With reference to our general framework, let
 $\Omega:= J(T)  \  , \  J:= I \ , \ \alpha :=
(y,n) \in I , \ \ma := C \, |\Tnd(y)| \ {\rm and \ }
  \ra:=  y. $ Thus, the  family $\mathcal{ R}$ of resonant sets $\ra$ consists
  of pre-images of the point $z_o$. Furthermore,
  define  $\De(\ra,\varphi(\ma)) :=
B(y ,\varphi( C \, |\Tnd(y)| ) \, ) $ and let
 $$\De^u_l(\varphi,n) :=   \bigcup_{(y,m) \in J^u_l(n): }
 \!\!\!\!\!  B\left(y,\varphi(  C \, |\Tmd(y)| ) \,  \right) \ , $$
 where $J^u_l(n) := \{ (y,m) \in I : k^{n-1} < \,   C \, |\Tmd(y)| \leq
 k^n \} $ and $\varphi (r) :=   r^{-1} \, \psi(r) $. Here $k>1$ is
 a constant.
In view of the discussion above,     $$   \Wt \, \supset \,
\Lambda(\varphi) \, := \, \limsup_{n \to \infty}
\De^u_l(\varphi,n) \ \ . $$

\noindent Now let $m$  be  Sullivan measure,  $ \delta = \dim
J(T)$ and $\gamma=0$. Thus $m$ is a non-atomic,
$\delta$--conformal probability measure supported on $J(T)$ and
since $T$ is expanding it is of type (M2).  We have the following
statement concerning local ubiquity.

\medskip

\begin{proposition} The pair $(\mathcal{ R}, \beta) $ is a local
$m$--ubiquitous system relative to $(\r,l,u)$ where for $k \geq
k_o $ -- a positive constant dependent only on the rational map
$T$,  $$l_{t+1} = u_t := k^t   \ \ (t \in \N)  \ \ \ \ { and } \ \
\ \ \r : r \to \r(r) \ := \ {\rm constant} \ \times \ r^{-1} \ .$$
\label{julia}
\end{proposition}

The proposition follows from the following two facts which can be
found in \cite{ihes}. For ease of reference we keep the same
notation and numbering of constants as in \cite{ihes}. For $X \in
\R^+$, let $I(X)$ denote the set of pairs $(y,n) \in I$ such that
$$ f_n(y) - C_8  \ \leq \ X \ \leq \ f_{n+1}(y) + C_8 \ , $$ where
$f_n(y) := \log |\Tnd(y)|$; i.e. the $n$-th ergodic sum of  $f=
\log|T'|$.

\noindent {\bf $\bullet$} {\em Constant Multiplicity: \ } Let $z
\in J(T)$. Then there are no more than $C_9$ pairs $(y,n) \in
I(X)$ such that $$ z \, \in \, B \left( y, C_{\!10} \;
|\Tnd(y)|^{-1} \right) \ . $$

 This is the second part of the statement of Lemma 8 in
\cite{ihes}.

\noindent {\bf $\bullet$} {\em Local Counting Result: \ } Let $B$
be an arbitrary Euclidean ball  centred on a point of $J(T)$. Then
there exists a constant $X_o(B)$ such that for $X \geq X_o(B)$ $$
\# \{ (y,n) \in I(X) : y \in B \} \ \asymp \ m(B) \  e^{\delta X}
\ , $$ where the implied constants are independent of $B$.  This
statement is the last line of the proof of Theorem 4 in
\cite{ihes}.

\vspace{2mm}

To prove the proposition, let $\r(r) := C_{\!10} \ r^{-1}$ and $B$
be an arbitrary ball centred on a point of $J(T)$. Then for  $k >
e^{2\, C_8} \  |T'(z_o)| $ and $n$ sufficiently large  we have
that \begin{eqnarray*}
m (  \ B  \ \  \cap  \!\!\!\!\! \bigcup_{ (y,m) \in
J^u_l(n) }
 \!\!\!\!\!  B\left(y,\r( k^n ) \right)   \ )  & \geq &
 m  (   \bigcup_{\substack{(y,m) \in J^u_l(n): \\ y \in \frac{1}{2} B  }}
 \!\!\!\!\!  B\left(y,\r( k^n ) \right)   \ )  \\
&&\\ & \geq  &m  (   \bigcup_{\substack{(y,m) \in I(X): \\ y \in
\frac{1}{2} B } }
 \!\!\!\!\!  B\left(y,\r( k^n ) \right)   \ ) \\
&&\\
&\geq & C_9^{-1} \!\!\!\!\!\!\!  \sum_{(y,m) \in I(X): \, y \in
 \frac{1}{2} B}
 \!\!\!\!\!  m \left( B\left(y,\r( k^n ) \right)  \,  \right)
 \\
&&\\
& \gg & \r(k^n)^{\delta}  \ \  \# \{ (y,m) \in I(X):  y \in
 \mbox{\small $\frac{1}{2}$} B \}
 \\
&&\\
& \gg &\displaystyle k^{-n \, \delta}
 \ e^{\delta  X} \ m (\mbox{\small $\frac{1}{2}$} B)
\ \gg \ m(B)  \ \ ,
\end{eqnarray*}
where $X := n\log k - \log C - C_8$.
The dependence of $z_o$ from the size of $k$ can be removed by
choosing $k
> e^{2\, C_8} \ \sup\{|T'(z)| : z \in J(T)\}  $.

In view of Proposition \ref{julia},   the fact that the measure
$m$ is of type (M2) and that $\r$ is $u$-regular for any $k>1$,
Corollary \ref{cor2} of Theorem \ref{THM1} and Theorem \ref{THM3}
yield the divergent parts of the following statements.

\medskip

\begin{theorem}
Let $\psi$  be  a real, positive decreasing function.  Then $$ m(
\Wt ) = \left\{
\begin{array}{ll}
 0 & {\rm if} \;\;\; \sum_{r=1}^\infty \;
 \psi\left(r\right)^{\delta} \;\; r^{-1} <\infty \; ,\\[2ex]
 1 & {\rm if} \;\;\; \sum_{r=1}^\infty  \; \psi\left(r\right)^{\delta} \;\;
r^{-1} =\infty \;  .
\end{array}
\right. $$ \label{J1}
\end{theorem}

\medskip

\begin{theorem}
 Let $f$ be a dimension function such that
$r^{-\delta} \, f(r)\to \infty$ as $r\to 0 \, $ and  $r^{-\delta}
\, f(r) $ is decreasing.  Let $\p$    be a real, positive
decreasing function. Then $$ \hf\left( \Wt \right)= \left\{
\begin{array}{ll}
 0 & {\rm if} \;\;\;
 \sum_{r=1}^\infty \;  f\left(\p(r)/r\right) \;\;
 r^{\delta-1} <\infty \; , \\[2ex]
 \infty & {\rm if} \;\;\; \sum_{r=1}^\infty \;
 f\left(\p(r)/r\right) \;\;  r^{\delta-1} =\infty \; .
\end{array}
\right. $$ \label{J2}
\end{theorem}

In \cite{inv}, Hausdorff dimension results for  $ \Wt$ were
established. For example,  in the case $\psi(r) = r^{-\tau}$ write
$W_{z_0} (\tau)$ for the set $\Wt$. Then for $\tau \geq 0$,  $\dim
W_{z_0} (\tau) = \delta/(1+ \tau) $. Recently, it has been shown
that for $\tau > 0 $ the $\delta/(1+ \tau)$--dimensional Hausdorff
measure of $W_{z_0} (\tau)$ is either zero or infinity
\cite{Hvel}.  Clearly, Theorem \ref{J2} implies the dimension
statement and shows that the Hausdorff measure  at the critical
exponent  is actually infinite.

Theorem \ref{J1} enables us to deduce the following `logarithmic
law' for orbit approximation. Fix a point $z_0$ in $J$ and for any
other point $z$ in $J$ let $d_n(z)$ denote the distance of
$T^n(z)$ from $z_0$. Then for $m$-almost all points $z$ in $J$ $$
\limsup_{n \to \infty} \frac{-\log d_n(z)}{\log \log |\Tnd(z)|} =
\frac{1}{\delta} \ . $$ This statement can be viewed as the
(expanding) rational map analogue of Sullivan's logarithmic law
for geodesics 
\cite{Strap,SDS}.

\subsection{Diophantine approximation with restrictions}
In a series of papers, G. Harman has considered the problem of
obtaining zero-one laws for hybrids of the  classical set
$W(\psi)$ in which the  numerator and denominator of the rational
approximates are  restricted to various sets of number theoretic
interest.  We refer the reader to Chapter 6 in \cite{Har} for a
full exposition.

To illustrate the diversity of our main theorems, we consider a
specific case in which the number theoretic set is the set of
prime numbers $\mathcal{ P}$. Thus, let $\psi$ be an approximating
function and let
$$
W_\mathcal{ P}(\psi):= \{x \in [0,1] :|x - p/q| < \p(q) \ {\rm for\
i.m.\ } \; (p,q) \in \mathcal{ P} \times \mathcal{ P} \} \ . $$ With
reference to our general framework, let $\Omega:= [0,1] \ , \ J:=
\{ (p,q) \in \mathcal{ P} \times \mathcal{ P} : p \leq q \} \ , \ \alpha
:= (p,q) \in J \ , \ \ma := q \ {\rm and \ } \ra:=  p/q. $ Thus,
the  family $\mathcal{ R}$ of resonant sets $\ra$ consists of
rationals $p/q$ with both numerator and denominator prime.
Furthermore, define $\De(\ra,\p(\ma)) := B(p/q ,\p( q )   ) $ and
let
 $$\De^u_l(\p, n) :=  \!\!\!\! \bigcup_{(p,q) \in J^u_l(n): }
  B(p/q,\p(q) )   = \!\!\!\!\! \bigcup_{\substack{q \in \mathcal{ P} :  \\ 2^{n-1} < q \leq
  2^n}}
  \ \bigcup_{p \in \mathcal{ P}: \,  p \leq q} B(p/q, \p(q)) \ , $$
 where $J^u_l(n) := \{ (p,q) \in J : 2^{n-1} < \,  q \leq
 2^n \} \, .$ Then
$$W_\mathcal{ P}(\psi) = \Lambda(\p) :=  \limsup_{n \to \infty}
\De^u_l(\p, n) \,  . $$

Now let $m$ be one--dimensional Lebesgue measure, $\delta = 1$ and
$\gamma=0 $. Then a relatively standard  analytic argument making
use of sieve theory leads to the following  local $m$--ubiquity
statement.

\medskip

\begin{proposition} The pair $(\mathcal{ R}, \beta) $ is a local
$m$--ubiquitous system relative to $(\r,l,u)$, where $l_{t+1} =
u_t := 2^t$, $(t \in \N)$  and $\r : r \to \r(r)  := \, {\rm
constant}
 \times  (\log  r)^{2} r^{-2}$. \label{harman}
\end{proposition}

The $\log$ term in the function $\r$  is, of course,  a
consequence of the prime number theorem. The proposition can be
deduced from Lemma 6.3 in \cite{Har}. However, for the details see
\cite{becky}. Thus, in view of the proposition and the fact that
the measure $m$ is of type (M2) and that $\r$ is $u$-regular,
Corollary \ref{cor2} of Theorem \ref{THM1} and Theorem \ref{THM3}
yield the divergent parts of the following statements.

\medskip

\begin{theorem}\label{Re1}
Let $\psi$  be  a real, positive decreasing function. Then $$
W_\mathcal{ P}(\psi)  = \left\{\begin{array}{ll}
 0 & {\rm if} \;\;\; \sum_{r=1}^\infty \;  \psi(r)  \;\; r \; (\log r)^{-2}  \; < \infty \;
 ,\\[1.5ex]
 1 & {\rm if} \;\;\; \sum_{r=1}^\infty \; \psi(r)  \;\; r \;\; (\log r)^{-2} =\infty \; \; .
\end{array}\right.
$$
\end{theorem}

This zero-one law was first established by Harman. Theorem
\ref{THM1} shows that the statement is in fact a simple
consequence of local $m$-ubiquity.

\medskip

\begin{theorem}\label{Re2}
 Let $f$ be a dimension function such that
$r^{-1} \, f(r)\to \infty$ as $r\to 0 \, $ and  $r^{-1} \, f(r) $
is decreasing. Let $\psi$  be a real,  positive decreasing
function. Then $$ \hf\left(W_\mathcal{ P}(\psi)\right)= \left\{
\begin{array}{ll}
 0 & {\rm
if} \;\;\; \sum_{r=1}^\infty \;  f\left(\p(r)\right) \;\;  r \;
(\log r)^{-2} <\infty\\[1.5ex]
 \infty & {\rm if} \;\;\; \sum_{r=1}^\infty \;
f\left(\p(r)\right) \;\;  r \; (\log r)^{-2} =\infty \; .
\end{array}
\right. $$
\end{theorem}

Consider the case when $\psi(r) = r^{-\tau}$. Then, Theorem
\ref{Re2} implies that $\dim W_\mathcal{ P}(\psi) = 2/\tau \ \
(\tau \geq 2) $ just as in the classical case. However, the
$2/\tau$--dimensional Hausdorff measure is zero unlike the
classical case in which it is infinite. Thus, restricting to prime
numerators and denominators has no effect on the dimension but
drastically effects the Hausdorff measure. The above mentioned
dimension result is not new and can be found in  \cite[Theorem
10.8]{Har}.

For analogous results associated with problems in which the
numerator and denominator of the rational approximates are
restricted to other sets of number theoretic interest see
\cite{becky}.

\subsection{Diophantine approximation in $\Q_p$}
For a  prime $p$, let  $|\ |_p$ denote the $p$--adic metric and
let $\Q_p$ denote the $p$--adic field. Furthermore, let $\Z_p$
denote the ring of $p$--adic integers and for a vector $\x$ in
$\Q_p^m$ let $|\x|_p:=\max \{  |x_1|_p \, , \, \dots \, , \,
|x_m|_p \} \,$ -- the p--adic norm on the vector space $\Q_p^m$.
We now consider the $p$--adic analogue of the `classical' set  $W
\left(m,n;\psi\right)$  -- see \S\ref{lfapp}.  For an
approximating function $\psi$, let $$
\begin{array}{l}
\hspace*{-0ex} W_p  \left(m,n;\psi\right) := \{X\in\Z_p^{mn}:|\bq
X+\bp|_p<
\psi(\max(|\bp|,|\bq|))\\[1ex] \hspace*{41ex} {\rm for \ i.m.}\ (\bp , \bq
)\in\Z^n \times \Z^m \}  ,
\end{array}
$$ where  $|\x|:=\max \{  |x_1| \, , \, \dots \, , \, |x_m|\} \,$
is the usual supremum norm of a vector $\x $ in $ \Z^m $. By
definition,
 $ |\bq X+\bp|_p \,
= \, \max_{1\leq j \leq n} |\bq . X^{(j)} \, + \, p_j|_p $ where $
X^{(j)} $ is the $j$'th column vector of the $m \times n $ matrix
$X  \in  \Z_p^{mn}$.

  There are two points worth making when comparing the
above set with the `classical'  set $W  \left(m,n;\psi\right)$.
Firstly, the approximating function in the $p$--adic setup is a
function of $ \max(|\bp|,|\bq|)$ rather than simply $|\bq|$.  This
is due to the fact that within the $p$--adic setup for any
$X\in\Z_p^{mn}$ and $\bq\in\Z^m$ there exists $\bp\in\Z^n$ such
that the quantity  $|\bq X+\bp|_p$ can be made arbitrarily small.
Thus, the set of $X\in\Z_p^{mn}$ for which    $|\bq X+\bp|_p<
\psi(|\q|)$ for  i.m.  $ (\bp , \bq )\in\Z^n \times \Z^m $ is in
fact  the whole space  $\Z_p^{mn}$ and there is nothing to prove.
Secondly,  in the $p$-adic setup the `normalizing' factor $|\bq|$
does not appear on the right hand side of $|\bq X+\bp|_p$. This is
due to the fact that the $p$--adic metric is an ultra metric.
Thus, if $|\bq X+\bp|_p<\psi(\max(|\bp|,|\bq|)) $ then $X$ lies in
the $\psi(\max(|\bp|,|\bq|)) $ neighborhood of $ \{
X\in\Z_p^{mn}:\bq X+\bp=0\} $ -- the resonant set associated with
the pair $(\bp , \bq )$ (see below).

With reference to our general framework, let   $\Omega:=
\Z_p^{mn}, \ J:= \{ ({\bf p},{\bf q}) \in \Z^n \times \Z^m
\backslash \{ {\bf 0}\}: |\bp| \leq |\bq| \} \ , \ \alpha :=
(\bp,\bq) \in J \ , \ \ma := |\bq| \ {\rm and \ }
  \ra:=  \{ X\in\Z_p^{mn}:\bq
X+\bp=0\}. $ Thus, the  family $\mathcal{ R}$ of resonant sets $\ra$
consists of $(m-1)n$--dimensional sets. Furthermore,
$\De(\ra,\p(\ma)) := \{ X \in \Z_p^{mn}:  {\rm dist} (X, \ra) \leq
\psi(|\bq|) \} \  $ and let $$\De^u_l(\p, n) :=   \bigcup_{2^{n-1}
< |\bq| \leq 2^n} \ \bigcup_{|\bp| \leq |\bq| }   \De(\ra,\p(\ma))
\ . $$ The metric $d$ is of course the $p$--adic norm.  Then $$W
\, \left(m,n;\psi\right) \supset \Lambda(\p) := \limsup_{n \to
\infty} \De^u_l(\p, n) \ \ . $$

\noindent Now let $m$ be the standard  Haar measure on
$\Q_p^{mn}$, $\de=mn$ and $\ga=(m-1)n$. Thus $m(\Z_p^{mn})=1$ and
$ m(B(x, p^{-t})) = (p^{-t})^{mn}$ for $t \in \N$. Then a
probabilistic argument involving mean and variance techniques
yields the following statement.

\medskip

\begin{proposition} The pair $(\mathcal{ R}, \beta) $ is a local
$m$--ubiquitous system relative to $(\r,l,u)$ with $\r, l$ and $u$
as in Proposition \ref{proplf}. \label{proplfp}
\end{proposition}

The  ideas necessary for the proof can be found  in
 \cite[Chp. 4]{lutz}. Briefly, for $N \in \N $ define the set
$\fN$ to be empty if  $|N|_p<1$  and if $|N|_p =1 $ let  $\fN$ be
the set of $({\bf p},{\bf q}) \in \Z^n \times \Z^m $ such that
\begin{enumerate}
\item $|\br|=|\q|=q_1=N$
\item    $p_j$ is co-prime to $q_1$ with
$0<p_j< N \, \omega(N)^{-1/2}$ for $j=1,\dots,n$.
\end{enumerate}

\noindent The argument on page 86 together with the mean--variance
argument beginning on  page 94 and   Lemma 4.1 on page 70 imply
that for $n$ large enough $
 m\left(  \De^u_l(\r,n) \right)  \ge  \ka  > 0 $, where $\rho $
 is as in Proposition \ref{proplfp}. This proves
global $m$-ubiquity. However, much more is true.  A simple
applications of Lemma 1.1 on page 14 enables us to conclude that $
 m\left(  \De^u_l(\r,n) \right)  \to   1  $
 as $ n \to   \infty  $.
The required  local $m$-ubiquity statement immediately follows --
see \S\ref{thesystems}.

In view of the proposition and the fact that the measure $m$ is of
type (M2) and that $\rho$ is $u$-regular, Corollary \ref{cor2} of
Theorem \ref{THM1}  and Theorem \ref{THM3} yield the divergent
parts of the following statements. They are the $p$-adic analogues
of the theorems stated in \S\ref{lfapp}.

\medskip

\begin{theorem}\label{lfpkg}
Let $\psi $ be a real, positive decreasing function. Then $$
m(W_p(m,n;\psi) ) = \left\{
\begin{array}{ll}
 0 & \ {\rm if} \;\;\;
\sum_{r=1}^\infty \;  \psi\left(r\right)^{n} \;\; r^{m+n-1}
<\infty \; ,\\[1.5ex]
 1 & \ {\rm if} \;\;\; \sum_{r=1}^\infty \; \psi\left(r\right)^{n}
\;\; r^{m+n-1} =\infty \; .
\end{array}\right.
$$
\end{theorem}

\medskip

\begin{theorem}\label{lfphm}
Let $f$ be a dimension function such that $r^{-mn} \, f(r)\to
\infty$ as $r\to 0 \, $ and  $r^{-mn} \, f(r) $ is decreasing.
Furthermore suppose that $r^{-(m-1)n}f(r)$ is increasing. Let $\p
$ be  a real, positive  decreasing function. Then $$
\hf\left(W_p(m,n;\p)\right)= \left\{
\begin{array}{ll}
 0 & {\rm \ if} \ \sum_{r=1}^\infty  \  f\left(\p(r)\right)
\  \p\left(r\right)^{-\left(m-1\right)n}  r^{m+n-1}
 <\infty  ,\\[1.5ex]
 \infty & {\rm \ if} \  \sum_{r=1}^\infty \
f\left(\p(r)\right)  \ \p\left(r\right)^{-\left(m-1\right)n}
r^{m+n-1} =\infty .
\end{array}\right.
$$
\end{theorem}

Theorem \ref{lfpkg} is not new and first appeared (in full
generality) in \cite{lutz}; the specific case $m=1$  is due to
Jarn\'{\i}k \cite{Jarpa}. Theorem \ref{lfphm} is new.  Previously
\cite{Abercrombie}, it had been shown that $ \dim  W_p(m,n;\tau) =
(m-1)n + (m+n)/\tau $ for $ \tau  > (m+n)/n $. As usual,
$W_p(m,n;\tau)$ corresponds to the set  $W_p(m,n;\p)$  with  $\p:
r \to r^{-\tau} $. Obviously, Theorem \ref{lfphm} implies this
dimension result and shows that $\mathcal{ H}^s (W_p(m,n;\tau)) =
\infty $ at the critical exponent $s =   \dim W_p(m,n;\tau) $.


\subsection{Diophantine approximation on manifolds}
Let $M$ denote  an $m$--dimensional submanifold in $\R^n$ with
$n\ge2$. Given an approximating function $\p$, the problem is to
determine measure theoretic laws for points ${\bf x} \in \R^n$
resticted to the manifold $M$. This restriction means that the
points ${\bf x}= \{x_1, ...,x_n\} $ of interest consist of
dependent variables reflecting the fact that ${\bf x} \in M $. The
fact that the variables are dependent unlike in the classical
setup where $M = \R^n$, introduces major difficulties in
attempting to describe the measure theoretic structure of
$\limsup$ sets restricted to $M$. There are two main types of
$\limsup$ sets that can be considered.

\subsubsection{Dual/linear approximation on manifolds.}
For an approximating function $\p$, let $$W(M;\psi):= \{\x\in M:
|\q.\x-p|<|\q|\psi(|\q|)\ \mathrm{for\ i.m.\
}(p,\q)\in\Z\times\Z^n\},$$ where $\bx.\by=x_1y_1+\dots+x_ny_n$ is
the normal scalar product of two vectors $\bx,\by$ in $\R^n$ and
$|\bx|$ is the usual supremum norm of the vector $\bx\in\R^n$. To
make any reasonable progress we impose the condition that the
$m$--dimensional manifold $M$ arises from a non--degenerate map
$\f:U\to \R^n$ where $U$ is an open subset of $\R^m$ and
$M:=\f(U)$. The  map $\f:U\to \R^n:\bu\mapsto
\f(\bu)=(f_1(\bu),\dots,f_n(\bu))$ is said to be
\emph{non--degenerate} at $\bu\in U$ if there exists some $l\in\N$
such that $\f$  is $l$ times continuously differentiable on some
sufficiently small ball centred at $\bu$ and the partial
derivatives of $\f$ at $\bu$ of orders up to $l$  span $\R^n$. The
map $\f$ is \emph{non--degenerate} if it is non--degenerate at
almost every (in terms of $m$--dimensional Lebesgue measure) point
in $U$; in turn the manifold $M=\f(U)$ is also said to be
non--degenerate. Geometrically, the non--degeneracy of $M$ at
$\by_0\in M$ means that $M$ deviates from any hyperplane in $\R^n$
that contains $\by_0$.

 Consider any ball $B_{\uo}$ centred at $\uo$ in $U$. Then, with
reference to our general framework let  $\Omega:= B_{\uo}   , J:=
\{ ( p,{\bf q}) \in \Z \times \Z^n \backslash \{ {\bf 0}\}: |p|
\leq |\bq| \} \ , \ \alpha := (p,\bq) \in J \ , \ \ma := |\bq| \
{\rm and \ }
  \ra:=  \{ \bu \in B_{\uo}:  \bq.f(\bu) = p \}. $
Thus, the  family $\mathcal{ R}$ of resonant sets $\ra$ have dimension
$m-1$ and arise from the intersection of $M$ with the
$(n-1)$--dimensional hyperplanes given by $\{\bx\in\R^n:\bq.\bx=p$
with $(p,\bq)\in J\}$.
  Furthermore,
$\De(\ra,\p(\ma)) := \{ \bu \in B_{\uo}:  {\rm dist} (\bu, \ra)
\leq \psi\left(|\bq|\right) \} \  $ and $$\De^u_l(\p, n) :=
\bigcup_{k^{n-1} < |\bq| \leq k^n} \ \bigcup_{|\bp| \leq |\bq| }
\De(\ra,\p(\ma)) \ . $$ Here $k > 1$ is a constant.  Then
$$W(M;\psi)\supset \Lambda(\p) := \limsup_{n \to \infty}
\De^u_l(\p, n) \ \ . $$ Note that $$\Lambda(\p)=\{\bu\in
B_{\uo}:|\bq.f(\bu)-p|<|\bq|\p(|\bq|)\ \mathrm{for\ i.m.\
}(p,\bq)\in\Z\times\Z^n\}. $$ Now let the measure  $m$  be
normalised $m$--dimensional Lebesgue measure on $B_{\uo}$, $
\delta = m$ and  $\gamma= m-1 $. Then for almost all $\uo\in U$
there exists a corresponding
 ball $B_{\uo}$ for which  the following local $m$--ubiquity statement
holds.

\medskip

\begin{proposition}\label{mani}
The pair $(\mathcal{ R}, \beta) $ is a local $m$--ubiquitous system
relative to $(\r,l,u)$ where for $k > k_0(B_{\uo})$ -- a positive
absolute constant, $l_{t+1} = u_t := k^t$ $(t \in \N)$ and $\r : r
\to \r(r)  := \, {\rm constant}
 \times  \r^{-(n+1)}$.
\end{proposition}

This proposition can be deduced from  Proposition 3.3 in
\cite{bbkm}. In view of the proposition and the fact that the
measure $m$ is of type (M2) and that $\rho$ is $u$-regular,
Corollary \ref{cor2} of Theorem \ref{THM1}  and Theorem \ref{THM3}
yield the divergent parts of the following statements.

\medskip

\begin{theorem}\label{manithm1}\label{manikg}
Let $\p$ be a real, positive  decreasing function, $M$ be a
non-degenerate manifold, $m$ be the induced Lebesgue measure on
$M$. Then $$
 m(W(M;\psi) ) = \left\{
\begin{array}{cl}
 0 \ \ & {\rm if} \;\;\;
 \sum_{r=1}^\infty \;  \psi\left(r\right)  \;\; r^{n} <\infty\;
 ,\\[1.5ex]
  m(M) \ \  & {\rm if} \;\;\; \sum_{r=1}^\infty \; \psi\left(r\right) \;\;
r^{n} =\infty \; .
\end{array}\right.
$$
\end{theorem}

\medskip

\begin{theorem}\label{manithm2}
Let $f$ be  a dimension function such that $r^{-m}f(r)\to\infty $
as $r\to0$ and $r^{-m}f(r)$ is decreasing. Furthermore, suppose
that $ r^{-(m-1)} f(r) $ is  increasing. Let $\p$ be a real,
positive decreasing function. Then $$ \mathcal{H}^f(W(M;\psi))=
\infty\   \ \ \ \  \mathrm{if} \ \ \ \ \sum_{r=1}^\infty \fp r))
\, \p(r)^{-(m-1) } \, r^n =\infty. $$
\end{theorem}

 Theorem \ref{manithm1} is not new. The convergence
half was independently proved in \cite{vb2} and \cite{bkm} and the
divergence half in \cite{bbkm}. The convergence part  is not at
all obvious and requires delicate covering and counting arguments
to make use of the  `natural' cover of $W(M;\psi)$.
 In fact it implies that any non--degenerate manifold is extremal which
was a longstanding conjecture of Baker--Sprindzuk. This conjecture
was proved independently in \cite{KM}.

 Theorem \ref{manithm2} is  new. It shows that any non-degenerate
 manifold is of \Ja type for divergence in the case of dual
 approximation where \Ja type  is the Hausdorff measure analogue of the
 notion of Khintchine/Groshev type \cite[pg 29]{BD}.
However, unlike in previous applications, we are currently
 unable to show that
$\mathcal{H}^f(W(M;\psi))=0$ when the sum in Theorem
\ref{manithm2} converges.
 Previously \cite{DD}, it had been
shown that $\dim W(M;\tau)\ge m-1+(n+1)/\tau$ for $\tau> n+1$.
Currently there is no  general upper bound. As usual $W(M;\tau)$
corresponds to the set $W(M;\psi)$ with $\psi: r \mapsto
r^{-\tau}$. Obviously Theorem \ref{manithm2} implies this
dimension result and shows that $\mathcal{H}^s(W(M;\tau))=\infty$
for $s= m-1+(n+1)/\tau$ which is almost certainly the critical
exponent.


\subsubsection{Simultaneous approximation on
manifolds.}\label{circleapp}
For an approximating function $\p$, let $$ W(M;\psi)=\{\bx\in M:
|q\bx-\bp|<|q|\psi(|q|)\ \mathrm{for\ i.m.\
}(\bp,q)\in\Z^n\times\Z\}. $$ Recall that $M$ is an
$m$--dimensional submanifold embedded in $\R^n$. Even under the
restriction that $M$ is non--degenerate, results analogous to
those described above  for dual approximation currently seem out
of reach. However in \cite{us3}, we have recently made  advances
in the case that $M$ is a  $C^3$ planar curve. Moreover, for
particular planar curves such as the unit circle and the parabola
we are able to establish reasonably complete measure theoretic
laws. To this end, let $M= S^1:= \{(x_1,x_2) \in \R^2 : x_1^2 +
x_2^2 =1 \} $ denote the unit circle. Thus, $W(S^1;\p)$ consists
of points $(x_1,x_2) \in S^1$ for which there exist i.m. rational
pairs $(\frac{p_1}{q}, \frac{p_2}{q})$ such that the following
pair of inequalities are simultaneously satisfied:
\begin{equation}
|x_i - p_i/q|  \  < \  \psi(|q|) \hspace{1cm} (i=1,2) \ .
\label{cirapp}
\end{equation}

\noindent {\em Throughout the following discussion, assume  that
$r^2\p(r) \to 0 $ as $r \to \infty$.} Also, without loss of
generality assume that $q \in \N$. The following fact shows that
under the above assumption on $\psi$,  there is a one to one
correspondence between the rational approximates $(p_1/q, \,
p_2/q)$ satisfying (\ref{cirapp}) and the Pythagorean triples $s^2
+t^2 = q^2$.

\noindent {\bf $\bullet$} For $q$ large, any rational pair
$(p_1/q, \, p_2/q)$ satisfying (\ref{cirapp}) lies on $S^1$.

\smallskip

This is trivial, since $(1- \p(q))^2 \leq (p_1/q)^2 + (p_2/q)^2 <
(1+3\p(q))^2$ and so for $q$ sufficiently large  $ |p_1^2 + p_2^2
- q^2| < 1$. Now  notice that the left hand side is an integer.

\noindent Now with reference to our general framework, let
$\Omega:= S^1, J:= \{ (\bp,q) \in \Z^2 \times \N : \bp/q \in S^1
\} , \ \alpha := (\bp,q) \in J \ , \ \ma := q \ {\rm and \ } \ra:=
\bp/q. $ Thus, the family $\mathcal{ R}$ of resonant sets $\ra$
consists of rational pairs $(p_1/q, \, p_2/q)$  lying on $S^1$.
Furthermore, define $\De(\ra,\p(\ma)) := B(\bp/q ,\p( q )   ) $
and let
$$\De^u_l(\p, n) :=   \bigcup_{2^{n-1} < \,  q \leq
 2^n  }
  \ \bigcup_{\substack{\bp \in \Z^2 : \\ p_1^2 + p_2^2 = q^2 }   } B(\bp/q, \p(q)) \ ,
$$ Then $$W(S^1;\p)  \supset \Lambda(\p) :=  \limsup_{n \to
\infty} \De^u_l(\p, n) \,  . $$ Now let $m$ be normalized
 Lebesgue measure on $S^1$, $\delta = 1$ and
$\gamma=0 $. Then we have  the following local $m$--ubiquity
statement.

\medskip

\begin{proposition}\label{circle}
The pair $(\mathcal{ R}, \beta) $ is a local $m$--ubiquitous system
relative to $(\r,l,u)$, where $l_{t+1} = u_t := 2^t$   $(t \in
\N)$  and $\r : r \to \r(r)  :=  \, {\rm constant}
 \times r^{-1}.$
\end{proposition}
Although this statement can be found explicitly in \cite{Dick}, we
shall give an alternative proof which is  shorter and probably
more adaptable to analogous problems. For a point  $a$  on $S^1$,
let $A$ be the arc with centre (mid-point) $a$ and radius $r(A)$.
Clearly $r(A) \asymp m(A)$. For $N \in \N$, let $ \mathcal{ C}(N,A) :=
\{ (\bp,q) \in J : \bp/q \in A {\rm \ with \  } N <   q \leq
 2N   \}  $.  The proposition is a simple consequence of
 the following two facts.

\noindent {\bf $\bullet$} For $N \geq N_o(A)$, $ \#\mathcal{ C}(N,A)
\asymp r(A) \, N \asymp m(A) \, N $.

This follows from standard results on the distribution of
Pytha\-go\-rean triples.

 \noindent {\bf $\bullet$} For $N \geq N_o$, if
 $r(A) < 2^{-4/3} N^{-1}$ then $ \#\mathcal{ C}(N,A) \leq 2 $.

This is the key result. Let $b$ and $b'$ be the end points of $A$.
Then for $N$ large enough, the triangle $T$ subtended by the three
points $ a,b,b'$ has area  less than $ r(A)^3 $. Now suppose there
are three rational points $ {\bf s}/q, {\bf t}/q', {\bf u}/q'' \in
\mathcal{ C}(N,A)$ and let $\triangle$ be the triangle subtended
by them. Clearly, ${\rm area}(\triangle) \leq {\rm area}(T)$. Thus
$$ 2 \, r(A)^3  \ \geq \ 2 \, {\rm area}(\triangle) \ = \  \left|
 \begin{array}{ccc} 1 & s/q & t/q \\   1 & s'/q' & t'/q'
\\ 1 & s''/q'' & t''/q''
\end{array}  \right|  \ \geq \  \frac{1}{q\, q'q''}   \ \geq \
 \frac{1}{(2N)^3}\ . $$
Hence if $r(A) < 2^{-4/3} N^{-1}$, the triangle $\triangle$ cannot
exist so the three rational points 
must lie on a straight line. However this is impossible
since they lie on $S^1$.

\medskip

 Thus, in view of the proposition and the fact that
the measure $m$ is of type (M2) and that $\r$ is $u$-regular,
Theorem \ref{THM3} yields the divergent parts of the following
statement.

\medskip

\begin{theorem}\label{cirthm}
 Let $f$ be a dimension function such that
$r^{-1} \, f(r)\to \infty$ as $r\to 0 \, $ and  $r^{-1} \, f(r) $
is decreasing. Let $\psi$  be a real,  positive decreasing
function such that $r^2\p(r) \to 0 $ as $r \to \infty$. Then $$
\hf\left(W(S^1;\p)\right)= \left\{
\begin{array}{ll}
 0 \ \  & {\rm if}
\;\;\; \sum_{r=1}^\infty \;  f\left(\p(r)\right)  <\infty \;
,\\[1.5ex]
 \infty \ \ & {\rm if} \;\;\; \sum_{r=1}^\infty \;
f\left(\p(r)\right)  =\infty \; .
\end{array}\right.
$$
\end{theorem}

Let us consider the case when $\psi(r) = r^{-\tau}$ and write
$W(S^1;\tau)$ for $W(S^1;\p)$. In \cite{Dick}, it was shown that
 $\dim W(S^1\tau) = 1/\tau $ when $\tau > 2 $. Clearly, Theorem
\ref{cirthm} implies this dimension result and shows that $\mathcal{
H} ^{1/\tau}(W(S^1;\tau)) = \infty$. Note that the condition that
$\tau > 2 $, ensures that $r^2\p(r) \to 0 $ as $r \to \infty$ and
so that the  rational points of interest are forced to lie on the
circle. This is not the case when $\tau \leq 2$ and even the
problem of determining  the dimension of $W(S^1;\tau)$ becomes
highly non-trivial. In view of Dirchlet's theorem on simultaneous
Diophantine approximation one knows that $W(S^1;\tau)= S^1$ for
$\tau \leq 3/2$. On the other hand, it is not difficult to show
that $m(W(S^1;\tau)) =0 $ for $\tau > 3/2$ and moreover that $\dim
W(S^1;\tau) \leq  (3-\tau) /\tau  $ when $3/2 \leq \tau \leq 2 $.

In a forthcoming paper \cite{us3}, general measure theoretic laws
for  non-degenerate $C^3$ planar curves are established. A simple
consequence of these results is that:  $$ \dim W(S^1;\tau) =
(3-\tau) /\tau  \  \hspace{14mm} (3/2 \leq \tau \leq 2)  \ .  $$

\noindent Another consequence of the results in  \cite{us3} is the
following law with respect to the measure $m$ on $S^1$.

\medskip

\begin{thcs}  Let $\psi $ be  a
real, positive decreasing function. Then $$ m
\left(W(S^1;\p)\right)= \left\{
\begin{array}{cl}
 0 \ \  & {\rm if} \;\;\; \sum_{r=1}^\infty \; \left(r\, \psi(r)\right)^2 <
\infty \; ,\\[1.5ex]
 1 \ \ & {\rm if} \;\;\; \sum_{r=1}^\infty \;
\left(r\, \psi(r)\right)^2 = \infty \; .
\end{array}
\right. $$
\end{thcs}

\subsection{Sets of exact order}
With reference to our general framework, given two approximating
functions $\varphi$ and $\psi$ with $\varphi$ in some sense
`smaller' than $\psi$, consider the set $E(\psi,\varphi)$ of
points $x$ in $\Omega$ for which
 $$x\in\De(\ra,\p(\ma)) \hspace{4mm} {\rm \ for \ i.m. \ } \alpha \in
J  ; $$ and that $$
 \ \ \ \ \ \ \ \ x \notin \De(\ra,\varphi(\ma))  \hspace{5mm} { \rm \ for \ all \ but \
 finitely \ many \,
} \alpha \in J   \ . $$ In short, $E(\psi,\varphi) :=
\Lambda(\psi) \setminus \Lambda(\varphi)$.  Thus the approximation
properties of points $x$ in $E(\psi,\varphi)$ are `sandwiched'
between the  functions $\varphi$ and $\psi$.  In \cite{exactus3},
under the classical linear forms setup (cf. \S\ref{lfapp}), we
have shown that the measure theoretic laws for $W(m,n;\psi)$ with
respect to the measures $m$ (Theorem (Khintchine--Groshev)) and
$\mathcal{ H}^f$ (Theorem DV) give rise to precise metric results for
the corresponding `exact order' sets $E(m,n;\psi,\varphi)$. In
short, the key idea is to construct an appropriate dimension
function $f$ for which  $\mathcal{ H}^f(W(m,n;\psi))=\infty$ and
$\mathcal{ H}^f(W(m,n;\varphi))=0$ and so $\mathcal{
H}^f(E(m,n;\psi,\varphi)) = \infty$.

Regarding our general framework and the exact order sets
$E(\psi,\varphi)$, the arguments used in \cite{exactus3} can  be
carried over to obtain analogous statements of the theorems in
\cite{exactus3} provided analogues of both the
Khin\-tchine--Groshev Theorem and Theorem DV hold. The point is
that both the divergence and convergence halves of these theorems
are required.  In any case, this is the case  for all the
applications considered in this section except the previous one
and the Kleinian groups application when the group is of the
second kind with parabolic elements. Thus, for each of the other
applications it is possible to  describe the measure theoretic
structure of the associated exact order sets $E(\psi,\varphi)$. We
shall leave the details to the energetic reader.



\vspace{6mm}

{\bf Acknowledgements.} SV would like to thank  his new friends,
Ayesha and Iona for their permanent smiles and  wonderfully
positive `outlook' -- long may it last !!  Finally and most
importantly of all  he would like to thank Bridget Bennett for ...
just about everything. 

 We would like to thank the reviewer for her/his comments
which have made the early part of the paper more accessible and
for the reference to Knopp \cite{knopp}. We would also like to
thank the editor for his professionalism and `speed of response'
during the various stages of the reviewing  process.

\bibliographystyle{amsalpha}

\end{document}